%param.2
\magnification=\magstep1
\def\firstpage{1}
\pageno=\firstpage
%fonts.6
\font\fiverm=cmr5
\font\sevenrm=cmr7
\font\sevenbf=cmbx7
\font\eightrm=cmr8
\font\eightbf=cmbx8
\font\ninerm=cmr9
\font\ninebf=cmbx9
\font\tenbf=cmbx10
\font\magtenbf=cmbx10 scaled\magstep1

\font\magnineeufm=eufm9 scaled\magstep1

%
%smallfonts.tex
%
\newskip\ttglue
\font\fiverm=cmr5
\font\fivei=cmmi5
\font\fivesy=cmsy5
\font\fivebf=cmbx5
\font\sixrm=cmr6
\font\sixi=cmmi6
\font\sixsy=cmsy6
\font\sixbf=cmbx6
\font\sevenrm=cmr7
\font\eightrm=cmr8
\font\eighti=cmmi8
\font\eightsy=cmsy8
\font\eightit=cmti8
\font\eightsl=cmsl8
\font\eighttt=cmtt8
\font\eightbf=cmbx8
\font\ninerm=cmr9
\font\ninei=cmmi9
\font\ninesy=cmsy9
\font\nineit=cmti9
\font\ninesl=cmsl9
\font\ninett=cmtt9
\font\ninebf=cmbx9
\font\twelverm=cmr12
\font\twelvei=cmmi12
\font\twelvesy=cmsy12
\font\twelveit=cmti12
\font\twelvesl=cmsl12
\font\twelvett=cmtt12
\font\twelvebf=cmbx12

%% EIGHT POINT FONT FAMILY

\def\eightpoint{\def\rm{\fam0\eightrm}  
  \textfont0=\eightrm \scriptfont0=\sixrm \scriptscriptfont0=\fiverm
  \textfont1=\eighti  \scriptfont1=\sixi  \scriptscriptfont1=\fivei
  \textfont2=\eightsy  \scriptfont2=\sixsy  \scriptscriptfont2=\fivesy
  \textfont3=\tenex  \scriptfont3=\tenex  \scriptscriptfont3=\tenex
  \textfont\itfam=\eightit  \def\it{\fam\itfam\eightit}
  \textfont\slfam=\eightsl  \def\sl{\fam\slfam\eightsl}
  \textfont\ttfam=\eighttt  \def\tt{\fam\ttfam\eighttt}
  \textfont\bffam=\eightbf  \scriptfont\bffam=\sixbf
    \scriptscriptfont\bffam=\fivebf  \def\bf{\fam\bffam\eightbf}
  \tt  \ttglue=.5em plus.25em minus.15em
  \normalbaselineskip=9pt
  \setbox\strutbox=\hbox{\vrule height7pt depth2pt width0pt}
  \let\sc=\sixrm  \let\big=\eightbig \normalbaselines\rm}

\def\eightbig#1{{\hbox{$\textfont0=\ninerm\textfont2=\ninesy
        \left#1\vbox to6.5pt{}\right.$}}}

%% NINE POINT FONT FAMILY

\def\ninepoint{\def\rm{\fam0\ninerm}  
  \textfont0=\ninerm \scriptfont0=\sixrm \scriptscriptfont0=\fiverm
  \textfont1=\ninei  \scriptfont1=\sixi  \scriptscriptfont1=\fivei
  \textfont2=\ninesy  \scriptfont2=\sixsy  \scriptscriptfont2=\fivesy
  \textfont3=\tenex  \scriptfont3=\tenex  \scriptscriptfont3=\tenex
  \textfont\itfam=\nineit  \def\it{\fam\itfam\nineit}
  \textfont\slfam=\ninesl  \def\sl{\fam\slfam\ninesl}
  \textfont\ttfam=\ninett  \def\tt{\fam\ttfam\ninett}
  \textfont\bffam=\ninebf  \scriptfont\bffam=\sixbf
    \scriptscriptfont\bffam=\fivebf  \def\bf{\fam\bffam\ninebf}
  \tt  \ttglue=.5em plus.25em minus.15em
  \normalbaselineskip=11pt
  \setbox\strutbox=\hbox{\vrule height8pt depth3pt width0pt}
  \let\sc=\sevenrm  \let\big=\ninebig \normalbaselines\rm}

\def\ninebig#1{{\hbox{$\textfont0=\tenrm\textfont2=\tensy
        \left#1\vbox to7.25pt{}\right.$}}}

%% TWELVE POINT FONT FAMILY --- not really small

\def\twelvepoint{\def\rm{\fam0\twelverm}  
  \textfont0=\twelverm \scriptfont0=\eightrm \scriptscriptfont0=\sixrm
  \textfont1=\twelvei  \scriptfont1=\eighti  \scriptscriptfont1=\sixi
  \textfont2=\twelvesy  \scriptfont2=\eightsy  \scriptscriptfont2=\sixsy
  \textfont3=\tenex  \scriptfont3=\tenex  \scriptscriptfont3=\tenex
  \textfont\itfam=\twelveit  \def\it{\fam\itfam\twelveit}
  \textfont\slfam=\twelvesl  \def\sl{\fam\slfam\twelvesl}
  \textfont\ttfam=\twelvett  \def\tt{\fam\ttfam\twelvett}
  \textfont\bffam=\twelvebf  \scriptfont\bffam=\eightbf
    \scriptscriptfont\bffam=\sixbf  \def\bf{\fam\bffam\twelvebf}
  \tt  \ttglue=.5em plus.25em minus.15em
  \normalbaselineskip=11pt
  \setbox\strutbox=\hbox{\vrule height8pt depth3pt width0pt}
  \let\sc=\sevenrm  \let\big=\twelvebig \normalbaselines\rm}

\def\twelvebig#1{{\hbox{$\textfont0=\tenrm\textfont2=\tensy
        \left#1\vbox to7.25pt{}\right.$}}}
\catcode`\@=11
%
%  Include all definitions related to the fonts msam, msbm and eufm, so that
%  when this file is used by itself, the results with respect to those fonts
%  are equivalent to what they would have been using AMS-TeX.
%  Most symbols in fonts msam and msbm are defined using \newsymbol;
%  however, a few symbols that replace composites defined in plain must be
%  defined with \mathchardef.

\def\undefine#1{\let#1\undefined}
\def\newsymbol#1#2#3#4#5{\let\next@\relax
 \ifnum#2=\@ne\let\next@\msafam@\else
 \ifnum#2=\tw@\let\next@\msbfam@\fi\fi
 \mathchardef#1="#3\next@#4#5}
\def\mathhexbox@#1#2#3{\relax
 \ifmmode\mathpalette{}{\m@th\mathchar"#1#2#3}%
 \else\leavevmode\hbox{$\m@th\mathchar"#1#2#3$}\fi}
\def\hexnumber@#1{\ifcase#1 0\or 1\or 2\or 3\or 4\or 5\or 6\or 7\or 8\or
 9\or A\or B\or C\or D\or E\or F\fi}

\font\tenmsa=msam10
\font\sevenmsa=msam7
\font\fivemsa=msam5
\newfam\msafam
\textfont\msafam=\tenmsa
\scriptfont\msafam=\sevenmsa
\scriptscriptfont\msafam=\fivemsa
\edef\msafam@{\hexnumber@\msafam}
\mathchardef\dabar@"0\msafam@39
\def\dashrightarrow{\mathrel{\dabar@\dabar@\mathchar"0\msafam@4B}}
\def\dashleftarrow{\mathrel{\mathchar"0\msafam@4C\dabar@\dabar@}}

\def\ulcorner{\delimiter"4\msafam@70\msafam@70 }
\def\urcorner{\delimiter"5\msafam@71\msafam@71 }
\def\llcorner{\delimiter"4\msafam@78\msafam@78 }
\def\lrcorner{\delimiter"5\msafam@79\msafam@79 }
%    Note that there should not be a final space after the digits for a
%    \mathhexbox@.
\def\yen{{\mathhexbox@\msafam@55}}
\def\checkmark{{\mathhexbox@\msafam@58}}
\def\circledR{{\mathhexbox@\msafam@72}}
\def\maltese{{\mathhexbox@\msafam@7A}}

\font\tenmsb=msbm10
\font\sevenmsb=msbm7
\font\fivemsb=msbm5
\newfam\msbfam
\textfont\msbfam=\tenmsb
\scriptfont\msbfam=\sevenmsb
\scriptscriptfont\msbfam=\fivemsb
\edef\msbfam@{\hexnumber@\msbfam}
\def\Bbb#1{{\fam\msbfam\relax#1}}
\def\widehat#1{\setbox\z@\hbox{$\m@th#1$}%
 \ifdim\wd\z@>\tw@ em\mathaccent"0\msbfam@5B{#1}%
 \else\mathaccent"0362{#1}\fi}
\def\widetilde#1{\setbox\z@\hbox{$\m@th#1$}%
 \ifdim\wd\z@>\tw@ em\mathaccent"0\msbfam@5D{#1}%
 \else\mathaccent"0365{#1}\fi}
\font\teneufm=eufm10
\font\seveneufm=eufm7
\font\fiveeufm=eufm5
\newfam\eufmfam
\textfont\eufmfam=\teneufm
\scriptfont\eufmfam=\seveneufm
\scriptscriptfont\eufmfam=\fiveeufm

\catcode`\@=11
%%  Load amssym.def if necessary: If \newsymbol is undefined, do nothing
%%  and the following \input statement will be executed; otherwise
%%  change \input to a temporary no-op.
%#\ifx\undefined\newsymbol \else \begingroup\def\input#1 {\endgroup}\fi
%#\input amssym.def \relax
%%  Most symbols in fonts msam and msbm are defined using \newsymbol.  A few
%%  that are delimiters or otherwise require special treatment have already
%%  been defined as soon as the fonts were loaded.  Finally, a few symbols
%%  that replace composites defined in plain must be undefined first.
\newsymbol\boxdot 1200
\newsymbol\boxplus 1201
\newsymbol\boxtimes 1202
\newsymbol\square 1003
\newsymbol\blacksquare 1004
\newsymbol\centerdot 1205
\newsymbol\lozenge 1006
\newsymbol\blacklozenge 1007
\newsymbol\circlearrowright 1308
\newsymbol\circlearrowleft 1309
\undefine\rightleftharpoons
\newsymbol\rightleftharpoons 130A
\newsymbol\leftrightharpoons 130B
\newsymbol\boxminus 120C
\newsymbol\Vdash 130D
\newsymbol\Vvdash 130E
\newsymbol\vDash 130F
\newsymbol\twoheadrightarrow 1310
\newsymbol\twoheadleftarrow 1311
\newsymbol\leftleftarrows 1312
\newsymbol\rightrightarrows 1313
\newsymbol\upuparrows 1314
\newsymbol\downdownarrows 1315
\newsymbol\upharpoonright 1316
 
\newsymbol\downharpoonright 1317
\newsymbol\upharpoonleft 1318
\newsymbol\downharpoonleft 1319
\newsymbol\rightarrowtail 131A
\newsymbol\leftarrowtail 131B
\newsymbol\leftrightarrows 131C
\newsymbol\rightleftarrows 131D
\newsymbol\Lsh 131E
\newsymbol\Rsh 131F
\newsymbol\rightsquigarrow 1320
\newsymbol\leftrightsquigarrow 1321
\newsymbol\looparrowleft 1322
\newsymbol\looparrowright 1323
\newsymbol\circeq 1324
\newsymbol\succsim 1325
\newsymbol\gtrsim 1326
\newsymbol\gtrapprox 1327
\newsymbol\multimap 1328
\newsymbol\therefore 1329
\newsymbol\because 132A
\newsymbol\doteqdot 132B
 
\newsymbol\triangleq 132C
\newsymbol\precsim 132D
\newsymbol\lesssim 132E
\newsymbol\lessapprox 132F
\newsymbol\eqslantless 1330
\newsymbol\eqslantgtr 1331
\newsymbol\curlyeqprec 1332
\newsymbol\curlyeqsucc 1333
\newsymbol\preccurlyeq 1334
\newsymbol\leqq 1335
\newsymbol\leqslant 1336
\newsymbol\lessgtr 1337
\newsymbol\backprime 1038
\newsymbol\risingdotseq 133A
\newsymbol\fallingdotseq 133B
\newsymbol\succcurlyeq 133C
\newsymbol\geqq 133D
\newsymbol\geqslant 133E
\newsymbol\gtrless 133F
\newsymbol\sqsubset 1340
\newsymbol\sqsupset 1341
\newsymbol\vartriangleright 1342
\newsymbol\vartriangleleft 1343
\newsymbol\trianglerighteq 1344
\newsymbol\trianglelefteq 1345
\newsymbol\bigstar 1046
\newsymbol\between 1347
\newsymbol\blacktriangledown 1048
\newsymbol\blacktriangleright 1349
\newsymbol\blacktriangleleft 134A
\newsymbol\vartriangle 134D
\newsymbol\blacktriangle 104E
\newsymbol\triangledown 104F
\newsymbol\eqcirc 1350
\newsymbol\lesseqgtr 1351
\newsymbol\gtreqless 1352
\newsymbol\lesseqqgtr 1353
\newsymbol\gtreqqless 1354
\newsymbol\Rrightarrow 1356
\newsymbol\Lleftarrow 1357
\newsymbol\veebar 1259
\newsymbol\barwedge 125A
\newsymbol\doublebarwedge 125B
\undefine\angle
\newsymbol\angle 105C
\newsymbol\measuredangle 105D
\newsymbol\sphericalangle 105E
\newsymbol\varpropto 135F
\newsymbol\smallsmile 1360
\newsymbol\smallfrown 1361
\newsymbol\Subset 1362
\newsymbol\Supset 1363
\newsymbol\Cup 1264
 
\newsymbol\Cap 1265
 
\newsymbol\curlywedge 1266
\newsymbol\curlyvee 1267
\newsymbol\leftthreetimes 1268
\newsymbol\rightthreetimes 1269
\newsymbol\subseteqq 136A
\newsymbol\supseteqq 136B
\newsymbol\bumpeq 136C
\newsymbol\Bumpeq 136D
\newsymbol\lll 136E
 
\newsymbol\ggg 136F
 
\newsymbol\circledS 1073
\newsymbol\pitchfork 1374
\newsymbol\dotplus 1275
\newsymbol\backsim 1376
\newsymbol\backsimeq 1377
\newsymbol\complement 107B
\newsymbol\intercal 127C
\newsymbol\circledcirc 127D
\newsymbol\circledast 127E
\newsymbol\circleddash 127F
\newsymbol\lvertneqq 2300
\newsymbol\gvertneqq 2301
\newsymbol\nleq 2302
\newsymbol\ngeq 2303
\newsymbol\nless 2304
\newsymbol\ngtr 2305
\newsymbol\nprec 2306
\newsymbol\nsucc 2307
\newsymbol\lneqq 2308
\newsymbol\gneqq 2309
\newsymbol\nleqslant 230A
\newsymbol\ngeqslant 230B
\newsymbol\lneq 230C
\newsymbol\gneq 230D
\newsymbol\npreceq 230E
\newsymbol\nsucceq 230F
\newsymbol\precnsim 2310
\newsymbol\succnsim 2311
\newsymbol\lnsim 2312
\newsymbol\gnsim 2313
\newsymbol\nleqq 2314
\newsymbol\ngeqq 2315
\newsymbol\precneqq 2316
\newsymbol\succneqq 2317
\newsymbol\precnapprox 2318
\newsymbol\succnapprox 2319
\newsymbol\lnapprox 231A
\newsymbol\gnapprox 231B
\newsymbol\nsim 231C
\newsymbol\ncong 231D
\newsymbol\diagup 201E
\newsymbol\diagdown 201F
\newsymbol\varsubsetneq 2320
\newsymbol\varsupsetneq 2321
\newsymbol\nsubseteqq 2322
\newsymbol\nsupseteqq 2323
\newsymbol\subsetneqq 2324
\newsymbol\supsetneqq 2325
\newsymbol\varsubsetneqq 2326
\newsymbol\varsupsetneqq 2327
\newsymbol\subsetneq 2328
\newsymbol\supsetneq 2329
\newsymbol\nsubseteq 232A
\newsymbol\nsupseteq 232B
\newsymbol\nparallel 232C
\newsymbol\nmid 232D
\newsymbol\nshortmid 232E
\newsymbol\nshortparallel 232F
\newsymbol\nvdash 2330
\newsymbol\nVdash 2331
\newsymbol\nvDash 2332
\newsymbol\nVDash 2333
\newsymbol\ntrianglerighteq 2334
\newsymbol\ntrianglelefteq 2335
\newsymbol\ntriangleleft 2336
\newsymbol\ntriangleright 2337
\newsymbol\nleftarrow 2338
\newsymbol\nrightarrow 2339
\newsymbol\nLeftarrow 233A
\newsymbol\nRightarrow 233B
\newsymbol\nLeftrightarrow 233C
\newsymbol\nleftrightarrow 233D
\newsymbol\divideontimes 223E
\newsymbol\varnothing 203F
\newsymbol\nexists 2040
\newsymbol\Finv 2060
\newsymbol\Game 2061
\newsymbol\mho 2066
\newsymbol\eth 2067
\newsymbol\eqsim 2368
\newsymbol\beth 2069
\newsymbol\gimel 206A
\newsymbol\daleth 206B
\newsymbol\lessdot 236C
\newsymbol\gtrdot 236D
\newsymbol\ltimes 226E
\newsymbol\rtimes 226F
\newsymbol\shortmid 2370
\newsymbol\shortparallel 2371
\newsymbol\smallsetminus 2272
\newsymbol\thicksim 2373
\newsymbol\thickapprox 2374
\newsymbol\approxeq 2375
\newsymbol\succapprox 2376
\newsymbol\precapprox 2377
\newsymbol\curvearrowleft 2378
\newsymbol\curvearrowright 2379
\newsymbol\digamma 207A
\newsymbol\varkappa 207B
\newsymbol\Bbbk 207C
\newsymbol\hslash 207D
\undefine\hbar
\newsymbol\hbar 207E
\newsymbol\backepsilon 237F
%  Restore the catcode value for @ that was previously saved.
%#\catcode`\@=\csname pre amssym.tex at\endcsname

%\endinput
% links.1
% adapted from http://insti.physics.sunysb.edu/~siegel/tex.shtml
%
% postscript/pdf
\newcount\marknumber	\marknumber=1
\newcount\countdp \newcount\countwd \newcount\countht 
%
% for ordinary tex
%
\ifx\pdfoutput\undefined
\def\rgboo#1{}
\def\postscript#1{\special{" #1}}		%% for dvips
\postscript{
	/bd {bind def} bind def
	/fsd {findfont exch scalefont def} bd
	/sms {setfont moveto show} bd
	/ms {moveto show} bd
	/pdfmark where		% printers ignore pdfmarks
	{pop} {userdict /pdfmark /cleartomark load put} ifelse
	[ /PageMode /UseOutlines		% bookmark window open
	/DOCVIEW pdfmark}
\def\bookmark#1#2{\postscript{		% #1=subheadings (if not 0)
	[ /Dest /MyDest\the\marknumber /View [ /XYZ null null null ] /DEST pdfmark
	[ /Title (#2) /Count #1 /Dest /MyDest\the\marknumber /OUT pdfmark}%
	\advance\marknumber by1}
\def\pdfclink#1#2#3{%
	\hskip-.25em\setbox0=\hbox{#2}%
		\countdp=\dp0 \countwd=\wd0 \countht=\ht0%
		\divide\countdp by65536 \divide\countwd by65536%
			\divide\countht by65536%
		\advance\countdp by1 \advance\countwd by1%
			\advance\countht by1%
		\def\linkdp{\the\countdp} \def\linkwd{\the\countwd}%
			\def\linkht{\the\countht}%
	\postscript{
		[ /Rect [ -1.5 -\linkdp.0 0\linkwd.0 0\linkht.5 ] 
		/Border [ 0 0 0 ]
		/Action << /Subtype /URI /URI (#3) >>
		/Subtype /Link
		/ANN pdfmark}{\rgb{#1}{#2}}}
%
% for pdftex
%
\else
\def\rgboo#1{\pdfliteral{#1 rg #1 RG}}
\pdfcatalog{/PageMode /UseOutlines}		% bookmark window open
\def\bookmark#1#2{
	\pdfdest num \marknumber xyz
	\pdfoutline goto num \marknumber count #1 {#2}
	\advance\marknumber by1}
\def\pdfklink#1#2{%
	\noindent\pdfstartlink user
		{/Subtype /Link
		/Border [ 0 0 0 ]
		/A << /S /URI /URI (#2) >>}{\rgb{1 0 0}{#1}}%
	\pdfendlink}
\fi

\def\rgbo#1#2{\rgboo{#1}#2\rgboo{0 0 0}}
\def\rgb#1#2{\mark{#1}\rgbo{#1}{#2}\mark{0 0 0}}
\def\pdfklink#1#2{\pdfclink{1 0 0}{#1}{#2}}
\def\pdflink#1{\pdfklink{#1}{#1}}
%
% examples:
% \bookmark{0}{look here}
% \pdfclink{0 0 1}{testlink}{http://www.google.com/}
% \pdfklink{testlink}{http://www.google.com/}
% \pdflink{http://www.google.com/}
%titles.8
% requires fonts.5 or higher and smallfonts.tex
% uses links.* if included
% enumerates \demo consecutively (no section number)
%
\newcount\seccount  %% sections
\newcount\subcount  %% subsection
\newcount\clmcount  %% claim
\newcount\equcount  %% equation
\newcount\refcount  %% reference
\newcount\demcount  %% example
\newcount\execount  %% exercise
\newcount\procount  %% problem
\seccount=0
\equcount=1
\clmcount=1
\subcount=1
\refcount=1
\demcount=0
\execount=0
\procount=0
%
%% MISC STUFF
\def\proof{\medskip\noindent{\bf Proof.\ }}
\def\proofof(#1){\medskip\noindent{\bf Proof of \csname c#1\endcsname.\ }}
\def\qed{\hfill{\sevenbf QED}\par\medskip}
\def\references{\bigskip\noindent\hbox{\bf References}\medskip
                \ifx\pdflink\undefined\else\bookmark{0}{References}\fi}
\def\addref#1{\expandafter\xdef\csname r#1\endcsname{\number\refcount}
    \global\advance\refcount by 1}

\def\nextremark #1\par{\item{$\circ$} #1}
\def\firstremark #1\par{\bigskip\noindent{\bf Remarks.}
     \smallskip\nextremark #1\par}
\def\abstract#1\par{{\baselineskip=10pt
    \eightpoint\narrower\noindent{\eightbf Abstract.} #1\par}}
%
%% EQUATION
\def\equtag#1{\expandafter\xdef\csname e#1\endcsname{(\number\seccount.\number\equcount)}
              \global\advance\equcount by 1}
\def\equation(#1){\equtag{#1}\eqno\csname e#1\endcsname}
\def\equ(#1){\hskip-0.03em\csname e#1\endcsname}
%
%% CLAIMS (theorems etc)
\def\clmtag#1#2{\expandafter\xdef\csname cn#2\endcsname{\number\seccount.\number\clmcount}
                \expandafter\xdef\csname c#2\endcsname{#1~\number\seccount.\number\clmcount}
                \global\advance\clmcount by 1}
\def\claim #1(#2) #3\par{\clmtag{#1}{#2}
    \vskip.1in\medbreak\noindent
    {\bf \csname c#2\endcsname .\ }{\sl #3}\par
    \ifdim\lastskip<\medskipamount
    \removelastskip\penalty55\medskip\fi}
\def\clm(#1){\csname c#1\endcsname}
\def\clmno(#1){\csname cn#1\endcsname}
%
%% SECTION
\def\sectag#1{\global\advance\seccount by 1
              \expandafter\xdef\csname sectionname\endcsname{\number\seccount. #1}
              \equcount=1 \clmcount=1 \subcount=1 \execount=0 \procount=0}
\def\section#1\par{\vskip0pt plus.1\vsize\penalty-40
    \vskip0pt plus -.1\vsize\bigskip\bigskip
    \sectag{#1}
    \message{\sectionname}\leftline{\magtenbf\sectionname}
    \nobreak\smallskip\noindent
    \ifx\pdflink\undefined
    \else
      \bookmark{0}{\sectionname}
    \fi}
%
%% SUBSECTION
\def\subtag#1{\expandafter\xdef\csname subsectionname\endcsname{\number\seccount.\number\subcount. #1}
              \global\advance\subcount by 1}
\def\subsection#1\par{\vskip0pt plus.05\vsize\penalty-20
    \vskip0pt plus -.05\vsize\medskip\medskip
    \subtag{#1}
    \message{\subsectionname}\leftline{\tenbf\subsectionname}
    \nobreak\smallskip\noindent
    \ifx\pdflink\undefined
    \else
      \bookmark{0}{.... \subsectionname}  %% can get a bit cluttered
    \fi}
%
%% DEMO (examples etc)
\def\demtag#1#2{\global\advance\demcount by 1
              \expandafter\xdef\csname de#2\endcsname{#1~\number\demcount}}
\def\demo #1(#2) #3\par{
  \demtag{#1}{#2}
  \vskip.1in\medbreak\noindent
  {\bf #1 \number\demcount.\enspace}
  {\rm #3}\par
  \ifdim\lastskip<\medskipamount
  \removelastskip\penalty55\medskip\fi}
\def\dem(#1){\csname de#1\endcsname}
%
%% EXERCISE
\def\exetag#1{\global\advance\execount by 1
              \expandafter\xdef\csname ex#1\endcsname{Exercise~\number\seccount.\number\execount}}
\def\exercise(#1) #2\par{
  \exetag{#1}
  \vskip.1in\medbreak\noindent
  {\bf Exercise \number\execount.}
  {\rm #2}\par
  \ifdim\lastskip<\medskipamount
  \removelastskip\penalty55\medskip\fi}
\def\exe(#1){\csname ex#1\endcsname}
%
%% PROBLEM
\def\protag#1{\global\advance\procount by 1
              \expandafter\xdef\csname pr#1\endcsname{\number\seccount.\number\procount}}
\def\problem(#1) #2\par{
  \ifnum\procount=0
    \parskip=6pt
    \vbox{\bigskip\centerline{\bf Problems \number\seccount}\nobreak\medskip}
  \fi
  \protag{#1}
  \item{\number\procount.} #2}
\def\pro(#1){Problem \csname pr#1\endcsname}
%
%macros.21
%
% requires fonts.5 or later
% also defines mathds (double strike) family
%
\def\rightheadline{\hfil}
\def\leftheadline{\sevenrm\hfil HANS KOCH\hfil}
\headline={\ifnum\pageno=\firstpage\hfil\else
\ifodd\pageno{{\fiverm\rightheadline}\number\pageno}
\else{\number\pageno\fiverm\leftheadline}\fi\fi}
\footline={\ifnum\pageno=\firstpage\hss\tenrm\folio\hss\else\hss\fi}
\let\ov=\overline
\let\cl=\centerline

\let\eps=\varepsilon
\let\sss=\scriptscriptstyle

\def\AA{{\cal A}}
\def\BB{{\cal B}}

\def\DD{{\cal D}}

\def\FF{{\cal F}}

\def\II{{\cal I}}
\def\JJ{{\cal J}}

\def\LL{{\cal L}}

\def\NN{{\cal N}}

\def\XX{{\cal X}}
\def\YY{{\cal Y}}

\def\rmC{{\rm C}}
\def\rmL{{\rm L}}
\def\id{{\rm I}}

\def\Im{\mathop{\rm Im}\nolimits}
%
%%%%%%%%%%%%%%
\newfam\dsfam
\def\mathds #1{{\fam\dsfam\tends #1}}

\font\tends=dsrom10
\font\eightds=dsrom8
\textfont\dsfam=\tends
\scriptfont\dsfam=\eightds
%%%%%%%%%%%%%%
%

\def\integer{{\mathds Z}}

\def\real{{\mathds R}}
\def\complex{{\mathds C}}

\def\mean{{\Bbb E}}
\def\proj{{\Bbb P}}
\def\torus{{\Bbb T}}
\def\iso{{\Bbb J}}

\def\bdot{\hbox{\bf .}}

\def\defeq{\mathrel{\mathop=^{\sss\rm def}}}
\def\half{{1\over 2}}

\def\quarter{{1\over 4}}
\def\thalf{{\textstyle\half}}

%

%

%

%

%
% from TeX book: used for commutative diagram
% in math mode, before using matrix, do
% \def\normalbaselines{\baselineskip20pt\lineskip3pt\lineskiplimit3pt}

%% modification of graphicx.tex by Nathan Goldschmidt
\input miniltx

\ifx\pdfoutput\undefined
  \def\Gin@driver{dvips.def}  % we are not running PDFTeX
\else
  \def\Gin@driver{pdftex.def} % we are running PDFTeX
\fi
 
\input graphicx.sty
\resetatcatcode
%% table macros from opmac.tex
%% Petr Olsak, 2012 -- 2016
%% http://petr.olsak.net/opmac.html

\newcount\tmpnum % auxiliary count
\newdimen\tmpdim % auxiliary dimen
\def\opwarning#1{\immediate\write16{l.\the\inputlineno\space OPmac WARNING: #1.}}
\long\def\addto#1#2{\expandafter\def\expandafter#1\expandafter{#1#2}}
\long\def\isinlist#1#2#3{\begingroup \long\def\tmp##1#2##2\end{\def\tmp{##2}%
   \ifx\tmp\empty \endgroup \csname iffalse\expandafter\endcsname \else
                  \endgroup \csname iftrue\expandafter\endcsname \fi}% end of \def\tmp
   \expandafter\tmp#1\endlistsep#2\end
}
\def\tabstrut{\strut}     % strut in the \table
\def\tabiteml{\enspace}   % left material before each \table item
\def\tabitemr{\enspace}   % right material after each \table item
\def\vvkern{1pt}          % space between vertical lines
\def\hhkern{1pt}          % space between horizontal lines

%%%%%%%%%%%%%% \table -- sec. 3.19 in opmac-d.pdf

\newtoks\tabdata
\def\tabstrutA{\tabstrut}
\newcount\colnum
\def\ddlinedata{}
\def\vvleft{}

\def\table{\vbox\bgroup \catcode`\|=12 \tableA}
\def\tableA#1#2{\offinterlineskip \colnum=0 \def\tmpa{}\tabdata={}\scantabdata#1\relax
   \halign\expandafter{\the\tabdata\cr#2\crcr}\egroup}

\def\scantabdata#1{\let\next=\scantabdata
   \ifx\relax#1\let\next=\relax
   \else\ifx|#1\addtabvrule
      \else\isinlist{123456789}#1\iftrue \def\next{\scantabdataC#1}%
          \else \expandafter\ifx\csname tabdeclare#1\endcsname \relax
                \expandafter\ifx\csname paramtabdeclare#1\endcsname \relax
                   \opwarning{tab-declarator "#1" unknown, ignored}%
                \else \def\next{\expandafter \scantabdataB \csname paramtabdeclare#1\endcsname}\fi
             \else \def\next{\expandafter\scantabdataA \csname tabdeclare#1\endcsname}%
   \fi\fi\fi\fi \next
}
\def\scantabdataA#1{\addtabitem \expandafter\addtabdata\expandafter{#1\tabstrutA}\scantabdata}
\def\scantabdataB#1#2{\addtabitem\expandafter\addtabdata\expandafter{#1{#2}\tabstrutA}\scantabdata}
\def\scantabdataC {\def\tmpb{}\afterassignment\scantabdataD \tmpnum=}
\def\scantabdataD#1{\loop \ifnum\tmpnum>0 \advance\tmpnum by-1 \addto\tmpb{#1}\repeat
   \expandafter\scantabdata\tmpb
}

\def\unsskip{\ifdim\lastskip>0pt \unskip\fi}
\def\addtabitem{\ifnum\colnum>0 \addtabdata{&}\addto\ddlinedata{&\dditem}\fi
    \advance\colnum by1 \let\tmpa=\relax}
\def\addtabdata#1{\tabdata\expandafter{\the\tabdata#1}}
\def\addtabvrule{%
    \ifx\tmpa\vrule \addtabdata{\kern\vvkern}%
       \ifnum\colnum=0 \addto\vvleft{\vvitem}\else\addto\ddlinedata{\vvitem}\fi
    \else \ifnum\colnum=0 \addto\vvleft{\vvitemA}\else\addto\ddlinedata{\vvitemA}\fi\fi
    \let\tmpa=\vrule \addtabdata{\vrule}}

\def\crll{\crcr\noalign{\hrule\kern\hhkern\hrule}}

\def\crli{\crcr \omit
   \gdef\dditem{\omit\tablinefil}\gdef\vvitem{\kern\vvkern\vrule}\gdef\vvitemA{\vrule}%
   \vvleft\tablinefil\ddlinedata\crcr}
\def\crlli{\crli\noalign{\kern\hhkern}\crli}
\def\tablinefil{\leaders\hrule\hfil}

\def\crlp#1{\crcr \noalign{\kern-\drulewidth}%
   \omit \xdef\crlplist{#1}\xdef\crlplist{,\expandafter}\expandafter\crlpA\crlplist,\end,%
   \global\tmpnum=0 \gdef\dditem{\omit\crlpD}%
   \gdef\vvitem{\kern\vvkern\kern\drulewidth}\gdef\vvitemA{\kern\drulewidth}%
   \vvleft\crlpD\ddlinedata \global\tmpnum=0 \crcr}
\def\crlpA#1,{\ifx\end#1\else \crlpB#1-\end,\expandafter\crlpA\fi}
\def\crlpB#1#2-#3,{\ifx\end#3\xdef\crlplist{\crlplist#1#2,}\else\crlpC#1#2-#3,\fi}
\def\crlpC#1-#2-#3,{\tmpnum=#1\relax
   \loop \xdef\crlplist{\crlplist\the\tmpnum,}\ifnum\tmpnum<#2\advance\tmpnum by1 \repeat}
\def\crlpD{\global\advance\tmpnum by1
   \edef\tmpa{\noexpand\isinlist\noexpand\crlplist{,\the\tmpnum,}}%
   \tmpa\iftrue \kern-\drulewidth \tablinefil \kern-\drulewidth\else\hfil \fi}

\def\tskip{\afterassignment\tskipA \tmpdim}
\def\tskipA{\gdef\dditem{}\gdef\vvitem{}\gdef\vvitemA{}\gdef\tabstrutA{}%
    \vbox to\tmpdim{}\ddlinedata \crcr \noalign{\gdef\tabstrutA{\tabstrut}}}

\def\mspan{\omit \tabdata={\tabstrut}\let\tmpa=\relax \afterassignment\mspanA \mscount=}
\def\mspanA[#1]{\loop \ifnum\mscount>1 \csname span\endcsname \omit \advance\mscount by-1 \repeat
   \mspanB#1\relax}
\def\mspanB#1{\ifx\relax#1\def\tmpa{\def\tmpa####1}%
   \expandafter\tmpa\expandafter{\the\tabdata\ignorespaces}\expandafter\tmpa\else
   \ifx |#1\ifx\tmpa\vrule\addtabdata{\kern\vvkern}\fi \addtabdata{\vrule}\let\tmpa=\vrule
   \else \let\tmpa=\relax
      \ifx c#1\addtabdata{\tabiteml\hfil\ignorespaces##1\unsskip\hfil\tabitemr}\fi
      \ifx l#1\addtabdata{\tabiteml\ignorespaces##1\unsskip\hfil\tabitemr}\fi
      \ifx r#1\addtabdata{\tabiteml\hfil\ignorespaces##1\unsskip\tabitemr}\fi
   \fi \expandafter\mspanB \fi}

\newdimen\drulewidth  \drulewidth=0.4pt
\let\orihrule=\hrule  \let\orivrule=\vrule
\def\rulewidth{\afterassignment\rulewidthA \drulewidth}
\def\rulewidthA{\edef\hrule{\orihrule height\the\drulewidth}%
                \edef\vrule{\orivrule width\the\drulewidth}}

\long\def\frame#1{%
   \hbox{\vrule\vtop{\vbox{\hrule\kern\vvkern
      \hbox{\kern\hhkern\relax#1\kern\hhkern}%
   }\kern\vvkern\hrule}\vrule}}
%% boldmath

%=================================================================
% necessary fonts
%-----------------------------------------------------------------
   \font\Fivebf  =cmbx10  scaled 500 % five point bold
   \font\Sevenbf =cmbx10  scaled 700 % seven point bold
   \font\Tenbf   =cmbx10             % ten point bold
   \font\Fivemb  =cmmib10 scaled 500 % five point math bold
   \font\Sevenmb =cmmib10 scaled 700 % seven point math bold
   \font\Tenmb   =cmmib10            % ten point math bold
%=================================================================
% Math in Bold
% example $\boldmath{..}$   { math characters will be bold }
%-----------------------------------------------------------------
\def\boldmath{\textfont0=\Tenbf           \scriptfont0=\Sevenbf 
              \scriptscriptfont0=\Fivebf  \textfont1=\Tenmb
              \scriptfont1=\Sevenmb       \scriptscriptfont1=\Fivemb}
%=================================================================

%% twelveboldmath

%=================================================================
% necessary fonts
%-----------------------------------------------------------------
%  \font\Sevenbf   =cmbx10  scaled 700

%  \font\Sevenmb   =cmmib10 scaled 700

%=================================================================
% Math in Bold
% example $\twelveboldmath{...}$   { math characters will be bold }
%-----------------------------------------------------------------

%=================================================================

%\input param.2
%\input fonts.6
%\input smallfonts.tex
%\input symbols.1
%\input links.1
%\input titles.9
%\input macros.21
%\input boldmath.tex
%\input mygraphicx.tex
%
\newdimen\savedparindent
\savedparindent=\parindent
\font\tenamsb=msbm10 \font\sevenamsb=msbm7 \font\fiveamsb=msbm5
\newfam\bbfam
\textfont\bbfam=\tenamsb
\scriptfont\bbfam=\sevenamsb
\scriptscriptfont\bbfam=\fiveamsb
\def\field{{\mathds F}}
\def\bbb{\fam\bbfam}
\def\oldinteger{{\bbb Z}}
\def\oldnatural{{\bbb N}}
\def\buB{{\hbox{\magnineeufm B}}}
\let\Beth\theta
\def\gapi{\hskip 2pt}
\def\gapii{\hskip 2pt}
\def\rmc{{\rm c}}
\def\rme{{\rm e}}
\def\id{{\rm I}}
\def\LOP{{\mathds L}}
\def\nsv{\mathop{\not\!\partial}\nolimits}
\def\cosi{\mathop{\rm cosi}\nolimits}
\def\even{{\rm e}}
\def\odd{{\rm o}}
\def\bdot{\hbox{\bf .}}
\def\stwomat#1#2#3#4{{\eightpoint\left[\matrix{#1&#2\cr#3&#4\cr}\right]}}
\addref{Hopf}
\addref{Serrin}
\addref{CrAb}
\addref{MarMcC}
\addref{RT}
\addref{KlWe}
\addref{ChIo}
\addref{Dum}
\addref{ChIoo}
\addref{NWYNK}
\addref{AKxii}
\addref{AKxv}
\addref{Galdi}
\addref{HJNS}
\addref{AKxviii}
\addref{AKxix}
\addref{NPW}
\addref{GoSe}
\addref{WiZg}
\addref{BBLV}
\addref{AGK}
\addref{BLQ}
\addref{BQ}
\addref{Progs}
\addref{Ada}
\addref{Gnat}
\addref{IEEE}
\addref{MPFR}
\def\leftheadline{\sixrm\hfil ARIOLI \& KOCH\hfil}
\def\rightheadline{\sevenrm\hfil Hopf bifurcation in Navier-Stokes\hfil}
%
%%%%%%%%%%%%%%%%%%%%%%%%%%%%%%%%%%%%%%%%%%%%%%%%%%%%%%%
\cl{\magtenbf A Hopf bifurcation in the planar Navier-Stokes equations}
\bigskip

\cl{
Gianni Arioli
\footnote{$^1$}
{\eightpoint\hskip-2.9em
Department of Mathematics, Politecnico di Milano,
Piazza Leonardo da Vinci 32, 20133 Milano.
}
$^{\!\!\!,\!\!}$
\footnote{$^2$}
{\eightpoint\hskip-2.6em
Supported in part by the PRIN project
``Equazioni alle derivate parziali e
disuguaglianze analitico-geometriche associate''.}
and Hans Koch
\footnote{$^3$}
{\eightpoint\hskip-2.7em
Department of Mathematics, The University of Texas at Austin,
Austin, TX 78712.}
}

\bigskip
\abstract
We consider the Navier-Stokes equation
for an incompressible viscous fluid on a square,
satisfying Navier boundary conditions
and being subjected to a time-independent force.
As the kinematic viscosity is varied,
a branch of stationary solutions is shown
to undergo a Hopf bifurcation,
where a periodic cycle branches from the stationary solution.
Our proof is constructive
and uses computer-assisted estimates.

%%%%%%%%%%%%%%%%%%%%%%%%%%%%%%%%%%%%%
%%%%%%%%%%%%%%%%%%%%%%%%%%%%%%%%%%%%%
\section Introduction and main result
%%%%%%%%%%%%%%%%%%%%%%%%%%%%%%%%%%%%%
%%%%%%%%%%%%%%%%%%%%%%%%%%%%%%%%%%%%%

We consider the Navier-Stokes equations
$$
\partial_t u-\nu\Delta u+(u\cdot\nabla)u+\nabla p=f\ ,\quad
\nabla\cdot u=0\quad{\rm on~}\Omega\,,
\equation(NavierStokes)
$$
for the velocity $u=u(t,x,y)$
of an incompressible fluid on a planar domain $\Omega$,
satisfying suitable boundary conditions for $(x,y)\in\partial\Omega$
and initial conditions at $t=0$.
Here, $p$ denotes the pressure,
and $f=f(x,y)$ is a fixed time-independent external force.

Our focus is on solution curves and bifurcations
as the kinematic velocity $\nu$ is being varied.
In order to reduce the complexity of the problem,
the domain $\Omega$ is chosen to be as simple as possible,
namely the square $\Omega=(0,\pi)^2$.
Following [\rAGK], we impose Navier boundary conditions on $\partial\Omega$,
which are given by
$$
\eqalign{
u_1&=\partial_x u_2=0\quad{\rm on~}\{0,\pi\}\times(0,\pi)\,,\cr
u_2&=\partial_y u_1=0\quad{\rm on~}(0,\pi)\times\{0,\pi\}\,.\cr}
\equation(NavierBC)
$$
A fair amount is known about the (non)uniqueness of stationary solutions
in this case [\rAGK].
This includes the existence of a bifurcation between curves of stationary solutions
with different symmetries.

Here we prove the existence of a Hopf bifurcation
for the equation \equ(NavierStokes) with boundary conditions \equ(NavierBC),
and with a forcing function $f$ that satisfies
$$
(\partial_xf_2-\partial_yf_1)(x,y)=5\sin(x)\sin(2y)-13\sin(3x)\sin(2y)\,.
\equation(specificNScurlf)
$$
In a Hopf bifurcation, a stationary solution loses stability
and a small-amplitude limit cycle branches from the stationary solution [\rHopf,\rCrAb,\rMarMcC].
Among other things, this introduces a time scale in the system
and increases its complexity.
In this capacity, Hopf bifurcations in the Navier-Stokes equation
constitute an important first step
in the transition to turbulence in fluids,
as was described in the seminal work [\rRT].

Numerically, there is plenty of evidence that Hopf bifurcations occur
in the Navier-Stokes equation, but proofs are still very scarce.
An explicit example of a Hopf bifurcation
was given in [\rKlWe] for the rotating B\'enard problem.
A proof exists also for the Couette-Taylor problem [\rChIo,\rChIoo].
Sufficient conditions for the existence of a Hopf bifurcation
in a Navier-Stokes setting are presented in [\rGaldi].

Before giving a precise statement of our result, let us replace
the vector field $u$ in the equation \equ(NavierStokes) by $\nu^{-1}u$.
The equation for the rescaled function $u$ is
$$
\alpha\partial_t u-\Delta u+\gamma(u\cdot\nabla)u+\nabla p=f\ ,\quad
\nabla\cdot u=0\quad{\rm on~}\Omega\,,
\equation(NS)
$$
where $\gamma=\nu^{-2}$.
The value of $\alpha$ that corresponds to \equ(NavierStokes) is $\nu^{-1}$,
but this can be changed to any positive value by rescaling time.

Numerically, it is possible to find stationary solutions of \equ(NS)
for a wide range of values of the parameter $\gamma$.
At a value $\gamma_0\approx 83.1733117\ldots$ we observe
a Hopf bifurcation that leads to a branch of periodic solutions for $\gamma>\gamma_0$.

For a fixed value of $\alpha$,
the time period $\tau$ of the solution varies with $\gamma$.
Instead of looking for $\tau$-periodic solution of \equ(NS)
for fixed $\alpha$, we look for $2\pi$-periodic solutions,
where $\alpha=2\pi/\tau$ has to be determined.
To simplify notation, a $2\pi$-periodic function
will be identified with a function on the circle $\torus=\real/(2\pi\integer)$.
Our main theorem is the following.

\claim Theorem(NSHopfBif)
There exists a real number $\gamma_0=83.1733117\ldots$,
an open interval $I$ including $\gamma_0$,
and a real analytic function $(\gamma,x,y)\mapsto u_\gamma(x,y)$
from $I\times\Omega$ to $\real^2$,
such that $u_\gamma$ is a stationary solution of \equ(NS) and \equ(NavierBC)
for each $\gamma\in I$.
In addition, there exists a real number $\alpha_0=4.66592275\ldots$,
an open interval $J$ centered at the origin,
two real analytic functions $\gamma$ and $\alpha$ on $J$
that satisfy $\gamma(0)=\gamma_0$ and $\alpha(0)=\alpha_0$, respectively,
as well as two real analytic functions
$(s,t,x,y)\mapsto u_{s,\even}(t,x,y)$ and $(s,t,x,y)\mapsto u_{s,\odd}(t,x,y)$
from $J\times\torus\times\Omega$
to $\real^2$, such that the following holds.
For any given $\beta\in\complex$ satisfying $\beta^2\in J$,
the vector field $u=u_{s,\even}+\beta u_{s,\odd}$ with $s=\beta^2$
is a solution of \equ(NS) and \equ(NavierBC)
with $\gamma=\gamma(s)$ and $\alpha=\alpha(s)$.
Furthermore, $u_{0,\even}(t,\bdot\,,\bdot)=u_{\gamma_0}$
and $\partial_t u_{0,\odd}(t,\bdot\,,\bdot)\ne 0$.

To our knowledge, this is the first result establishing
the existence of a Hopf bifurcation for the Navier-Stokes equation
in a stationary environment.

Our proof of this theorem is computer-assisted.
The solutions are obtained by rewriting \equ(NS) and \equ(NavierBC)
as a suitable fixed point equation for scalar vorticity of $u$.
Here we take advantage of the fact that the domain is two-dimensional.
We isolate the periodic branch from the stationary branch
by using a scaling that admits two distinct limits at the bifurcation point.
This approach is also known as the blow-up method,
which is a common tool in the study of singularities and bifurcations [\rDum].

\smallskip
Computer-assisted methods have been applied successfully
to many different problems in analysis,
mostly in the areas of dynamical systems and partial differential equations.
Here we will just mention work that
concerns the Navier-Stokes equation or Hopf bifurcations.
For the Navier-Stokes equation,
the existence of symmetry-breaking bifurcations among stationary solutions
has been established in [\rNWYNK,\rAGK].
Periodic solutions for the Navier-Stokes flow in a stationary environment
have been obtained in [\rBBLV].
In the case of periodic forcing,
the problem of existence and stability of periodic orbits
has been investigated in [\rHJNS].
Concerning the existence of Hopf bifurcations,
a computer-assisted proof was given recently in [\rBLQ]
for a finite-dimensional dynamical system;
and an extension of their method to the Kuramoto-Sivashinsky PDE
is presented in [\rBQ].
For other recent computer-assisted proofs
we refer to [\rAKxix,\rNPW,\rGoSe,\rWiZg] and references therein.

Figure 1 depicts snapshots at $t=0$ and $t=\pi$
of a solution $u:\torus\times\Omega\to\real^2$ of the
equations \equ(NS) with boundary conditions \equ(NavierBC)
and forcing \equ(specificNScurlf),
obtained numerically for the parameter value $\gamma\approx 84.00\ldots$.

%%%%%%%%%%%%%%%%%%%%%%%%%%%%%%%%%%%%%%%%%%%%%%%%%%%%%%%%%%%%%
\vskip0.15in
\hbox{\hskip 55pt
\includegraphics[height=1.5in,width=1.5in]{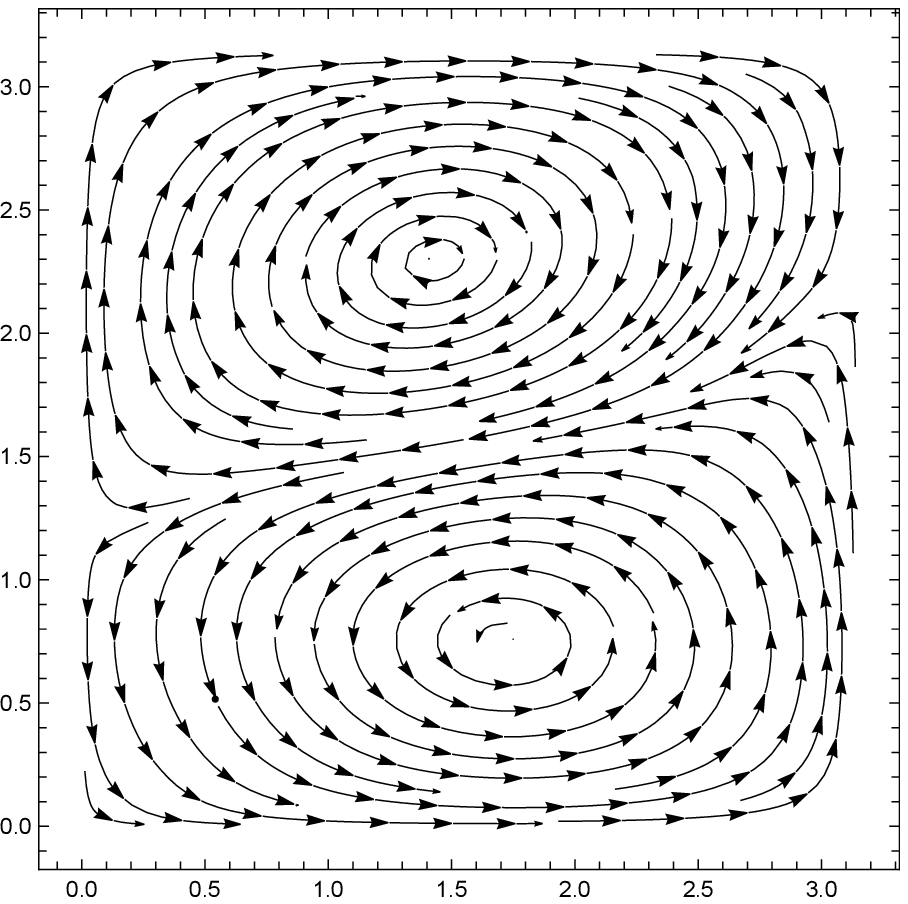}
\hskip 35pt
\includegraphics[height=1.5in,width=1.5in]{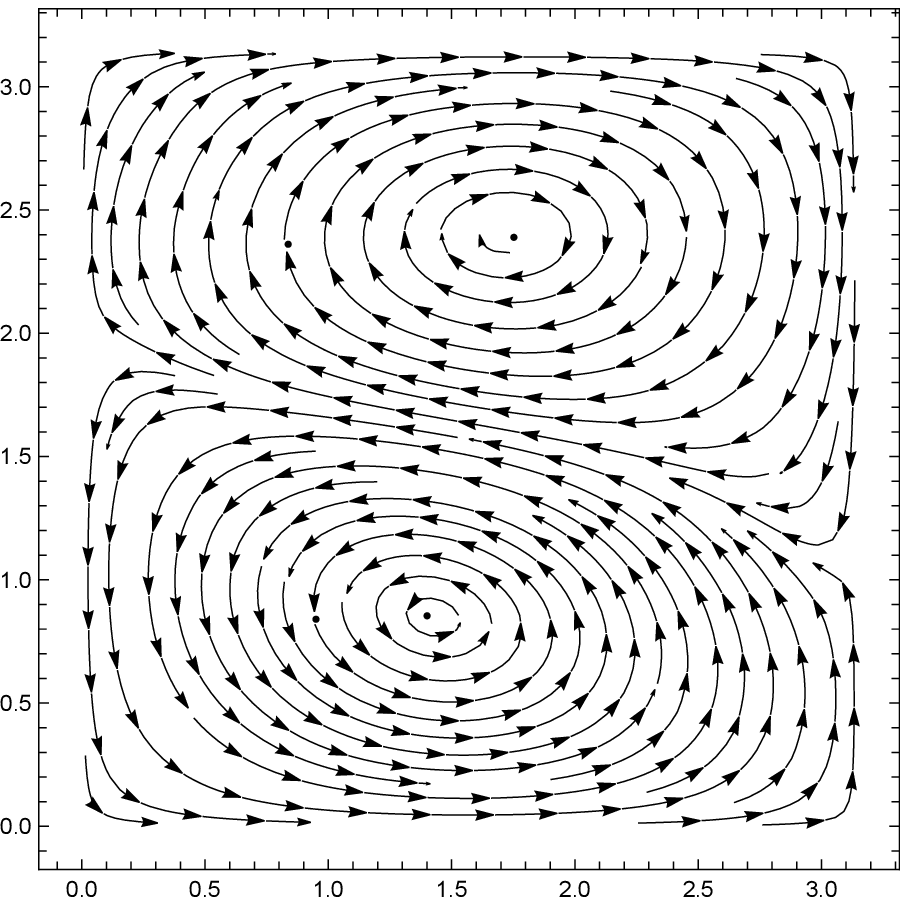}}
\vskip0.1in
\centerline{\eightpoint {\noindent\bf Figure 1.}
Snapshots at two distinct times of a time-periodic solution
for $\gamma\approx 84.00\ldots$
}
\vskip0.15in
%%%%%%%%%%%%%%%%%%%%%%%%%%%%%%%%%%%%%%%%%%%%%%%%%%%%%%%%%%%%%

As mentioned earlier, a system similar to the one
considered here is known to exhibit a symmetry-breaking bifurcation
within the class of stationary solutions [\rAGK].
The broken symmetry is $y\mapsto\pi/2-y$.
Based on a numerical computation of eigenvalues,
we expect an analogous bifurcation to occur here at $\gamma\approx 1450$.
Interestingly, the Hopf bifurcation described here occurs
at a significantly smaller value of $\gamma$.
We have not tried to prove the existence
of a symmetry-breaking bifurcation for the forcing \equ(specificNScurlf),
since such an analysis would duplicate the work in [\rAGK]
and go beyond the scope of the present paper.

\smallskip
The remaining part of this paper is organized as follows.
In Section 2, we first rewrite \equ(NS) as an equation for
the function $\Phi=\partial_y u_1-\partial_x u_2$,
which is the scalar vorticity of $-u$.
After a suitable scaling $\Phi=U_\beta\phi$,
the problem of constructing the solution branches
described in \clm(NSHopfBif) is reduced to three fixed point problems
for the function $\phi$.
These fixed point equations are solved in Section 3,
based on estimates described in Lemmas 3.3, 3.4, and 3.6.
Section 4 is devoted to the proof of these estimates,
which involves reducing them to a large number of trivial bounds
that can be (and have been) verified with the aid of a computer [\rProgs].

%%%%%%%%%%%%%%%%%%%%%%%%%%%%%%
%%%%%%%%%%%%%%%%%%%%%%%%%%%%%%
\section Fixed point equations
%%%%%%%%%%%%%%%%%%%%%%%%%%%%%%
%%%%%%%%%%%%%%%%%%%%%%%%%%%%%%

The goal here is to rewrite the equation \equ(NS)
with boundary conditions \equ(NavierBC) as a fixed point problem.
Applying the operator $\nsv:(u_1,u_2)\mapsto\partial_2 u_1-\partial_1 u_2$
on both sides of the equation \equ(NS), we obtain
$$
\alpha\partial_t\Phi
-\Delta\Phi+\gamma u\cdot\nabla\Phi=\nsv f\,,\qquad\Phi=\nsv u\,.
\equation(curlNS)
$$
Here, we have used that $\nsv\,(u\cdot\nabla)u=u\cdot\nabla\Phi$.
Using the divergence-free condition $\nabla\cdot u=0$,
one also finds that
$$
\Delta u=\iso\nabla\Phi\,,\qquad
\iso=\stwomat{0}{1}{-1}{0}\,.
\equation(Lapu)
$$
If $\Phi$ vanishes on the boundary of $\partial\Omega$,
then the equation \equ(Lapu) can be inverted to yield
$$
u=\nsv^{-1}\Phi\defeq\iso\nabla\Delta^{-1}\Phi\,,
\equation(curlInv)
$$
where $\Delta$ denotes the Dirichlet Laplacean on $\Omega$.

In Section 3
we will define a space of real analytic functions $\Phi$
that admit a representation
$$
\Phi(t,x,y)=\sum_{j,k\in\oldnatural_1}\Phi_{j,k}(t)\sin(jx)\sin(ky)\,,
\equation(PhixyExpansion)
$$
with the series converging uniformly on a complex
open neighborhood of $\torus^3$.
Here, and in what follows, $\oldnatural_1$ denotes the set of all positive integers.
If $\Phi$ admits such an expansion,
then the equation \equ(curlInv) yields
$$
\eqalign{
u_1(t,x,y)&=\sum_{j,k\in\oldnatural_1}\,{-k\over j^2+k^2}\Phi_{j,k}(t)\sin(jx)\cos(ky)\,,\cr
u_2(t,x,y)&=\sum_{j,k\in\oldnatural_1}\,{j\over j^2+k^2}\Phi_{j,k}(t)\cos(jx)\sin(ky)\,.\cr}
\equation(uxyExpansion)
$$
It is straightforward to check that the corresponding
vector field $u=(u_1,u_1)$ satisfies the Navier boundary conditions \equ(NavierBC).
So a solution $u$ of \equ(NS) and \equ(NavierBC)
can be obtained via \equ(uxyExpansion) from a solution $\Phi$
of the equation \equ(curlNS).
For convenience, we write \equ(curlNS) as
$$
(\alpha\partial_t-\Delta)\Phi+\thalf\gamma\LOP(\Phi)\Phi=\nsv f\,,
\equation(scNS)
$$
where $\LOP$ is the symmetric bilinear form defined by
$$
\LOP(\phi)\psi=
(\nabla\phi)\cdot\nsv^{-1}\psi
+(\nabla\psi)\cdot\nsv^{-1}\phi\,.
\equation(LOPDef)
$$
The coefficients $\Phi_{j,k}$ in the series \equ(PhixyExpansion)
are $2\pi$-periodic functions and thus admit an expansion
$$
\Phi_{j,k}=\sum_{n\in\oldinteger}\Phi_{n,j,k}\cosi_n\,,\qquad
\cosi_n(t)=\cases{\cos(nt) &if $n\ge 0$,\cr\sin(-nt) &if $n<0$.\cr}
\equation(PhijkSeries)
$$
Denote by $\oldnatural_0$ the set of all nonnegative integers.
For any subset $N\subset\oldnatural_0$ we define
$$
\mean_N\Phi=\sum_{n\in\oldinteger\atop|n|\in N}\sum_{j,k\in\oldnatural_1}
\Phi_{n,j,k}\cosi_n\times\sin_j\times\sin_k\,,
\equation(FreqProj)
$$
where $\sin_m(z)=\sin(mz)$.
In particular, the even frequency part $\Phi_\even$
(odd frequency part $\Phi_\odd$) of $\Phi$ is defined
to be the function $\mean_N\Phi$,
where $N$ is the set of all even (odd) nonnegative integers.
This leads to the decomposition $\Phi=\Phi_\even+\Phi_\odd$
that will be used below.

To simplify the discussion,
consider first non-stationary periodic solutions.
For $\gamma$ near the bifurcation point $\gamma_0$,
we expect $\Phi$ to be nearly time-independent.
So in particular, $\Phi_\odd$ is close to zero.
Consider the function $\phi=\phi_\even+\phi_\odd$
obtained by setting $\phi_\even=\Phi_\even$ and $\phi_\odd=\beta^{-1}\Phi_\odd$.
The scaling factor $\beta\ne 0$ will be chosen below,
in such a way that $\phi_\even$ and $\phi_\odd$ are of comparable size.
Substituting
$$
\Phi=U_\beta\phi\defeq\phi_\even+\beta\phi_\odd
\equation(PhiDecomp)
$$
into \equ(scNS) yields the equation
$$
(\alpha\partial_t-\Delta)\phi+\thalf\gamma\LOP_s(\phi)\phi=\nsv f\,,
\equation(scNSs)
$$
where $s=\beta^2$ and
$$
\LOP_s(\phi)\psi
=\LOP(\phi_\even)\psi_\even+\LOP(\phi_\even)\psi_\odd
+\LOP(\phi_\odd)\psi_\even+s\LOP(\phi_\odd)\psi_\odd\,.
\equation(LOPsDef)
$$
Finally, we convert \equ(scNSs) to a fixed point equation
by applying the inverse of $\alpha\partial_t-\Delta$ to both sides.
Setting $g=(-\Delta)^{-1}\nsv f$,
the resulting equation is $\tilde\phi=\phi$, where
$$
\tilde\phi=g
-\thalf\gamma|\Delta|^{1/2}(\alpha\partial_t-\Delta)^{-1}\hat\phi\,,\qquad
\hat\phi\defeq|\Delta|^{-1/2}\LOP_s(\phi)\phi\,.
\equation(scNSsFix)
$$

One of the features of the equation \equ(scNSs)
is that  the time-translate of a solution is again a solution.
We eliminate this symmetry by imposing the condition $\phi_{1,1,1}=0$.
In addition, we choose $\beta=\Beth^{-1}\Phi_{-1,1,1}$,
where $\Beth$ is some fixed constant that will be specified later.
This leads to the normalization conditions
$$
A\phi\defeq\phi_{1,1,1}=0\,,\qquad
B\phi\defeq\phi_{-1,1,1}=\Beth\,.
\equation(ABNormaliz)
$$
Notice that $\beta$ enters our main equation $\tilde\phi=\phi$
only via its square $s=\beta^2$.
It is convenient to regard $s$ to be the independent parameter
and express $\gamma$ as a function of $s$.
The functions $\gamma=\gamma(s)$ and $\alpha=\alpha(s)$
are determined by the condition that $\tilde\phi$
satisfies the normalization conditions \equ(ABNormaliz).
Applying the functionals $A$ and $B$ to both sides of \equ(scNSs),
using the identities
$A\Delta=-2A$, $A\partial_t=B$, $B\Delta=-2B$, $B\partial_t=-A$,
and imposing the conditions $A\tilde\phi=0$ and $B\tilde\phi=\Beth$,
we find that
$$
\gamma=-2^{3/2}{\Beth\over B\hat\phi}\,,\qquad
\alpha=2{A\hat\phi\over B\hat\phi}\,.
\equation(gammaalphaSol)
$$
For a fixed value of $s$, define $\FF_s(\phi)=\tilde\phi$,
where $\tilde\phi$ is given by \equ(scNSsFix),
with $\gamma=\gamma(s,\phi)$ and $\alpha=\alpha(s,\phi)$
determined by \equ(gammaalphaSol).
The fixed point equation for $\FF_s$ is used
to find non-stationary time-periodic solutions of \equ(scNSs).

\demo Remark(NotNormalized)
The choice \equ(gammaalphaSol) guarantees that
$A\tilde\phi=0$ and $B\tilde\phi=\Beth$,
even if $\phi$ does not satisfy the normalization conditions \equ(ABNormaliz).
Thus, the domain of the map $\FF_s$
can include non-normalized function $\phi$.
(The same is true for the map $\FF_\gamma$ described below.)
But a fixed point of this map will be normalized by construction.

In order to determine the bifurcation point $\gamma_0$
and the corresponding frequency $\alpha_0$,
we consider the map $\FF:\phi\mapsto\tilde\phi$ given by \equ(scNSsFix) with $s=0$.
The values of $\gamma$ and $\alpha$ are again given by \equ(gammaalphaSol),
so that $A\tilde\phi=0$ and $B\tilde\phi=\Beth$.
We will show that this map $\FF$ has a fixed point
$\phi$ with the property that $\phi_{n,j,k}=0$ whenever $|n|>1$.
The values of $\gamma$ and $\alpha$ for this fixed point
define $\gamma_0$ and $\alpha_0$.

A similar map $\FF_\gamma:\phi\mapsto\tilde\phi$,
given by \equ(scNSsFix) with $s=0$,
is used to find stationary solutions of the equation \equ(scNS).
In this case, the value of $\gamma$ is being fixed,
and $\phi_\odd$ is taken to be zero.
The goal is to show that this map $\FF_\gamma$
has a fixed point $\phi_\gamma$ that is independent of time $t$.
Then $\Phi=\phi_\gamma$ is a stationary solution of \equ(scNS).

\medskip
We finish this section by computing
the derivative of the map $\FF_s$ described after \equ(gammaalphaSol).
The resulting expressions will be needed later.
Like some of the above, the following is purely formal.
A proper formulation will be given in the next section.
For simplicity, assume that $\phi$ depends on a parameter.
The derivative of a quantity $q$ with respect to this parameter
will be denoted by $\dot q$.
Define
$$
\LL_\alpha=|\Delta|^{1/2}(\alpha\partial_t-\Delta)^{-1}\,,\qquad
\LL_\alpha'=\partial_t(\alpha\partial_t-\Delta)^{-1}\,.
\equation(LLalphaDef)
$$
Using that $\FF_s(\phi)=g-\thalf\gamma\LL_\alpha\hat\phi$
with $\hat\phi=|\Delta|^{-1/2}\LOP_s(\phi)\phi$,
the parameter-derivative of $\FF_s(\phi)$ is given by
$$
D\FF_s(\phi)\dot\phi
=-\thalf\LL_\alpha\Bigl[\bigl(\dot\gamma-\gamma\dot\alpha\LL_\alpha'\bigr)\hat\phi
+\gamma\dot{\hat\phi}\Bigr]\,,\qquad
\dot{\hat\phi}=2|\Delta|^{-1/2}\LOP_s(\phi)\dot\phi\,,
\equation(DFFsPhiDotPhi)
$$
where
$$
\dot\gamma
=2^{-3/2}{\gamma^2\over\Beth} B\dot{\hat\phi}\,,\qquad
\dot\alpha
=2^{-3/2}{\alpha\gamma\over\Beth}B\dot{\hat\phi}
-2^{-1/2}{\gamma\over\Beth}A\dot{\hat\phi}\,.
\equation(DotgammaDotalpha)
$$
The above expressions for $\dot\gamma$ and $\dot\alpha$
are obtained by differentiating \equ(gammaalphaSol).

%%%%%%%%%%%%%%%%%%%%%%%%%%%%%%%%%%%%
%%%%%%%%%%%%%%%%%%%%%%%%%%%%%%%%%%%%
\section The associated contractions
%%%%%%%%%%%%%%%%%%%%%%%%%%%%%%%%%%%%
%%%%%%%%%%%%%%%%%%%%%%%%%%%%%%%%%%%%

In this section, we formulate the fixed point problems
for the maps $\FF$, $\FF_\gamma$, and $\FF_s$
in a suitable functional setting.
The goal is to reduce the problems to a point
where we can invoke the contraction mapping theorem.
After describing the necessary estimates,
we give a proof of \clm(NSHopfBif) based on these estimates.

We start by defining suitable function spaces.
Given a real number $\rho>1$, denote by $\AA$
the space of all functions $h\in\rmL^2(\torus)$
that have a finite norm $\|h\|$, where
$$
\|h\|=|h_0|
+\sum_{n\in\oldnatural_1}\sqrt{|h_n|^2+|h_{-n}|^2}\rho^n\,,\qquad
h=\sum_{n\in\oldinteger}h_n\cosi_n\,.
\equation(AADef)
$$
Here $\cosi_n$ are the trigonometric function defined in \equ(PhijkSeries).
It is straightforward to check that $\AA$
is a Banach algebra under the pointwise product of functions.
That is, $\|gh\|\le\|g\|\|h\|$ for any two functions $g,h\in\AA$.
We also identify functions on $\torus$ with $2\pi$-periodic functions on $\real$.
In this sense, a function in $\AA$
extends analytically to the strip $T(\rho)=\{z\in\complex:|\Im z|<\log\rho\}$.

Given in addition $\varrho>1$, we denote by $\buB$
the space of all function $\Phi:\torus^2\to\AA$
that admit a representation \equ(PhixyExpansion)
and have a finite norm
$$
\|\Phi\|
=\sum_{j,k\in\oldnatural_1}\|\Phi_{j,k}\|\varrho^{j+k}\,.
\equation(GGrhoepsNorm)
$$
A function $(x,y)\mapsto(t\mapsto\Phi(t,x,y))$ in this space
will also be identified with a function
$(t,x,y)\mapsto\Phi(t,x,y)$ on $\torus^3$,
or with a function on $\real^3$ that is $2\pi$-periodic in each argument.
In this sense, every function in $\buB$
extends analytically to $T(\rho)\times T(\varrho)^2$.

We consider $\AA$ and $\buB$ to be Banach spaces
over $\field\in\{\real,\complex\}$.
In the case $\field=\real$, the functions in these spaces
are assumed to take real values for real arguments.

\smallskip
Clearly, a function $\Phi\in\buB$
admits an expansion \equ(FreqProj) with $N=\oldnatural_0$.
The sequence of Fourier coefficients $\Phi_{n,k,j}$
converges to zero exponentially as $|n|+j+k$ tends to infinity.
If all but finitely many of these coefficients vanish,
then $\Phi$ is called a Fourier polynomial.
The equation \equ(FreqProj) with $N\subset\oldnatural_0$ non-empty
defines a continuous projection $\mean_N$ on $\buB$
whose operator norm is $1$.
Using Fourier series, it is straightforward
to see that the equation \equ(LLalphaDef)
defines two bounded linear operators $\LL_\alpha$ and $\LL_\alpha'$
on $\buB$, for every $\alpha\in\complex$.
The operator $\LL_\alpha$ is in fact compact.
Specific estimates will be given in Section 4.
The following will be proved in Section 4 as well.

\claim Proposition(LOPBound)
If $\Phi$ and $\phi$ belong to $\buB$,
then so does $|\Delta|^{-1/2}\LOP(\Phi)\phi$, and
$$
\bigl\||\Delta|^{-1/2}\LOP(\Phi)\phi\bigr\|
\le\bigl\||\Delta|^{-1/2}\Phi\bigr\|\|\phi\|
+\|\Phi\|\bigl\||\Delta|^{-1/2}\phi\bigr\|\,.
\equation(LOPBound)
$$

This estimate implies e.g.~that the transformation $\phi\mapsto\tilde\phi$,
given by \equ(scNSsFix) for fixed values of $s$, $\gamma$ and $\alpha$,
is well-defined and compact as a map from $\buB$ to $\buB$.

As is common in computer-assisted proofs,
we reformulate the fixed point equation for the map $\phi\mapsto\tilde\phi$
as a fixed point problem for an associated quasi-Newton map.
Since we need three distinct versions of this map,
let us first describe a more general setting.

\medskip
Let $\FF:\DD\to\BB$ be a $\rmC^1$ map
defined on an open domain $\DD$ in a Banach space $\BB$.
Let $h\mapsto\varphi+Lh$ be a continuous affine map on $\BB$.
We define quasi-Newton map $\NN$ for $(\DD,\FF,\varphi,L)$
by setting
$$
\NN(h)=\FF(\varphi+Lh)-\varphi+(\id-L)h\,.
\equation(quasiNewtonMap)
$$
The domain of $\NN$ is defined to be the set of of all $h\in\BB$
with the property that $\varphi+Lh\in\DD$.
Notice that, if $h$ is a fixed point of $\NN$,
then $\varphi+Lh$ is a fixed point of $\FF$.
In our applications, $\varphi$ is an approximate fixed point of $\FF$
and $L$ is an approximate inverse of $\id-D\FF(\varphi)$.

The following is an immediate consequence of the contraction mapping theorem.

\claim Proposition(ContrMappingThm)
Let $\FF:\DD\to\BB$ be a $\rmC^1$ map
defined on an open domain in a Banach space $\BB$.
Let $h\mapsto\varphi+Lh$ be a continuous affine map on $\BB$.
Assume that the quasi-Newton map \equ(quasiNewtonMap)
includes a non-empty ball $B_\delta=\{h\in\BB: \|h\|<\delta\}$ in its domain,
and that
$$
\|\NN(0)\|<\eps\,,\qquad\|D\NN(h)\|<K\,,\qquad h\in B_\delta\,,
\equation(ContrMappingThm)
$$
where $\eps,K$ are positive real numbers that satisfy $\eps+K\delta<\delta$.
Then $\FF$ has a fixed point in $\varphi+LB_\delta$.
If $L$ is invertible, then this fixed point is unique in $\varphi+LB_\delta$.

In our applications below, $\BB$ is always a subspace of $\buB$.
The domain parameter $\rho$ and the constant $\Beth$
that appears in the normalization condition \equ(ABNormaliz)
are chosen to have the fixed values
$$
\rho=2^5\,,\qquad\Beth=2^{-12}\,.
\equation(ParamValues)
$$
The domain parameter $\varrho$ is defined implicitly in our proofs.
That is, the lemmas below
hold for $\varrho>1$ sufficiently close to $1$.

Consider first the problem of determining
the bifurcation point $\gamma_0$ and the associated frequency $\alpha_0$.
Let $\BB=\mean_{\sss\{0,1\}}\buB$ over $\real$.
For every $\delta>0$ define $B_\delta=\{h\in\BB: \|h\|<\delta\}$.
Let $s=0$, and denote by $\DD$ the set of all functions $\phi\in\BB$
with the property that $B\hat\phi\ne 0$.
Define $\FF:\DD\to\BB$ to be the map $\phi\mapsto\tilde\phi$ given by \equ(scNSsFix),
with $\gamma=\gamma(\phi)$ and $\alpha=\alpha(\phi)$
defined by the equation \equ(gammaalphaSol).
Clearly, $\FF$ is not only $\rmC^1$ but real analytic on $\DD$.

\claim Lemma(BifPoint)
With $\FF$ as described above,
there exists an affine isomorphism $h\mapsto\varphi+L_1h$ of $\BB$
and real numbers $\eps,\delta,K>0$ satisfying $\eps+K\delta<\delta$,
such that the following holds.
The quasi-Newton map $\NN$ associated with $(\BB,\FF,\varphi,L_1)$
includes the ball $B_\delta$ in its domain
and satisfies the bounds \equ(ContrMappingThm).
The domain of $\FF$ includes the ball in $\BB$
of radius $r=\delta\|L_1\|$, centered at $\varphi$.
For every function $\phi$ in this ball,
$\gamma(\phi)=83.1733117\ldots$ and $\alpha(\phi)=4.66592275\ldots$.

Our proof of this lemma is computer-assisted
and will be described in Section 4.

By \clm(ContrMappingThm),
the map $\FF$ has a unique fixed point $\phi^\ast\in\varphi+L_1 B_\delta$.
We define $\gamma_0=\gamma(\phi^\ast)$ and $\alpha_0=\alpha(\phi^\ast)$.

Our next goal is to construct
a branch of periodic solutions for the equation \equ(scNS).
Consider $\BB=\buB$ over $\field\in\{\real,\complex\}$.
By continuity, there exists an open ball $\JJ_0\subset\field$ centered at the origin,
and an open neighborhood $\DD$ of $\phi^\ast$ in $\BB$, such that
$B\hat\phi=B|\Delta|^{-1/2}\LOP_s(\phi)\phi$ is nonzero for all $s\in\JJ_0$
and all $\phi\in\DD$.
For every $s\in\JJ_0$,
define $\FF_s:\DD\to\BB$ to be the map $\phi\mapsto\tilde\phi$
given by \equ(scNSsFix), with $\gamma=\gamma(s,\phi)$
and $\alpha=\alpha(s,\phi)$ defined by the equation \equ(gammaalphaSol).

\claim Lemma(PerSolutions)
Let $\field=\real$.
There exists a isomorphism $L$ of $\buB$
such that the following holds.
If $\NN_0$ denotes the the quasi-Newton map
associated with $(\DD,\FF_0,\phi^\ast,L)$,
then the derivative $D\NN_0(0)$ of $\NN_0$
at the origin is a contraction.

Our proof of this lemma is computer-assisted
and will be described in Section 4.
As a consequence we have the following.

\claim Corollary(PerBranch)
Consider $\field=\complex$.
There exists an open disk $\JJ\subset\complex$,
centered at the origin,
and an analytic curve $s\mapsto\phi_s$ on $\JJ$ with values in $\DD$,
such that $\FF_s(\phi_s)=\phi_s$ for all $s\in\JJ$.
If $s$ belongs to the real interval $\JJ\cap\real$, then $\phi_s$ is real.
Furthermore, $\phi_0=\phi^\ast$.

\proof
Consider still $\field=\complex$.
For $s\in\II_0$, the derivative of $\NN_s$ on its domain is given by
$$
D\NN_s(h)=D\FF_s(\phi^\ast+Lh)L+\id-L\,.
\equation(DNNsh)
$$
Assume that some function $\psi\in\buB$ satisfies
$D\FF_0(\phi^\ast)\psi=\psi$.
We may assume that $\psi$ takes real values for real arguments.
A straightforward computation shows that
$D\NN_0(0)L^{-1}\psi=L^{-1}\psi$.
Since $D\NN_0(0)$ is a contraction in the real setting,
by \clm(PerSolutions), this implies that $\psi=0$.
So the operator $D\FF_0(\phi^\ast)$ does not have an eigenvalue $1$.
This operator is compact,
since it is the composition of a bounded linear operator
with the compact operator $\LL_\alpha$.
Thus, $D\FF_0(\phi^\ast)$ has no spectrum at $1$.
By the implicit function theorem,
there exists a complex open ball $\JJ$, centered at the origin,
such that the fixed point equation $\FF_s(\phi)=\phi$
has a solution $\phi=\phi_s$ for all $s\in\JJ$.
Furthermore, the curve $s\mapsto\phi_s$ is analytic,
passes through $\phi^\ast$ at $s=0$,
and there is a unique curve with this property.
By uniqueness, we also have $\ov{\phi_{\bar s}}=\phi_s$ for all $s\in\JJ$,
so $\phi_s$ is real for real values of $s\in\JJ$.
\qed

A branch of stationary periodic solutions for \equ(scNS)
is obtained similarly.
Consider $\BB=\mean_{\sss\{0\}}\buB$ over $\field\in\{\real,\complex\}$.
For every $\gamma\in\field$,
define $\FF_\gamma:\BB\to\BB$ to be the map $\phi\mapsto\tilde\phi$
given by \equ(scNSsFix), with $s=\alpha=0$.
Notice that $\phi^\ast_\even$ is a fixed point of $\FF_{\gamma_0}$.

\claim Lemma(StatSolutions)
Let $\field=\real$. There exists an isomorphism $L_0$ of $\BB$
such that the following holds.
If $\NN_{\gamma_0}$ denotes the the quasi-Newton map
associated with $(\BB,\FF_{\gamma_0},\phi^\ast_\even,L_0)$,
then the derivative $D\NN_{\gamma_0}(0)$ of $\NN_{\gamma_0}$
at the origin is a contraction.

Our proof of this lemma is computer-assisted
and will be described in Section 4.
As a consequence we have the following.

\claim Corollary(StatBranch)
Consider $\field=\complex$.
There exists an open disk $\II\subset\complex$,
centered at $\gamma_0$,
and an analytic curve $\gamma\mapsto\phi_\gamma$ on $\II$ with values in $\BB$,
such that $\FF_\gamma(\phi_\gamma)=\phi_\gamma$ for all $\gamma\in\II$.
If $\gamma$ belongs to the real interval $\II\cap\real$, then $\phi_\gamma$ is real.
Furthermore, $\phi_{\gamma_0}=\phi^\ast_\even$.

The proof of this corollary is analogous
to the proof of \clm(PerBranch).

We note that the disk $\II\ni\gamma_0$
is disjoint from the disk $\JJ\ni 0$ described in \clm(PerBranch).
So there is no ambiguity in using the notation
$\gamma\mapsto\phi_\gamma$ and $s\mapsto\phi_s$
for the curve of stationary and periodic solutions,
respectively, of the equation \equ(scNSs),

\smallskip
Based on the results stated in this section, we can now give a

\proofof(NSHopfBif)
As described in the preceding sections,
the curve $\gamma\mapsto\phi_\gamma$ for $\gamma\in\II$
yields a curve $\gamma\mapsto u_\gamma$ of stationary solutions of
the equation \equ(NS), where $u_\gamma=\nsv^{-1}\phi_\gamma$.
By our choice of function spaces,
the function $(\gamma,x,y)\mapsto u_\gamma(x,y)$
is real analytic on $I\times\torus^2$, where $I=\II\cap\real$.

Similarly, the curve $s\mapsto\phi_s$ for $s\in\JJ$
defines a family of of non-stationary periodic solutions
for \equ(NS), with $\gamma=\gamma_s$ and $\alpha=\alpha_s$
determined via the equation \equ(gammaalphaSol).
To be more precise,
the even frequency part $\phi_{s,\even}$ of $\phi_s$
determines a vector field $u_{s,\even}=\nsv^{-1}\phi_{s,\even}$,
and the odd frequency part $\phi_{s,\odd}$
determines a vector field $u_{s,\odd}=\nsv^{-1}\phi_{s,\odd}$.
If $\beta$ is a complex number such that $s=\beta^2\in\JJ$,
then $u=u_{s,\even}+\beta u_{s,\odd}$
is a periodic solution of \equ(NS), with $\gamma=\gamma_s$ and $\alpha=\alpha_s$.
Here, we have used the decomposition \equ(PhiDecomp).
By our choice of function spaces, the functions
$(s,t,x,y)\mapsto u_{s,\even}(t,x,y)$ and $(s,t,x,y)\mapsto u_{s,\odd}(t,x,y)$
are real analytic on $J\times\torus^3$, where $J=\JJ\cap\real$.
Clearly, $\partial_t u_{0,\odd}(t,\bdot\,,\bdot)\ne 0$,
due to the normalization condition $\phi_{-1,1,1}=\Beth$ imposed in \equ(ABNormaliz).
And by construction, we have $u=u_{\gamma_0}$ for $s=0$.
\qed

%%%%%%%%%%%%%%%%%%%%%%%%%%%%
%%%%%%%%%%%%%%%%%%%%%%%%%%%%
\section Remaining estimates
%%%%%%%%%%%%%%%%%%%%%%%%%%%%
%%%%%%%%%%%%%%%%%%%%%%%%%%%%

What remains to be proved are Lemmas
\clmno(BifPoint), \clmno(PerSolutions), and \clmno(StatSolutions).
Our method used in the proof of \clm(BifPoint)
can be considered perturbation theory
about the approximate fixed point $\varphi$ of $\FF$.
The function $\varphi$ is a Fourier polynomial with over $20000$ nonzero coefficients,
so a large number of estimates are involved.

We start by describing bounds on the bilinear function $\LOP$
and on the linear operators $\LL_\alpha$ and $\LL_\alpha'$.
These are the basic building blocks for our transformations
$\FF$, $\FF_s$, and $\FF_\gamma$.
The ``mechanical'' part of these estimates will be described in Subsection 4.4.

%%%%%%%%%%%%%%%%%%%%%%%%%%%%%%%%%%%%%%%%%%%%%%%%%%%%%%%%%%%%%%%%%%%%%%%%%%%
\subsection The bilinear form $\boldmath\LOP$ and a proof of \clm(LOPBound)
%%%%%%%%%%%%%%%%%%%%%%%%%%%%%%%%%%%%%%%%%%%%%%%%%%%%%%%%%%%%%%%%%%%%%%%%%%%

Consider the bilinear form $\LOP$ defined by \equ(LOPDef).
Using the identity \equ(curlInv), we have
$$
\eqalign{
\LOP(\Phi)\phi
&=(\nabla\Phi)\cdot\iso\nabla\Delta^{-1}\phi
+(\nabla\phi)\cdot\iso\nabla\Delta^{-1}\Phi\cr
&=\bigl[(\partial_x\Phi)\Delta^{-1}\partial_y\phi
-(\partial_y\Phi)\Delta^{-1}\partial_x\phi\bigr]
-\bigl[(\Delta^{-1}\partial_x\Phi)\partial_y\phi
-(\Delta^{-1}\partial_y\Phi)\partial_x\phi\bigr]\,.\cr}
\equation(LOPPhiphi)
$$
In order to obtain accurate estimates,
it is useful to have explicit expressions for $\LOP(\Phi)\phi$
in terms of the Fourier coefficients of $\Phi$ and $\phi$.
Given that $\LOP$ is bilinear,
and that the identity \equ(LOPPhiphi) holds pointwise in $t$,
it suffices to compute $\LOP(\Phi)\phi$
for the time-independent monomials
$$
\Phi=\sin_J\times\sin_K\,,\qquad
\phi=\sin_j\times\sin_k\,,
\equation(PhiphiModes)
$$
with $J,K,j,k>0$.
A straightforward computation shows that
$$
\eqalign{
\LOP(\Phi)\phi&=\Theta(Jk+jK)\bigl[
\sin_{J+j}\times\sin_{K-k}-\sin_{J-j}\times\sin_{K+k}
\bigr]\cr
&\quad+\Theta(Jk-jK)\bigl[
\sin_{J+j}\times\sin_{K+k}-\sin_{J-j}\times\sin_{K-k}
\bigr]\,,\cr}
\equation(LOPPhiphiOne)
$$
with $\Theta$ as defined below.
As a result we have
$$
|\Delta|^{-1/2}\LOP(\Phi)\phi
=\sum_{\sigma,\tau=\pm 1}N_{\sigma,\tau}\sin_{\sigma J+j}\times\sin_{\tau K+k}\,,
\equation(InvHalfLapLOPDecomp)
$$
where
$$
N_{\sigma,\tau}
=\Theta
{\sigma Jk-\tau Kj\over\sqrt{(\sigma J+j)^2+(\tau K+k)^2}}\,,\qquad
\Theta=\quarter\biggl({1\over J^2+K^2}-{1\over j^2+k^2}\biggr)\,.
\equation(NpmpmDef)
$$

\proofof(LOPBound)
Using the Cauchy-Schwarz inequality in $\real^2$, we find that
$$
|N_{\sigma,\tau}|
=|\Theta|{|(\sigma J+j)k-(\tau K+k)j|\over\sqrt{(\sigma J+j)^2+(\tau K+k)^2}}
\le|\Theta|\sqrt{j^2+k^2}\,.
\equation(NppBound)
$$
Since the absolute value of $N_{\sigma,\tau}$ is invariant
under an exchange of $(j,k)$ and $(J,K)$, this implies that
$$
|N_{\sigma,\tau}|\le{1/4\over\sqrt{j^2+k^2}}\vee{1/4\over\sqrt{J^2+K^2}}\,,
\equation(NpmpmBound)
$$
where $a\vee b=\max(a,b)$ for $a,b\in\real$.
As a result, we obtain the bound
$$
\bigl\||\Delta|^{-1/2}\LOP(\Phi)\phi\bigr\|
\le\bigl\||\Delta|^{-1/2}\Phi\bigr\|_{\varrho,\epsilon}\|\phi\|
+\|\Phi\|\bigl\||\Delta|^{-1/2}\phi\bigr\|\,.
\equation(LOPPhiphiBound)
$$
Using the nature of the norm \equ(GGrhoepsNorm),
and the fact that $\AA$ is a Banach algebra
for the pointwise product of functions,
this bound extends by bilinearity to arbitrary
functions $\Phi,\phi\in\buB$.
\qed

We note that the bound \equ(LOPPhiphiBound)
exploits the cancellations that lead to the expression \equ(LOPPhiphiOne).
A more straightforward estimate loses a factor of $2$
with respect to \equ(LOPPhiphiBound).
But it is not just this factor of $2$ that counts for us.
The expressions \equ(NpmpmDef) for the coefficients $N_{\sigma,\tau}$
and the bounds \equ(NpmpmBound) are used in our computations
and error estimates.
The expression on the right hand side of \equ(NpmpmBound)
is a decreasing function of the wavenumbers $j,k,J,K$,
so it can be used to estimate $\LOP(\Phi)\phi$
when $\Phi$ and/or $\phi$ are ``tails'' of Fourier series.

%%%%%%%%%%%%%%%%%%%%%%%%%%%%%%%%%%%%%%%%%%%%%%%%%%%%%%%%%%%%%%%%%%%%%%%%%%%%%%%
\subsection The linear operators $\boldmath\LL_\alpha$ and $\boldmath\LL_\alpha'$
%%%%%%%%%%%%%%%%%%%%%%%%%%%%%%%%%%%%%%%%%%%%%%%%%%%%%%%%%%%%%%%%%%%%%%%%%%%%%%%

Consider the linear operators $\LL_\alpha$
and $\LL_\alpha'$ defined in \equ(LLalphaDef), with $\alpha$ real.
A straightforward computation shows that
$$
\psi_{n,j,k}
=\sqrt{j^2+k^2}\,{(j^2+k^2)\phi_{n,j,k}-\alpha n\phi_{-n,j,k}
\over(j^2+k^2)^2+\alpha^2 n^2}\,,\qquad
\psi=\LL_\alpha\phi\,.
\equation(LLalphaphinjk)
$$
Using the Cauchy-Schwarz inequality in $\real^2$,
this yields the estimate
$$
\sqrt{|\psi_{n,j,k}|^2+|\psi_{-n,j,k}|^2}
\le C_{n,j,k}\sqrt{|\phi_{n,j,k}|^2+|\phi_{-n,j,k}|^2}\,,
\equation(LLalphaBound)
$$
with
$$
C_{n,j,k}=\sqrt{j^2+k^2\over(j^2+k^2)^2+\alpha^2 n^2}
\le{1\over\sqrt{2|\alpha n|}}\wedge {1\over\sqrt{j^2+k^2}}
\equation(LLalphaBoundConst)
$$
for $n\ne 0$, where $a\wedge b=\min(a,b)$ for $a,b\in\real$.
The last bound in \equ(LLalphaBoundConst)
is a decreasing function of $|n|,j,k$
and can be used to estimate $\LL_\alpha\phi$
when $\phi$ is the tail of a Fourier series.

For the operator $\LL_\alpha'$ we have
$$
\psi_{n,j,k}
=n\,{(j^2+k^2)\phi_{-n,j,k}+\alpha n\phi_{n,j,k}
\over(j^2+k^2)^2+\alpha^2 n^2}\,,\qquad
\psi=\LL_\alpha'\phi\,.
\equation(LLalphaPrimenjk)
$$
A bound analogous to \equ(LLalphaBound)
holds for $\psi=\LL_\alpha'\phi$, with
$$
C_{n,j,k}=\sqrt{n^2\over(j^2+k^2)^2+\alpha^2 n^2}\,.
\equation(LLalphaPrimeBoundConst)
$$
As can be seen from \equ(DFFsPhiDotPhi),
this bound is needed only for $n=\pm 1$,
since these are the only nonzero frequencies
of the function $\hat\phi=|\Delta|^{-1/2}\LOP_0(\phi)\phi$
with $\phi\in\mean_{\sss\{0,1\}}\buB$.
And for fixed $n$, the right hand side of \equ(LLalphaPrimeBoundConst)
is decreasing in $j$ and $k$.

%%%%%%%%%%%%%%%%%%%%%%%%%%%%%%%%%%%%%
\subsection Estimating operator norms
%%%%%%%%%%%%%%%%%%%%%%%%%%%%%%%%%%%%%

Recall that a function $\phi\in\buB$ admits a Fourier expansion
$$
\phi=\sum_{n\in\oldinteger}\;\sum_{j,k\in\oldnatural_1}\phi_{n,j,k}\theta_{n,j,k}\,,\qquad
\theta_{n,j,k}\defeq\cosi_n\times\sin_j\times\sin_k\,,
\equation(thetanjkDef)
$$
and that the norm of $\phi$ is given by
$$
\|\phi\|=\sum_{j,k\in\oldnatural_1}
\biggl[|\phi_{0,j,k}|+\sum_{n\in\oldnatural_1}
\sqrt{|\phi_{n,j,k}|^2+|\phi_{-n,j,k}|^2}\,\rho^n\biggr]\varrho^{j+k}\,.
\equation(GGrhoepsNormAgain)
$$
Let now $n\ge 0$.
A linear combination $c_{\sss+}\theta_{n,j,k}+c_{\sss-}\theta_{-n,j,k}$
will be referred to as a mode with frequency $n$ and wavenumbers $(j,k)$
or as a mode of type $(n,j,k)$.
We assume of course that $c_{\sss-}=0$ when $n=0$.
Since \equ(GGrhoepsNormAgain) is a weighted $\ell^1$ norm,
except for the $\ell^2$ norm used for modes,
we have a simple expression for the operator norm
of a continuous linear operator $\LL:\buB\to\buB$, namely
$$
\|\LL\|=\sup_{j,k\in\oldnatural_1}\;
\sup_{n\in\oldnatural_0}\;\sup_u\|\LL u\|/\|u\|\,,
\equation(OpNorm)
$$
where the third supremum is over all nonzero modes
$u$ of type $(n,j,k)$.

Let now $n,j,k\ge 1$ be fixed.
In computation where $\LL\theta_{\pm n,j,k}$ is known explicitly,
we use the following estimate.
Denote by $\LL_{n,j,k}$ the restriction of $\LL$
to the subspace spanned by the two functions $\theta_{\pm n,j,k}$.
For $q\ge 1$ define
$$
\|\LL_{n,j,k}\|_q
=\sup_{0\le p<q}\|\LL v_p\|\,,\qquad
v_p=\cos\Bigl({\pi p\over q}\Bigr){\theta_{n,j,k}\over\rho^n\varrho^{j+k}}
+\sin\Bigl({\pi p\over q}\Bigr){\theta_{-n,j,k}\over\rho^n\varrho^{j+k}}\,.
\equation(mNorm)
$$
Since every unit vector in the span of $\theta_{\pm n,j,k}$
lies within a distance less than ${\pi\over q}$
of one of the vectors $v_p$ or its negative, we have
$\|\LL_{n,j,k}\|\le\|\LL_{n,j,k}\|_q+{\pi\over q}\|\LL_{n,j,k}\|$.
Thus
$$
\|\LL_{n,j,k}\|\le{q\over q-\pi}\|\LL_{n,j,k}\|_m\,,\qquad q\ge 4\,.
\equation(LLmodeNormBound)
$$

Consider now the operator
$D\FF_s(\phi)$ described in \equ(DFFsPhiDotPhi),
with $\phi\in\mean_{\sss\{0,1\}}\buB$ fixed.
If $\dot\phi=u_n$ is a nonzero mode with frequency $n\ge 3$,
then $\dot{\hat\phi}=2|\Delta|^{-1/2}\LOP_0(\phi)\dot\phi$
belongs to $\mean_N\buB$ with $N=\{n-1,n,n+1\}$.
Thus, we have $\dot\gamma=\dot\alpha=0$, and
$$
D\FF_0(\phi)u_n
=-\gamma\LL_\alpha|\Delta|^{-1/2}\LOP_0(\phi)u_n\,.
\equation(DFFophiu)
$$
Due to the factor $\LL_\alpha$ in this equation,
if $u_n=c_{\sss+}\theta_{n,j,k}+c_{\sss-}\theta_{-n,j,k}$
with $(j,k)$ and $c_{\sss\pm}$ fixed,
then the ratios
$$
\|D\FF_0(\phi)u_n\|/\|u_n\|
\equation(DFFophiunRatios)
$$
are decreasing in $n$ for $n\ge 3$.
And the limit as $n\to\infty$ of this ratio is zero.

So for the operator $\LL=D\FF_0(\phi)$,
the supremum over $n\in\oldnatural_0$ in \equ(OpNorm)
reduces to a maximum over finitely many terms.
The same holds for the operator $\LL=D\NN_0(0)=D\FF_0(\phi^\ast)L+\id-L$
that is described in \clm(PerSolutions).
This is a consequence of the following choice.

\demo Remark(LMatrix)
The operator $L$ chosen in \clm(PerSolutions)
is a ``matrix perturbation'' of the identity,
in the sense that $L\theta_{n,j,k}=\theta_{n,j,k}$
for all but finitely many indices $(n,j,k)$.
The same is true for the operators $L_1$
and $L_0$ chosen in \clm(BifPoint) and \clm(StatSolutions),
respectively.

%%%%%%%%%%%%%%%%%%%%%%%%%%%%%%
\subsection Computer estimates
%%%%%%%%%%%%%%%%%%%%%%%%%%%%%%

Lemmas \clmno(BifPoint), \clmno(StatSolutions),
and \clmno(PerSolutions) assert the existence of certain objects
that satisfy a set of strict inequalities.
The goal here is to construct these objects,
and to verify the necessary inequalities
by combining the estimates that have been described so far.

The above-mentioned ``objects'' are real numbers,
real Fourier polynomials, and linear operators
that are finite-rank perturbations of the identity.
They are obtained via purely numerical computations.
Verifying the necessary inequalities is largely an organizational task,
once everything else has been set up properly.
Roughly speaking, the procedure follows that of a well-designed numerical program,
but instead of truncation Fourier series and ignoring rounding errors,
we determine rigorous enclosures at every step along the computation.
This part of the proof is written in the programming language Ada [\rAda].
The following is meant to be a rough guide for the reader who
wishes to check the correctness of our programs.
The complete details can be found in [\rProgs].

\smallskip
An enclosure for a function $\phi\in\buB$
is a set in $\buB$ that includes $\phi$
and is defined in terms of (bounds on) a Fourier polynomial
and finitely many error terms.
We define such sets hierarchically,
by first defining enclosures for elements in simpler spaces.
In this context, a ``bound'' on a map $f:\XX\to\YY$
is a function $F$ that assigns to a set $X\subset\XX$
of a given type ({\tt Xtype}) a set $Y\subset\YY$
of a given type ({\tt Ytype}), in such a way that
$y=f(x)$ belongs to $Y$ for all $x\in X$.
In Ada, such a bound $F$ can be implemented by defining
a {\tt procedure F(X:{\gapii}in Xtype; Y:{\gapii}out Ytype)}.

Our most basic enclosures are specified by pairs {\tt S=(S.C,S.R)},
where {\tt S.C} is a representable real number ({\tt Rep})
and {\tt S.R} a nonnegative representable real number ({\tt Radius}).
Given a Banach algebra $\XX$ with unit ${\bf 1}$,
such a pair {\tt S} defines a ball in $\XX$ which we denote by
$\langle{\tt S},\XX\rangle=\{x\in\XX:\|x-({\tt S.C}){\bf 1}\|\le{\tt S.R}\}$.

When $\XX=\real$,
then the data type described above is called {\tt Ball}.
Bounds on some standard functions involving the type {\tt Ball}
are defined in the package {\tt Flts\_Std\_Balls}.
Other basic functions are covered in the packages {\tt Vectors} and {\tt Matrices}.
Bounds of this type have been used in many computer-assisted proofs;
so we focus here on the more problem-specific aspects of our programs.

Consider now the space $\AA$
for a fixed domain radius $\varrho>1$ of type {\tt Radius}.
As mentioned before \dem(LMatrix),
we only need to consider Fourier polynomials in $\AA$.
Our enclosures for such polynomials are defined by
an {\tt array(-I$_{\rmc}${\gapi}..{\gapi}I$_{\rmc}$) of Ball}.
This data type is named {\tt NSPoly},
and the enclosure associated with data {\tt P} of this type is
$$
\langle{\tt P},\AA\rangle
\defeq\sum_{i=-I_{\rmc}}^{I_{\rmc}}
\bigl\langle{\tt P(i)},\real\bigr\rangle\cosi_{\nu(i)}\,,
\equation(NSPolyEnclosure)
$$
where $\nu$ is an increasing index function with the property that $\nu(-i)=-\nu(i)$.
The type {\tt NSPoly} is defined in the package {\tt NSP},
which also implements bounds on some basic operations
for Fourier polynomials in $\AA$.
Among the arguments to {\tt NSP} is a nonnegative integer $n$
(named {\tt NN}).
Our proof of \clm(StatSolutions) and \clm(BifPoint) uses
$I_c=n=0$ and $I_c=n=1$, respectively, and $\nu(i)=i$.
Values $n\ge 2$ are uses when estimating the norm of $\LL u$
for the operator $\LL=D\NN_0(0)$, with $u$ a mode of frequency $n$.
In this case, $\nu$ takes values in $\{-n,n\}$ or $\{-n-1,-n,-n+1,0,n-1,n,n+1\}$,
depending on whether $n$ is odd or even.
(The value $\nu=0$ is being used only for $n=2$.)
The package {\tt NSP} also defines a data type {\tt NSErr}
as an {\tt array(0{\gapi}..{\gapi}I$_{\rm c}$) of Radius}.
This type will be used below.

Given in addition a positive number $\varrho\ge 1$ of type {\tt Radius},
our enclosures for functions in $\buB$
are defined by pairs {\tt(F.C,F.E)},
where {\tt F.C} is an
{\tt array(1{\gapi}..{\gapi}J$_{\rmc}$,1{\gapi}..{\gapi}K$_{\rmc}$) of NSPoly}
and {\tt F.E} is an
{\tt array(1{\gapi}..{\gapi}J$_{\rme}$,1{\gapi}..{\gapi}K$_{\rme}$) of NSErr};
all for a fixed value of the parameter {\tt NN}.
This data type is named {\tt Fourier3},
and the enclosure associated with {\tt F=(F.C,F.E)} is
$$
\langle{\tt F},\buB\rangle
\defeq\sum_{j=1}^{J_\rmc}\sum_{k=1}^{K_\rmc}\bigl\langle{\tt F.C(j,k)},\AA\bigr\rangle
\times\sin_j\times\sin_k
+\sum_{J=1}^{J_\rme}\sum_{K=1}^{K_\rme}H_{\sss J,K}({\tt F.E(J,K)})\,.
\equation(FouThreeEnclosure)
$$
Here, $H_{\sss J,K}({\tt E})$ denotes the set of all functions
$\phi=\sum_{i=0}^{I_{\rmc}}\phi^i$ with $\|\phi^i\|\le{\tt E(i)}$,
where $\phi^i$ can be any function in $\buB$
whose coefficients $\phi^i_{n,j,k}$ vanish unless
$j\ge J$, $k\ge K$, and $|n|=\nu(i)$.

The type {\tt Fourier3} and bounds on some standard functions
involving this type are defined in the child package {\tt NSP.Fouriers}.
This package is a modified version of the
package {\tt Fouriers2} that was used earlier in [\rAKxii,\rAKxviii,\rAGK].
The procedure {\tt Prod} is now a bound on the bilinear map $|\Delta|^{-1/2}\LOP_0$.
The error estimates used in {\tt Prod} are based on the inequality \equ(NpmpmBound).
The package {\tt NSP.Fouriers}
also includes bounds {\tt InvLinear} and {\tt DtInvLinear}
on the linear operators $\LL_\alpha$ and $\LL_\alpha'$, respectively.
These bounds use the estimates described in Subsection 4.3.

As far as the proof of \clm(BifPoint) is concerned,
it suffices now to compose existing bounds
to obtain a bound on the map $\FF$ and its derivative $D\FF$.
This is done by the procedures {\tt GMap} and {\tt DGMap} in {\tt Hopf.Fix}.
Here we use enclosures of type {\tt NN=1}.

The type of quasi-Newton map $\NN$ defined by \equ(quasiNewtonMap)
has been used in several computer-assisted proof before.
So the process of constructing a bound on $\NN$
from a bound on $\FF$ has been automated
in the generic packages {\tt Linear} and {\tt Linear.Contr}.
(Changes compared to earlier versions are mentioned in the program text.)
This includes the computation of an approximate inverse $L_1$
for the operator $\id-D\FF(\varphi)$.
A bound on $\NN$ is defined (in essence) by the procedure {\tt Linear.Contr.Contr},
instantiated with {\tt Map => GMap}.
And a bound on $D\NN$ is defined by {\tt Linear.Contr.Contr},
with {\tt DMap => DGMap}.
Bounds on operator norms are obtained via {\tt Linear.OpNorm}.
Another problem-dependent ingredient in these procedures,
besides {\tt Map} and {\tt DMap}, are data of type {\tt Modes}.
These data are constructed by the procedure {\tt Make}
in the package {\tt Hopf}.
They define a splitting of the given space $\BB$ into a finite direct sum.
For details on how such a splitting is defined and used we refer to [\rAKxix].

If the parameter {\tt NN} has the value $0$,
then the procedures {\tt GMap} and {\tt DGMap} define
bounds on the map $\FF_\gamma$ and its derivative, respectively.
The operator $L_0$ used in \clm(StatSolutions)
has the property that $M_0=L_0-\id$
satisfies $M_0=P_0M_0P_0$, where $P_0=\mean_{\sss\{0\}}\proj_{m_0}$
for some positive integer $m_0$.
Here, and in what follows, $\proj_m$ denotes
the canonical projection in $\buB$
with the property that $\proj_m\phi$ is obtained
from $\phi$ by restricting the second sum in \equ(thetanjkDef)
to wavenumbers $j,k\le m$.

If {\tt NN} has a value $n\ge 2$, then the procedure {\tt DGMap}
defines a bound on the map $(\phi,\psi)\mapsto D\FF_0(\phi)\psi$,
restricted to the subspace $\mean_{\sss\{0,1\}}\buB\times\mean_{\sss\{n\}}\buB$.
The linear operator $L$ that is used in \clm(PerSolutions)
admits a decomposition $L=\id+M_1+M_2+\ldots+M_N$ of the following type.
After choosing a suitable sequence $n\mapsto m_n$ of positive integers,
we set $M_n=P_n(L-\id)P_n$,
where $P_1=\mean_{\sss\{0,1\}}\proj_{m_1}$
and $P_n=\mean_{\sss\{n\}}\proj_{m_n}$ for $n=2,3,\ldots,N$.
This structure of $L$ simplifies the use of \equ(OpNorm)
for estimating the norm of $\LL=D\NN_0(0)$.
Furthermore, to check that $L$ is invertible,
it suffices to verify that $\id+M_n$ is invertible
on the finite-dimensional subspace $P_n\buB$,
for each positive $n\le N$.

The linear operator $L_1$ that is used in \clm(BifPoint)
is of the form $L_1=\id+M_1$ with $M_1$ as described above.

All the steps required in the proofs
of Lemmas \clmno(BifPoint), \clmno(StatSolutions), and \clmno(PerSolutions)
are organized in the main program {\tt Check}.
As $n$ ranges from $0$ to $N=305$,
this program defines the parameters that are used in the proof
for {\tt NN} $=n$, instantiates the necessary packages,
computes the appropriate matrix $M_n$,
verifies that $\id+M_n$ is invertible,
reads $\varphi$ from the file {\tt BP.approx},
and then calls the procedure {\tt ContrFix} from the
(instantiated version of the) package {\tt Hopf.Fix} to verify the necessary inequalities.

The representable numbers ({\tt Rep}) used in our programs
are standard [\rIEEE] extended floating-point numbers (type {\tt LLFloat}).
High precision [\rMPFR] floating-point numbers (type {\tt MPFloat})
are used as well, but not in any essential way.
Both types support controlled rounding.
{\tt Radius} is always a subtype of {\tt LLFloat}.
Our programs were run successfully on a $20$-core workstation,
using a public version of the gcc/gnat compiler [\rGnat].
For further details,
including instruction on how to compile and run our programs,
we refer to [\rProgs].

\bigskip
%%%%%%%%%%%
\references
%%%%%%%%%%%

{\ninepoint

\item{[\rHopf]} E.~Hopf,
{\sl Abzweigung einer periodischen L\"osung
von einer station\"aren L\"osung eines Differentialsystems},
Ber. Math.-Phys. Kl. Siichs. Akad. Wiss. Leipzig, {\bf 94}, 3--22 (1942).

\item{[\rSerrin]} J.~Serrin,
{\sl A Note on the Existence of Periodic Solutions of the Navier-Stokes Equations},
Arch. Rational Mech. Anal. {\bf 3} 120--122, (1959).
%%% periodic forcing

\item{[\rCrAb]} M.G.~Crandall and P.H.~Rabinowitz,
{\sl The Hopf bifurcation theorem in infinite dimensions},
Arch. Rational Mech. Anal. {\bf 67}, 53--72 (1977).
%\pdfclink{0 0 1}{online here}
%{https://www.semanticscholar.org/paper/The-Hopf-Bifurcation-Theorem-in-infinite-dimensions-Crandall-Rabinowitz/c21e703fdf69882c64cff84eb976cc16ddf7345e}

\item{[\rMarMcC]} J.~Marsden, M.~McCracken,
{\sl The Hopf bifurcation and its applications},
Springer Applied Mathematical Sciences Lecture Notes Series, Vol. 19, 1976.

\item{[\rRT]} D.~Ruelle, F.~Takens, {\sl On the Nature of Turbulence},
Commun. Math. Phys. 20, 167--192 (1971)

\item{[\rKlWe]} P.~Kloeden, R.~Wells,
{\sl An explicit example of Hopf bifurcation in fluid mechanics},
Proc. Roy. Soc. London Ser. A {\bf 390}, 293--320 (1983).
%\pdfclink{0 0 1}{online here}
%{https://royalsocietypublishing.org/doi/10.1098/rspa.1983.0133}

\item{[\rChIo]} P.~Chossat, G.~Iooss,
{\sl Primary and secondary bifurcations in the Couette-Taylor problem},
Japan J. Appl. Math. {\bf 2}, 37--68 (1985).

\item{[\rDum]} F.~Dumortier,
{\sl Techniques in the theory of local bifurcations:
blow-up, normal forms, nilpotent bifurcations, singular perturbations};
in: {\sl Bifurcations and periodic orbits of vector fields},
(D.~Schlomiuk, ed., Kluwer Acad.~Pub.)
NATO ASI Ser. C Math. Phys. Sci. {\bf 408}, 10--73 (1993).

\item{[\rChIoo]} P.~Chossat, G.~Iooss,
{\sl The Couette-Taylor problem},
Applied Mathematical Sciences, 102. Springer-Verlag, New York, 1994

\item{[\rNWYNK]} M.T.~Nakao, Y.~Watanabe, N.~Yamamoto, T.~Nishida, M.-N.~Kim,
{\sl Computer assisted proofs of bifurcating solutions for nonlinear heat convection problems},
J. Sci. Comput. {\bf 43}, 388--401 (2010).
%\pdfclink{0 0 1}{online here}{https://doi.org/10.1007/s10915-009-9303-3}
%% Oberbeck-Boussinesq equation (reduction 3d --> 2d) for The Rayleigh-Benard problem.
%% (extension to 3d is discussed at the end)
%% Consider only bifurcations of stationary solutions
%% Similarly for the referencces 9 and 12 that they cite (on the same problem)

\item{[\rAKxii]} G.~Arioli, H.~Koch,
{\sl Non-symmetric low-index solutions for a symmetric boundary value problem},
J. Differ. Equations {\bf 252}, 448--458 (2012).

\item{[\rAKxv]} G.~Arioli, H.~Koch,
{\sl Some symmetric boundary value problems and non-symmetric solutions},
J. Differ. Equations {\bf 259}, 796--816 (2015).

\item{[\rGaldi]} G.P.~Galdi,
{\sl On bifurcating time-periodic flow of
a Navier-Stokes liquid past a cylinder},
Arch. Rational Mech. Anal. {\bf 222}, 285--315 (2016).
Digital Object Identifier (DOI) 10.1007/s00205-016-1001-3

\item{[\rHJNS]} C.-H.~Hsia, C.-Y.~Jung, T.B.~Nguyen, and M.-C.~Shiu,
{\sl On time periodic solutions,
asymptotic stability and bifurcations of Navier-Stokes equations},
Numer. Math. {\bf 135}, 607--638 (2017).
%%% time-periodic forcing, and reference
%%% found this ref in [\rBBLV]

\item{[\rAKxviii]} G.~Arioli, H.~Koch,
{\sl Spectral stability for the wave equation with periodic forcing},
J. Differ. Equations {\bf 265}, 2470--2501 (2018).

\item{[\rAKxix]} G.~Arioli, H.~Koch,
{\sl Non-radial solutions for some semilinear elliptic equations on the disk},
Nonlinear Analysis {\bf 179}, 294–308 (2019).

\item{[\rNPW]} M.T.~Nakao, M.~Plum, Y.~Watanabe,
{\sl Numerical verification methods and computer-assisted proofs
for partial differential equations},
Springer Series in Computational Mathematics, Vol. 53,
Springer Singapore, 2019

\item{[\rGoSe]} J.~G\'omez-Serrano,
{\sl Computer-assisted proofs in PDE: a survey},
SeMA {\bf 76}, 459--484 (2019).

\item{[\rWiZg]} D.~Wilczak, P.~Zgliczy\'nski,
{\sl A geometric method for infinite-dimensional chaos:
Symbolic dynamics for the Kuramoto-Sivashinsky PDE on the line},
J. Differ. Equations {\bf 269}, 8509--8548 (2020).

\item{[\rBBLV]} J.~B.~van den Berg, M.~Breden, J.-P.~Lessard, L.~van Veen,
{\sl Spontaneous periodic orbits in the Navier-Stokes flow},
Preprint 2019,
%\pdfclink{0 0 1}{online here}{https://arxiv.org/abs/1902.00384}

\item{[\rAGK]} G.~Arioli, F.~Gazzola, H.~Koch,
{\sl Uniqueness and bifurcation branches for planar steady
Navier-Stokes equations under Navier boundary conditions},
Preprint 2020.

\item{[\rBLQ]} J.~B.~van den Berg, J.-P.~Lessard, E.~Queirolo,
{\sl Rigorous verification of Hopf bifurcations
via desingularization and continuation},
Preprint 2020.

\item{[\rBQ]} J.~B.~van den Berg, E.~Queirolo,
{\sl Validating Hopf bifurcation in the Kuramoto-Sivashinky PDE},
in preparation.
% talk by Elena Queirolo, see
% https://researchseminars.org/seminar/CRM-CAMP

\item{[\rProgs]} G.~Arioli, H.~Koch,
{\sl Programs and data files for the proof of Lemmas \clmno(BifPoint),
\clmno(StatSolutions), \clmno(PerSolutions), and \clmno(JoinBranches)},
\pdfclink{0 0 1}{{\tt https://web.ma.utexas.edu/users/koch/papers/nshopf/}}
{https://web.ma.utexas.edu/users/koch/papers/nshopf/}

%%%%%%%%%%%%%%%%%%%%%%%%%%%%%%%%%%%%

\item{[\rAda]} Ada Reference Manual, ISO/IEC 8652:2012(E),
available e.g. at\hfil\break
\pdfclink{0 0 1}{{\tt www.ada-auth.org/arm.html}}
{http://www.ada-auth.org/arm.html}

\item{[\rGnat]}
A free-software compiler for the Ada programming language,
which is part of the GNU Compiler Collection; see
\pdfclink{0 0 1}{{\tt gnu.org/software/gnat/}}{http://gnu.org/software/gnat/}

\item{[\rIEEE]} The Institute of Electrical and Electronics Engineers, Inc.,
{\sl IEEE Standard for Binary Float\-ing--Point Arithmetic},
ANSI/IEEE Std 754--2008.

\item{[\rMPFR]} The MPFR library for multiple-precision floating-point computations
with correct rounding; see
\pdfclink{0 0 1}{{\tt www.mpfr.org/}}{http://www.mpfr.org/}

}

\bye
%%%%%%%%%%%%%%%%%%%%%%%%%%%%%%% including param.2 %%%%%%%%%%%%%%%%%%%%%%%%%%%%%%%
%param.2
\magnification=\magstep1
\def\firstpage{1}
\pageno=\firstpage
%%%%%%%%%%%%%%%%%%%%%%%%%%%%%%%%% end of param.2 %%%%%%%%%%%%%%%%%%%%%%%%%%%%%%%%
%%%%%%%%%%%%%%%%%%%%%%%%%%%%%%% including fonts.6 %%%%%%%%%%%%%%%%%%%%%%%%%%%%%%%
%fonts.6
\font\fiverm=cmr5
\font\sevenrm=cmr7
\font\sevenbf=cmbx7
\font\eightrm=cmr8
\font\eightbf=cmbx8
\font\ninerm=cmr9
\font\ninebf=cmbx9
\font\tenbf=cmbx10
\font\magtenbf=cmbx10 scaled\magstep1

\font\magnineeufm=eufm9 scaled\magstep1

%
%%%%%%%%%%%%%%%%%%%%%%%%%%%%%%%%% end of fonts.6 %%%%%%%%%%%%%%%%%%%%%%%%%%%%%%%%
%%%%%%%%%%%%%%%%%%%%%%%%%%%%%%% including smallfonts.tex %%%%%%%%%%%%%%%%%%%%%%%%%%%%%%%
%smallfonts.tex
%
\newskip\ttglue
\font\fiverm=cmr5
\font\fivei=cmmi5
\font\fivesy=cmsy5
\font\fivebf=cmbx5
\font\sixrm=cmr6
\font\sixi=cmmi6
\font\sixsy=cmsy6
\font\sixbf=cmbx6
\font\sevenrm=cmr7
\font\eightrm=cmr8
\font\eighti=cmmi8
\font\eightsy=cmsy8
\font\eightit=cmti8
\font\eightsl=cmsl8
\font\eighttt=cmtt8
\font\eightbf=cmbx8
\font\ninerm=cmr9
\font\ninei=cmmi9
\font\ninesy=cmsy9
\font\nineit=cmti9
\font\ninesl=cmsl9
\font\ninett=cmtt9
\font\ninebf=cmbx9
\font\twelverm=cmr12
\font\twelvei=cmmi12
\font\twelvesy=cmsy12
\font\twelveit=cmti12
\font\twelvesl=cmsl12
\font\twelvett=cmtt12
\font\twelvebf=cmbx12

%% EIGHT POINT FONT FAMILY

\def\eightpoint{\def\rm{\fam0\eightrm}  
  \textfont0=\eightrm \scriptfont0=\sixrm \scriptscriptfont0=\fiverm
  \textfont1=\eighti  \scriptfont1=\sixi  \scriptscriptfont1=\fivei
  \textfont2=\eightsy  \scriptfont2=\sixsy  \scriptscriptfont2=\fivesy
  \textfont3=\tenex  \scriptfont3=\tenex  \scriptscriptfont3=\tenex
  \textfont\itfam=\eightit  \def\it{\fam\itfam\eightit}
  \textfont\slfam=\eightsl  \def\sl{\fam\slfam\eightsl}
  \textfont\ttfam=\eighttt  \def\tt{\fam\ttfam\eighttt}
  \textfont\bffam=\eightbf  \scriptfont\bffam=\sixbf
    \scriptscriptfont\bffam=\fivebf  \def\bf{\fam\bffam\eightbf}
  \tt  \ttglue=.5em plus.25em minus.15em
  \normalbaselineskip=9pt
  \setbox\strutbox=\hbox{\vrule height7pt depth2pt width0pt}
  \let\sc=\sixrm  \let\big=\eightbig \normalbaselines\rm}

\def\eightbig#1{{\hbox{$\textfont0=\ninerm\textfont2=\ninesy
        \left#1\vbox to6.5pt{}\right.$}}}

%% NINE POINT FONT FAMILY

\def\ninepoint{\def\rm{\fam0\ninerm}  
  \textfont0=\ninerm \scriptfont0=\sixrm \scriptscriptfont0=\fiverm
  \textfont1=\ninei  \scriptfont1=\sixi  \scriptscriptfont1=\fivei
  \textfont2=\ninesy  \scriptfont2=\sixsy  \scriptscriptfont2=\fivesy
  \textfont3=\tenex  \scriptfont3=\tenex  \scriptscriptfont3=\tenex
  \textfont\itfam=\nineit  \def\it{\fam\itfam\nineit}
  \textfont\slfam=\ninesl  \def\sl{\fam\slfam\ninesl}
  \textfont\ttfam=\ninett  \def\tt{\fam\ttfam\ninett}
  \textfont\bffam=\ninebf  \scriptfont\bffam=\sixbf
    \scriptscriptfont\bffam=\fivebf  \def\bf{\fam\bffam\ninebf}
  \tt  \ttglue=.5em plus.25em minus.15em
  \normalbaselineskip=11pt
  \setbox\strutbox=\hbox{\vrule height8pt depth3pt width0pt}
  \let\sc=\sevenrm  \let\big=\ninebig \normalbaselines\rm}

\def\ninebig#1{{\hbox{$\textfont0=\tenrm\textfont2=\tensy
        \left#1\vbox to7.25pt{}\right.$}}}

%% TWELVE POINT FONT FAMILY --- not really small

\def\twelvepoint{\def\rm{\fam0\twelverm}  
  \textfont0=\twelverm \scriptfont0=\eightrm \scriptscriptfont0=\sixrm
  \textfont1=\twelvei  \scriptfont1=\eighti  \scriptscriptfont1=\sixi
  \textfont2=\twelvesy  \scriptfont2=\eightsy  \scriptscriptfont2=\sixsy
  \textfont3=\tenex  \scriptfont3=\tenex  \scriptscriptfont3=\tenex
  \textfont\itfam=\twelveit  \def\it{\fam\itfam\twelveit}
  \textfont\slfam=\twelvesl  \def\sl{\fam\slfam\twelvesl}
  \textfont\ttfam=\twelvett  \def\tt{\fam\ttfam\twelvett}
  \textfont\bffam=\twelvebf  \scriptfont\bffam=\eightbf
    \scriptscriptfont\bffam=\sixbf  \def\bf{\fam\bffam\twelvebf}
  \tt  \ttglue=.5em plus.25em minus.15em
  \normalbaselineskip=11pt
  \setbox\strutbox=\hbox{\vrule height8pt depth3pt width0pt}
  \let\sc=\sevenrm  \let\big=\twelvebig \normalbaselines\rm}

\def\twelvebig#1{{\hbox{$\textfont0=\tenrm\textfont2=\tensy
        \left#1\vbox to7.25pt{}\right.$}}}
\catcode`\@=11
%
%  Include all definitions related to the fonts msam, msbm and eufm, so that
%  when this file is used by itself, the results with respect to those fonts
%  are equivalent to what they would have been using AMS-TeX.
%  Most symbols in fonts msam and msbm are defined using \newsymbol;
%  however, a few symbols that replace composites defined in plain must be
%  defined with \mathchardef.

\def\undefine#1{\let#1\undefined}
\def\newsymbol#1#2#3#4#5{\let\next@\relax
 \ifnum#2=\@ne\let\next@\msafam@\else
 \ifnum#2=\tw@\let\next@\msbfam@\fi\fi
 \mathchardef#1="#3\next@#4#5}
\def\mathhexbox@#1#2#3{\relax
 \ifmmode\mathpalette{}{\m@th\mathchar"#1#2#3}%
 \else\leavevmode\hbox{$\m@th\mathchar"#1#2#3$}\fi}
\def\hexnumber@#1{\ifcase#1 0\or 1\or 2\or 3\or 4\or 5\or 6\or 7\or 8\or
 9\or A\or B\or C\or D\or E\or F\fi}

\font\tenmsa=msam10
\font\sevenmsa=msam7
\font\fivemsa=msam5
\newfam\msafam
\textfont\msafam=\tenmsa
\scriptfont\msafam=\sevenmsa
\scriptscriptfont\msafam=\fivemsa
\edef\msafam@{\hexnumber@\msafam}
\mathchardef\dabar@"0\msafam@39
\def\dashrightarrow{\mathrel{\dabar@\dabar@\mathchar"0\msafam@4B}}
\def\dashleftarrow{\mathrel{\mathchar"0\msafam@4C\dabar@\dabar@}}

\def\ulcorner{\delimiter"4\msafam@70\msafam@70 }
\def\urcorner{\delimiter"5\msafam@71\msafam@71 }
\def\llcorner{\delimiter"4\msafam@78\msafam@78 }
\def\lrcorner{\delimiter"5\msafam@79\msafam@79 }
%    Note that there should not be a final space after the digits for a
%    \mathhexbox@.
\def\yen{{\mathhexbox@\msafam@55}}
\def\checkmark{{\mathhexbox@\msafam@58}}
\def\circledR{{\mathhexbox@\msafam@72}}
\def\maltese{{\mathhexbox@\msafam@7A}}

\font\tenmsb=msbm10
\font\sevenmsb=msbm7
\font\fivemsb=msbm5
\newfam\msbfam
\textfont\msbfam=\tenmsb
\scriptfont\msbfam=\sevenmsb
\scriptscriptfont\msbfam=\fivemsb
\edef\msbfam@{\hexnumber@\msbfam}
\def\Bbb#1{{\fam\msbfam\relax#1}}
\def\widehat#1{\setbox\z@\hbox{$\m@th#1$}%
 \ifdim\wd\z@>\tw@ em\mathaccent"0\msbfam@5B{#1}%
 \else\mathaccent"0362{#1}\fi}
\def\widetilde#1{\setbox\z@\hbox{$\m@th#1$}%
 \ifdim\wd\z@>\tw@ em\mathaccent"0\msbfam@5D{#1}%
 \else\mathaccent"0365{#1}\fi}
\font\teneufm=eufm10
\font\seveneufm=eufm7
\font\fiveeufm=eufm5
\newfam\eufmfam
\textfont\eufmfam=\teneufm
\scriptfont\eufmfam=\seveneufm
\scriptscriptfont\eufmfam=\fiveeufm

\catcode`\@=11
%%  Load amssym.def if necessary: If \newsymbol is undefined, do nothing
%%  and the following \input statement will be executed; otherwise
%%  change \input to a temporary no-op.
%#\ifx\undefined\newsymbol \else \begingroup\def\input#1 {\endgroup}\fi
%#\input amssym.def \relax
%%  Most symbols in fonts msam and msbm are defined using \newsymbol.  A few
%%  that are delimiters or otherwise require special treatment have already
%%  been defined as soon as the fonts were loaded.  Finally, a few symbols
%%  that replace composites defined in plain must be undefined first.
\newsymbol\boxdot 1200
\newsymbol\boxplus 1201
\newsymbol\boxtimes 1202
\newsymbol\square 1003
\newsymbol\blacksquare 1004
\newsymbol\centerdot 1205
\newsymbol\lozenge 1006
\newsymbol\blacklozenge 1007
\newsymbol\circlearrowright 1308
\newsymbol\circlearrowleft 1309
\undefine\rightleftharpoons
\newsymbol\rightleftharpoons 130A
\newsymbol\leftrightharpoons 130B
\newsymbol\boxminus 120C
\newsymbol\Vdash 130D
\newsymbol\Vvdash 130E
\newsymbol\vDash 130F
\newsymbol\twoheadrightarrow 1310
\newsymbol\twoheadleftarrow 1311
\newsymbol\leftleftarrows 1312
\newsymbol\rightrightarrows 1313
\newsymbol\upuparrows 1314
\newsymbol\downdownarrows 1315
\newsymbol\upharpoonright 1316
 
\newsymbol\downharpoonright 1317
\newsymbol\upharpoonleft 1318
\newsymbol\downharpoonleft 1319
\newsymbol\rightarrowtail 131A
\newsymbol\leftarrowtail 131B
\newsymbol\leftrightarrows 131C
\newsymbol\rightleftarrows 131D
\newsymbol\Lsh 131E
\newsymbol\Rsh 131F
\newsymbol\rightsquigarrow 1320
\newsymbol\leftrightsquigarrow 1321
\newsymbol\looparrowleft 1322
\newsymbol\looparrowright 1323
\newsymbol\circeq 1324
\newsymbol\succsim 1325
\newsymbol\gtrsim 1326
\newsymbol\gtrapprox 1327
\newsymbol\multimap 1328
\newsymbol\therefore 1329
\newsymbol\because 132A
\newsymbol\doteqdot 132B
 
\newsymbol\triangleq 132C
\newsymbol\precsim 132D
\newsymbol\lesssim 132E
\newsymbol\lessapprox 132F
\newsymbol\eqslantless 1330
\newsymbol\eqslantgtr 1331
\newsymbol\curlyeqprec 1332
\newsymbol\curlyeqsucc 1333
\newsymbol\preccurlyeq 1334
\newsymbol\leqq 1335
\newsymbol\leqslant 1336
\newsymbol\lessgtr 1337
\newsymbol\backprime 1038
\newsymbol\risingdotseq 133A
\newsymbol\fallingdotseq 133B
\newsymbol\succcurlyeq 133C
\newsymbol\geqq 133D
\newsymbol\geqslant 133E
\newsymbol\gtrless 133F
\newsymbol\sqsubset 1340
\newsymbol\sqsupset 1341
\newsymbol\vartriangleright 1342
\newsymbol\vartriangleleft 1343
\newsymbol\trianglerighteq 1344
\newsymbol\trianglelefteq 1345
\newsymbol\bigstar 1046
\newsymbol\between 1347
\newsymbol\blacktriangledown 1048
\newsymbol\blacktriangleright 1349
\newsymbol\blacktriangleleft 134A
\newsymbol\vartriangle 134D
\newsymbol\blacktriangle 104E
\newsymbol\triangledown 104F
\newsymbol\eqcirc 1350
\newsymbol\lesseqgtr 1351
\newsymbol\gtreqless 1352
\newsymbol\lesseqqgtr 1353
\newsymbol\gtreqqless 1354
\newsymbol\Rrightarrow 1356
\newsymbol\Lleftarrow 1357
\newsymbol\veebar 1259
\newsymbol\barwedge 125A
\newsymbol\doublebarwedge 125B
\undefine\angle
\newsymbol\angle 105C
\newsymbol\measuredangle 105D
\newsymbol\sphericalangle 105E
\newsymbol\varpropto 135F
\newsymbol\smallsmile 1360
\newsymbol\smallfrown 1361
\newsymbol\Subset 1362
\newsymbol\Supset 1363
\newsymbol\Cup 1264
 
\newsymbol\Cap 1265
 
\newsymbol\curlywedge 1266
\newsymbol\curlyvee 1267
\newsymbol\leftthreetimes 1268
\newsymbol\rightthreetimes 1269
\newsymbol\subseteqq 136A
\newsymbol\supseteqq 136B
\newsymbol\bumpeq 136C
\newsymbol\Bumpeq 136D
\newsymbol\lll 136E
 
\newsymbol\ggg 136F
 
\newsymbol\circledS 1073
\newsymbol\pitchfork 1374
\newsymbol\dotplus 1275
\newsymbol\backsim 1376
\newsymbol\backsimeq 1377
\newsymbol\complement 107B
\newsymbol\intercal 127C
\newsymbol\circledcirc 127D
\newsymbol\circledast 127E
\newsymbol\circleddash 127F
\newsymbol\lvertneqq 2300
\newsymbol\gvertneqq 2301
\newsymbol\nleq 2302
\newsymbol\ngeq 2303
\newsymbol\nless 2304
\newsymbol\ngtr 2305
\newsymbol\nprec 2306
\newsymbol\nsucc 2307
\newsymbol\lneqq 2308
\newsymbol\gneqq 2309
\newsymbol\nleqslant 230A
\newsymbol\ngeqslant 230B
\newsymbol\lneq 230C
\newsymbol\gneq 230D
\newsymbol\npreceq 230E
\newsymbol\nsucceq 230F
\newsymbol\precnsim 2310
\newsymbol\succnsim 2311
\newsymbol\lnsim 2312
\newsymbol\gnsim 2313
\newsymbol\nleqq 2314
\newsymbol\ngeqq 2315
\newsymbol\precneqq 2316
\newsymbol\succneqq 2317
\newsymbol\precnapprox 2318
\newsymbol\succnapprox 2319
\newsymbol\lnapprox 231A
\newsymbol\gnapprox 231B
\newsymbol\nsim 231C
\newsymbol\ncong 231D
\newsymbol\diagup 201E
\newsymbol\diagdown 201F
\newsymbol\varsubsetneq 2320
\newsymbol\varsupsetneq 2321
\newsymbol\nsubseteqq 2322
\newsymbol\nsupseteqq 2323
\newsymbol\subsetneqq 2324
\newsymbol\supsetneqq 2325
\newsymbol\varsubsetneqq 2326
\newsymbol\varsupsetneqq 2327
\newsymbol\subsetneq 2328
\newsymbol\supsetneq 2329
\newsymbol\nsubseteq 232A
\newsymbol\nsupseteq 232B
\newsymbol\nparallel 232C
\newsymbol\nmid 232D
\newsymbol\nshortmid 232E
\newsymbol\nshortparallel 232F
\newsymbol\nvdash 2330
\newsymbol\nVdash 2331
\newsymbol\nvDash 2332
\newsymbol\nVDash 2333
\newsymbol\ntrianglerighteq 2334
\newsymbol\ntrianglelefteq 2335
\newsymbol\ntriangleleft 2336
\newsymbol\ntriangleright 2337
\newsymbol\nleftarrow 2338
\newsymbol\nrightarrow 2339
\newsymbol\nLeftarrow 233A
\newsymbol\nRightarrow 233B
\newsymbol\nLeftrightarrow 233C
\newsymbol\nleftrightarrow 233D
\newsymbol\divideontimes 223E
\newsymbol\varnothing 203F
\newsymbol\nexists 2040
\newsymbol\Finv 2060
\newsymbol\Game 2061
\newsymbol\mho 2066
\newsymbol\eth 2067
\newsymbol\eqsim 2368
\newsymbol\beth 2069
\newsymbol\gimel 206A
\newsymbol\daleth 206B
\newsymbol\lessdot 236C
\newsymbol\gtrdot 236D
\newsymbol\ltimes 226E
\newsymbol\rtimes 226F
\newsymbol\shortmid 2370
\newsymbol\shortparallel 2371
\newsymbol\smallsetminus 2272
\newsymbol\thicksim 2373
\newsymbol\thickapprox 2374
\newsymbol\approxeq 2375
\newsymbol\succapprox 2376
\newsymbol\precapprox 2377
\newsymbol\curvearrowleft 2378
\newsymbol\curvearrowright 2379
\newsymbol\digamma 207A
\newsymbol\varkappa 207B
\newsymbol\Bbbk 207C
\newsymbol\hslash 207D
\undefine\hbar
\newsymbol\hbar 207E
\newsymbol\backepsilon 237F
%  Restore the catcode value for @ that was previously saved.
%#\catcode`\@=\csname pre amssym.tex at\endcsname

%\endinput
%%%%%%%%%%%%%%%%%%%%%%%%%%%%%%%%% end of symbols.1 %%%%%%%%%%%%%%%%%%%%%%%%%%%%%%%%
%%%%%%%%%%%%%%%%%%%%%%%%%%%%%%% including links.1 %%%%%%%%%%%%%%%%%%%%%%%%%%%%%%%
% links.1
% adapted from http://insti.physics.sunysb.edu/~siegel/tex.shtml
%
% postscript/pdf
\newcount\marknumber	\marknumber=1
\newcount\countdp \newcount\countwd \newcount\countht 
%
% for ordinary tex
%
\ifx\pdfoutput\undefined
\def\rgboo#1{}
\def\postscript#1{\special{" #1}}		%% for dvips
\postscript{
	/bd {bind def} bind def
	/fsd {findfont exch scalefont def} bd
	/sms {setfont moveto show} bd
	/ms {moveto show} bd
	/pdfmark where		% printers ignore pdfmarks
	{pop} {userdict /pdfmark /cleartomark load put} ifelse
	[ /PageMode /UseOutlines		% bookmark window open
	/DOCVIEW pdfmark}
\def\bookmark#1#2{\postscript{		% #1=subheadings (if not 0)
	[ /Dest /MyDest\the\marknumber /View [ /XYZ null null null ] /DEST pdfmark
	[ /Title (#2) /Count #1 /Dest /MyDest\the\marknumber /OUT pdfmark}%
	\advance\marknumber by1}
\def\pdfclink#1#2#3{%
	\hskip-.25em\setbox0=\hbox{#2}%
		\countdp=\dp0 \countwd=\wd0 \countht=\ht0%
		\divide\countdp by65536 \divide\countwd by65536%
			\divide\countht by65536%
		\advance\countdp by1 \advance\countwd by1%
			\advance\countht by1%
		\def\linkdp{\the\countdp} \def\linkwd{\the\countwd}%
			\def\linkht{\the\countht}%
	\postscript{
		[ /Rect [ -1.5 -\linkdp.0 0\linkwd.0 0\linkht.5 ] 
		/Border [ 0 0 0 ]
		/Action << /Subtype /URI /URI (#3) >>
		/Subtype /Link
		/ANN pdfmark}{\rgb{#1}{#2}}}
%
% for pdftex
%
\else
\def\rgboo#1{\pdfliteral{#1 rg #1 RG}}
\pdfcatalog{/PageMode /UseOutlines}		% bookmark window open
\def\bookmark#1#2{
	\pdfdest num \marknumber xyz
	\pdfoutline goto num \marknumber count #1 {#2}
	\advance\marknumber by1}
\def\pdfklink#1#2{%
	\noindent\pdfstartlink user
		{/Subtype /Link
		/Border [ 0 0 0 ]
		/A << /S /URI /URI (#2) >>}{\rgb{1 0 0}{#1}}%
	\pdfendlink}
\fi

\def\rgbo#1#2{\rgboo{#1}#2\rgboo{0 0 0}}
\def\rgb#1#2{\mark{#1}\rgbo{#1}{#2}\mark{0 0 0}}
\def\pdfklink#1#2{\pdfclink{1 0 0}{#1}{#2}}
\def\pdflink#1{\pdfklink{#1}{#1}}
%
% examples:
% \bookmark{0}{look here}
% \pdfclink{0 0 1}{testlink}{http://www.google.com/}
% \pdfklink{testlink}{http://www.google.com/}
% \pdflink{http://www.google.com/}
%%%%%%%%%%%%%%%%%%%%%%%%%%%%%%%%% end of links.1 %%%%%%%%%%%%%%%%%%%%%%%%%%%%%%%%
%%%%%%%%%%%%%%%%%%%%%%%%%%%%%%% including titles.9 %%%%%%%%%%%%%%%%%%%%%%%%%%%%%%%
%titles.8
% requires fonts.5 or higher and smallfonts.tex
% uses links.* if included
% enumerates \demo consecutively (no section number)
%
\newcount\seccount  %% sections
\newcount\subcount  %% subsection
\newcount\clmcount  %% claim
\newcount\equcount  %% equation
\newcount\refcount  %% reference
\newcount\demcount  %% example
\newcount\execount  %% exercise
\newcount\procount  %% problem
\seccount=0
\equcount=1
\clmcount=1
\subcount=1
\refcount=1
\demcount=0
\execount=0
\procount=0
%
%% MISC STUFF
\def\proof{\medskip\noindent{\bf Proof.\ }}
\def\proofof(#1){\medskip\noindent{\bf Proof of \csname c#1\endcsname.\ }}
\def\qed{\hfill{\sevenbf QED}\par\medskip}
\def\references{\bigskip\noindent\hbox{\bf References}\medskip
                \ifx\pdflink\undefined\else\bookmark{0}{References}\fi}
\def\addref#1{\expandafter\xdef\csname r#1\endcsname{\number\refcount}
    \global\advance\refcount by 1}

\def\nextremark #1\par{\item{$\circ$} #1}
\def\firstremark #1\par{\bigskip\noindent{\bf Remarks.}
     \smallskip\nextremark #1\par}
\def\abstract#1\par{{\baselineskip=10pt
    \eightpoint\narrower\noindent{\eightbf Abstract.} #1\par}}
%
%% EQUATION
\def\equtag#1{\expandafter\xdef\csname e#1\endcsname{(\number\seccount.\number\equcount)}
              \global\advance\equcount by 1}
\def\equation(#1){\equtag{#1}\eqno\csname e#1\endcsname}
\def\equ(#1){\hskip-0.03em\csname e#1\endcsname}
%
%% CLAIMS (theorems etc)
\def\clmtag#1#2{\expandafter\xdef\csname cn#2\endcsname{\number\seccount.\number\clmcount}
                \expandafter\xdef\csname c#2\endcsname{#1~\number\seccount.\number\clmcount}
                \global\advance\clmcount by 1}
\def\claim #1(#2) #3\par{\clmtag{#1}{#2}
    \vskip.1in\medbreak\noindent
    {\bf \csname c#2\endcsname .\ }{\sl #3}\par
    \ifdim\lastskip<\medskipamount
    \removelastskip\penalty55\medskip\fi}
\def\clm(#1){\csname c#1\endcsname}
\def\clmno(#1){\csname cn#1\endcsname}
%
%% SECTION
\def\sectag#1{\global\advance\seccount by 1
              \expandafter\xdef\csname sectionname\endcsname{\number\seccount. #1}
              \equcount=1 \clmcount=1 \subcount=1 \execount=0 \procount=0}
\def\section#1\par{\vskip0pt plus.1\vsize\penalty-40
    \vskip0pt plus -.1\vsize\bigskip\bigskip
    \sectag{#1}
    \message{\sectionname}\leftline{\magtenbf\sectionname}
    \nobreak\smallskip\noindent
    \ifx\pdflink\undefined
    \else
      \bookmark{0}{\sectionname}
    \fi}
%
%% SUBSECTION
\def\subtag#1{\expandafter\xdef\csname subsectionname\endcsname{\number\seccount.\number\subcount. #1}
              \global\advance\subcount by 1}
\def\subsection#1\par{\vskip0pt plus.05\vsize\penalty-20
    \vskip0pt plus -.05\vsize\medskip\medskip
    \subtag{#1}
    \message{\subsectionname}\leftline{\tenbf\subsectionname}
    \nobreak\smallskip\noindent
    \ifx\pdflink\undefined
    \else
      \bookmark{0}{.... \subsectionname}  %% can get a bit cluttered
    \fi}
%
%% DEMO (examples etc)
\def\demtag#1#2{\global\advance\demcount by 1
              \expandafter\xdef\csname de#2\endcsname{#1~\number\demcount}}
\def\demo #1(#2) #3\par{
  \demtag{#1}{#2}
  \vskip.1in\medbreak\noindent
  {\bf #1 \number\demcount.\enspace}
  {\rm #3}\par
  \ifdim\lastskip<\medskipamount
  \removelastskip\penalty55\medskip\fi}
\def\dem(#1){\csname de#1\endcsname}
%
%% EXERCISE
\def\exetag#1{\global\advance\execount by 1
              \expandafter\xdef\csname ex#1\endcsname{Exercise~\number\seccount.\number\execount}}
\def\exercise(#1) #2\par{
  \exetag{#1}
  \vskip.1in\medbreak\noindent
  {\bf Exercise \number\execount.}
  {\rm #2}\par
  \ifdim\lastskip<\medskipamount
  \removelastskip\penalty55\medskip\fi}
\def\exe(#1){\csname ex#1\endcsname}
%
%% PROBLEM
\def\protag#1{\global\advance\procount by 1
              \expandafter\xdef\csname pr#1\endcsname{\number\seccount.\number\procount}}
\def\problem(#1) #2\par{
  \ifnum\procount=0
    \parskip=6pt
    \vbox{\bigskip\centerline{\bf Problems \number\seccount}\nobreak\medskip}
  \fi
  \protag{#1}
  \item{\number\procount.} #2}
\def\pro(#1){Problem \csname pr#1\endcsname}
%
%%%%%%%%%%%%%%%%%%%%%%%%%%%%%%%%% end of titles.9 %%%%%%%%%%%%%%%%%%%%%%%%%%%%%%%%
%%%%%%%%%%%%%%%%%%%%%%%%%%%%%%% including macros.21 %%%%%%%%%%%%%%%%%%%%%%%%%%%%%%%
%macros.21
%
% requires fonts.5 or later
% also defines mathds (double strike) family
%
\def\rightheadline{\hfil}
\def\leftheadline{\sevenrm\hfil HANS KOCH\hfil}
\headline={\ifnum\pageno=\firstpage\hfil\else
\ifodd\pageno{{\fiverm\rightheadline}\number\pageno}
\else{\number\pageno\fiverm\leftheadline}\fi\fi}
\footline={\ifnum\pageno=\firstpage\hss\tenrm\folio\hss\else\hss\fi}
\let\ov=\overline
\let\cl=\centerline

\let\eps=\varepsilon
\let\sss=\scriptscriptstyle

\def\AA{{\cal A}}
\def\BB{{\cal B}}

\def\DD{{\cal D}}

\def\FF{{\cal F}}

\def\II{{\cal I}}
\def\JJ{{\cal J}}

\def\LL{{\cal L}}

\def\NN{{\cal N}}

\def\XX{{\cal X}}
\def\YY{{\cal Y}}

\def\rmC{{\rm C}}
\def\rmL{{\rm L}}
\def\id{{\rm I}}

\def\Im{\mathop{\rm Im}\nolimits}
%
%%%%%%%%%%%%%%
\newfam\dsfam
\def\mathds #1{{\fam\dsfam\tends #1}}

\font\tends=dsrom10
\font\eightds=dsrom8
\textfont\dsfam=\tends
\scriptfont\dsfam=\eightds
%%%%%%%%%%%%%%
%

\def\integer{{\mathds Z}}

\def\real{{\mathds R}}
\def\complex{{\mathds C}}

\def\mean{{\Bbb E}}
\def\proj{{\Bbb P}}
\def\torus{{\Bbb T}}
\def\iso{{\Bbb J}}

\def\bdot{\hbox{\bf .}}

\def\defeq{\mathrel{\mathop=^{\sss\rm def}}}
\def\half{{1\over 2}}

\def\quarter{{1\over 4}}
\def\thalf{{\textstyle\half}}

%

%

%

%

%
% from TeX book: used for commutative diagram
% in math mode, before using matrix, do
% \def\normalbaselines{\baselineskip20pt\lineskip3pt\lineskiplimit3pt}

%%%%%%%%%%%%%%%%%%%%%%%%%%%%%%%%% end of macros.21 %%%%%%%%%%%%%%%%%%%%%%%%%%%%%%%%
%%%%%%%%%%%%%%%%%%%%%%%%%%%%%%% including mygraphicx.tex %%%%%%%%%%%%%%%%%%%%%%%%%%%%%%%
%% modification of graphicx.tex by Nathan Goldschmidt
\input miniltx

\ifx\pdfoutput\undefined
  \def\Gin@driver{dvips.def}  % we are not running PDFTeX
\else
  \def\Gin@driver{pdftex.def} % we are running PDFTeX
\fi
 
\input graphicx.sty
\resetatcatcode
%%%%%%%%%%%%%%%%%%%%%%%%%%%%%%%%% end of mygraphicx.tex %%%%%%%%%%%%%%%%%%%%%%%%%%%%%%%%
%%%%%%%%%%%%%%%%%%%%%%%%%%%%%%% including opmac.1 %%%%%%%%%%%%%%%%%%%%%%%%%%%%%%%
%% table macros from opmac.tex
%% Petr Olsak, 2012 -- 2016
%% http://petr.olsak.net/opmac.html

\newcount\tmpnum % auxiliary count
\newdimen\tmpdim % auxiliary dimen
\def\opwarning#1{\immediate\write16{l.\the\inputlineno\space OPmac WARNING: #1.}}
\long\def\addto#1#2{\expandafter\def\expandafter#1\expandafter{#1#2}}
\long\def\isinlist#1#2#3{\begingroup \long\def\tmp##1#2##2\end{\def\tmp{##2}%
   \ifx\tmp\empty \endgroup \csname iffalse\expandafter\endcsname \else
                  \endgroup \csname iftrue\expandafter\endcsname \fi}% end of \def\tmp
   \expandafter\tmp#1\endlistsep#2\end
}
\def\tabstrut{\strut}     % strut in the \table
\def\tabiteml{\enspace}   % left material before each \table item
\def\tabitemr{\enspace}   % right material after each \table item
\def\vvkern{1pt}          % space between vertical lines
\def\hhkern{1pt}          % space between horizontal lines

%%%%%%%%%%%%%% \table -- sec. 3.19 in opmac-d.pdf

\newtoks\tabdata
\def\tabstrutA{\tabstrut}
\newcount\colnum
\def\ddlinedata{}
\def\vvleft{}

\def\table{\vbox\bgroup \catcode`\|=12 \tableA}
\def\tableA#1#2{\offinterlineskip \colnum=0 \def\tmpa{}\tabdata={}\scantabdata#1\relax
   \halign\expandafter{\the\tabdata\cr#2\crcr}\egroup}

\def\scantabdata#1{\let\next=\scantabdata
   \ifx\relax#1\let\next=\relax
   \else\ifx|#1\addtabvrule
      \else\isinlist{123456789}#1\iftrue \def\next{\scantabdataC#1}%
          \else \expandafter\ifx\csname tabdeclare#1\endcsname \relax
                \expandafter\ifx\csname paramtabdeclare#1\endcsname \relax
                   \opwarning{tab-declarator "#1" unknown, ignored}%
                \else \def\next{\expandafter \scantabdataB \csname paramtabdeclare#1\endcsname}\fi
             \else \def\next{\expandafter\scantabdataA \csname tabdeclare#1\endcsname}%
   \fi\fi\fi\fi \next
}
\def\scantabdataA#1{\addtabitem \expandafter\addtabdata\expandafter{#1\tabstrutA}\scantabdata}
\def\scantabdataB#1#2{\addtabitem\expandafter\addtabdata\expandafter{#1{#2}\tabstrutA}\scantabdata}
\def\scantabdataC {\def\tmpb{}\afterassignment\scantabdataD \tmpnum=}
\def\scantabdataD#1{\loop \ifnum\tmpnum>0 \advance\tmpnum by-1 \addto\tmpb{#1}\repeat
   \expandafter\scantabdata\tmpb
}

\def\unsskip{\ifdim\lastskip>0pt \unskip\fi}
\def\addtabitem{\ifnum\colnum>0 \addtabdata{&}\addto\ddlinedata{&\dditem}\fi
    \advance\colnum by1 \let\tmpa=\relax}
\def\addtabdata#1{\tabdata\expandafter{\the\tabdata#1}}
\def\addtabvrule{%
    \ifx\tmpa\vrule \addtabdata{\kern\vvkern}%
       \ifnum\colnum=0 \addto\vvleft{\vvitem}\else\addto\ddlinedata{\vvitem}\fi
    \else \ifnum\colnum=0 \addto\vvleft{\vvitemA}\else\addto\ddlinedata{\vvitemA}\fi\fi
    \let\tmpa=\vrule \addtabdata{\vrule}}

\def\crll{\crcr\noalign{\hrule\kern\hhkern\hrule}}

\def\crli{\crcr \omit
   \gdef\dditem{\omit\tablinefil}\gdef\vvitem{\kern\vvkern\vrule}\gdef\vvitemA{\vrule}%
   \vvleft\tablinefil\ddlinedata\crcr}
\def\crlli{\crli\noalign{\kern\hhkern}\crli}
\def\tablinefil{\leaders\hrule\hfil}

\def\crlp#1{\crcr \noalign{\kern-\drulewidth}%
   \omit \xdef\crlplist{#1}\xdef\crlplist{,\expandafter}\expandafter\crlpA\crlplist,\end,%
   \global\tmpnum=0 \gdef\dditem{\omit\crlpD}%
   \gdef\vvitem{\kern\vvkern\kern\drulewidth}\gdef\vvitemA{\kern\drulewidth}%
   \vvleft\crlpD\ddlinedata \global\tmpnum=0 \crcr}
\def\crlpA#1,{\ifx\end#1\else \crlpB#1-\end,\expandafter\crlpA\fi}
\def\crlpB#1#2-#3,{\ifx\end#3\xdef\crlplist{\crlplist#1#2,}\else\crlpC#1#2-#3,\fi}
\def\crlpC#1-#2-#3,{\tmpnum=#1\relax
   \loop \xdef\crlplist{\crlplist\the\tmpnum,}\ifnum\tmpnum<#2\advance\tmpnum by1 \repeat}
\def\crlpD{\global\advance\tmpnum by1
   \edef\tmpa{\noexpand\isinlist\noexpand\crlplist{,\the\tmpnum,}}%
   \tmpa\iftrue \kern-\drulewidth \tablinefil \kern-\drulewidth\else\hfil \fi}

\def\tskip{\afterassignment\tskipA \tmpdim}
\def\tskipA{\gdef\dditem{}\gdef\vvitem{}\gdef\vvitemA{}\gdef\tabstrutA{}%
    \vbox to\tmpdim{}\ddlinedata \crcr \noalign{\gdef\tabstrutA{\tabstrut}}}

\def\mspan{\omit \tabdata={\tabstrut}\let\tmpa=\relax \afterassignment\mspanA \mscount=}
\def\mspanA[#1]{\loop \ifnum\mscount>1 \csname span\endcsname \omit \advance\mscount by-1 \repeat
   \mspanB#1\relax}
\def\mspanB#1{\ifx\relax#1\def\tmpa{\def\tmpa####1}%
   \expandafter\tmpa\expandafter{\the\tabdata\ignorespaces}\expandafter\tmpa\else
   \ifx |#1\ifx\tmpa\vrule\addtabdata{\kern\vvkern}\fi \addtabdata{\vrule}\let\tmpa=\vrule
   \else \let\tmpa=\relax
      \ifx c#1\addtabdata{\tabiteml\hfil\ignorespaces##1\unsskip\hfil\tabitemr}\fi
      \ifx l#1\addtabdata{\tabiteml\ignorespaces##1\unsskip\hfil\tabitemr}\fi
      \ifx r#1\addtabdata{\tabiteml\hfil\ignorespaces##1\unsskip\tabitemr}\fi
   \fi \expandafter\mspanB \fi}

\newdimen\drulewidth  \drulewidth=0.4pt
\let\orihrule=\hrule  \let\orivrule=\vrule
\def\rulewidth{\afterassignment\rulewidthA \drulewidth}
\def\rulewidthA{\edef\hrule{\orihrule height\the\drulewidth}%
                \edef\vrule{\orivrule width\the\drulewidth}}

\long\def\frame#1{%
   \hbox{\vrule\vtop{\vbox{\hrule\kern\vvkern
      \hbox{\kern\hhkern\relax#1\kern\hhkern}%
   }\kern\vvkern\hrule}\vrule}}
%%%%%%%%%%%%%%%%%%%%%%%%%%%%%%%%% end of opmac.1 %%%%%%%%%%%%%%%%%%%%%%%%%%%%%%%%
%%%%%%%%%%%%%%%%%%%%%%%%%%%%%%% including boldmath.tex %%%%%%%%%%%%%%%%%%%%%%%%%%%%%%%
%% boldmath

%=================================================================
% necessary fonts
%-----------------------------------------------------------------
   \font\Fivebf  =cmbx10  scaled 500 % five point bold
   \font\Sevenbf =cmbx10  scaled 700 % seven point bold
   \font\Tenbf   =cmbx10             % ten point bold
   \font\Fivemb  =cmmib10 scaled 500 % five point math bold
   \font\Sevenmb =cmmib10 scaled 700 % seven point math bold
   \font\Tenmb   =cmmib10            % ten point math bold
%=================================================================
% Math in Bold
% example $\boldmath{..}$   { math characters will be bold }
%-----------------------------------------------------------------
\def\boldmath{\textfont0=\Tenbf           \scriptfont0=\Sevenbf 
              \scriptscriptfont0=\Fivebf  \textfont1=\Tenmb
              \scriptfont1=\Sevenmb       \scriptscriptfont1=\Fivemb}
%=================================================================

%% twelveboldmath

%=================================================================
% necessary fonts
%-----------------------------------------------------------------
%  \font\Sevenbf   =cmbx10  scaled 700

%  \font\Sevenmb   =cmmib10 scaled 700

%=================================================================
% Math in Bold
% example $\twelveboldmath{...}$   { math characters will be bold }
%-----------------------------------------------------------------

%=================================================================

%%%%%%%%%%%%%%%%%%%%%%%%%%%%%%%%% end of boldmath.tex %%%%%%%%%%%%%%%%%%%%%%%%%%%%%%%%
%\input param.2
%\input fonts.6
%\input smallfonts.tex
%\input symbols.1
%\input links.1
%\input titles.9
%\input macros.21
%\input boldmath.tex
%\input mygraphicx.tex
%
\newdimen\savedparindent
\savedparindent=\parindent
\font\tenamsb=msbm10 \font\sevenamsb=msbm7 \font\fiveamsb=msbm5
\newfam\bbfam
\textfont\bbfam=\tenamsb
\scriptfont\bbfam=\sevenamsb
\scriptscriptfont\bbfam=\fiveamsb
\def\field{{\mathds F}}
\def\bbb{\fam\bbfam}
\def\oldinteger{{\bbb Z}}
\def\oldnatural{{\bbb N}}
\def\buB{{\hbox{\magnineeufm B}}}
\let\Beth\theta
\def\gapi{\hskip 2pt}
\def\gapii{\hskip 2pt}
\def\rmc{{\rm c}}
\def\rme{{\rm e}}
\def\id{{\rm I}}
\def\LOP{{\mathds L}}
\def\nsv{\mathop{\not\!\partial}\nolimits}
\def\cosi{\mathop{\rm cosi}\nolimits}
\def\even{{\rm e}}
\def\odd{{\rm o}}
\def\bdot{\hbox{\bf .}}
\def\stwomat#1#2#3#4{{\eightpoint\left[\matrix{#1&#2\cr#3&#4\cr}\right]}}
\addref{Hopf}
\addref{Serrin}
\addref{CrAb}
\addref{MarMcC}
\addref{RT}
\addref{KlWe}
\addref{ChIo}
\addref{Dum}
\addref{ChIoo}
\addref{NWYNK}
\addref{AKxii}
\addref{AKxv}
\addref{Galdi}
\addref{HJNS}
\addref{AKxviii}
\addref{AKxix}
\addref{NPW}
\addref{GoSe}
\addref{WiZg}
\addref{BBLV}
\addref{AGK}
\addref{BLQ}
\addref{BQ}
\addref{Progs}
\addref{Ada}
\addref{Gnat}
\addref{IEEE}
\addref{MPFR}
\def\leftheadline{\sixrm\hfil ARIOLI \& KOCH\hfil}
\def\rightheadline{\sevenrm\hfil Hopf bifurcation in Navier-Stokes\hfil}
%
%%%%%%%%%%%%%%%%%%%%%%%%%%%%%%%%%%%%%%%%%%%%%%%%%%%%%%%
\cl{\magtenbf A Hopf bifurcation in the planar Navier-Stokes equations}
\bigskip

\cl{
Gianni Arioli
\footnote{$^1$}
{\eightpoint\hskip-2.9em
Department of Mathematics, Politecnico di Milano,
Piazza Leonardo da Vinci 32, 20133 Milano.
}
$^{\!\!\!,\!\!}$
\footnote{$^2$}
{\eightpoint\hskip-2.6em
Supported in part by the PRIN project
``Equazioni alle derivate parziali e
disuguaglianze analitico-geometriche associate''.}
and Hans Koch
\footnote{$^3$}
{\eightpoint\hskip-2.7em
Department of Mathematics, The University of Texas at Austin,
Austin, TX 78712.}
}

\bigskip
\abstract
We consider the Navier-Stokes equation
for an incompressible viscous fluid on a square,
satisfying Navier boundary conditions
and being subjected to a time-independent force.
As the kinematic viscosity is varied,
a branch of stationary solutions is shown
to undergo a Hopf bifurcation,
where a periodic cycle branches from the stationary solution.
Our proof is constructive
and uses computer-assisted estimates.

%%%%%%%%%%%%%%%%%%%%%%%%%%%%%%%%%%%%%
%%%%%%%%%%%%%%%%%%%%%%%%%%%%%%%%%%%%%
\section Introduction and main result
%%%%%%%%%%%%%%%%%%%%%%%%%%%%%%%%%%%%%
%%%%%%%%%%%%%%%%%%%%%%%%%%%%%%%%%%%%%

We consider the Navier-Stokes equations
$$
\partial_t u-\nu\Delta u+(u\cdot\nabla)u+\nabla p=f\ ,\quad
\nabla\cdot u=0\quad{\rm on~}\Omega\,,
\equation(NavierStokes)
$$
for the velocity $u=u(t,x,y)$
of an incompressible fluid on a planar domain $\Omega$,
satisfying suitable boundary conditions for $(x,y)\in\partial\Omega$
and initial conditions at $t=0$.
Here, $p$ denotes the pressure,
and $f=f(x,y)$ is a fixed time-independent external force.

Our focus is on solution curves and bifurcations
as the kinematic velocity $\nu$ is being varied.
In order to reduce the complexity of the problem,
the domain $\Omega$ is chosen to be as simple as possible,
namely the square $\Omega=(0,\pi)^2$.
Following [\rAGK], we impose Navier boundary conditions on $\partial\Omega$,
which are given by
$$
\eqalign{
u_1&=\partial_x u_2=0\quad{\rm on~}\{0,\pi\}\times(0,\pi)\,,\cr
u_2&=\partial_y u_1=0\quad{\rm on~}(0,\pi)\times\{0,\pi\}\,.\cr}
\equation(NavierBC)
$$
A fair amount is known about the (non)uniqueness of stationary solutions
in this case [\rAGK].
This includes the existence of a bifurcation between curves of stationary solutions
with different symmetries.

Here we prove the existence of a Hopf bifurcation
for the equation \equ(NavierStokes) with boundary conditions \equ(NavierBC),
and with a forcing function $f$ that satisfies
$$
(\partial_xf_2-\partial_yf_1)(x,y)=5\sin(x)\sin(2y)-13\sin(3x)\sin(2y)\,.
\equation(specificNScurlf)
$$
In a Hopf bifurcation, a stationary solution loses stability
and a small-amplitude limit cycle branches from the stationary solution [\rHopf,\rCrAb,\rMarMcC].
Among other things, this introduces a time scale in the system
and increases its complexity.
In this capacity, Hopf bifurcations in the Navier-Stokes equation
constitute an important first step
in the transition to turbulence in fluids,
as was described in the seminal work [\rRT].

Numerically, there is plenty of evidence that Hopf bifurcations occur
in the Navier-Stokes equation, but proofs are still very scarce.
An explicit example of a Hopf bifurcation
was given in [\rKlWe] for the rotating B\'enard problem.
A proof exists also for the Couette-Taylor problem [\rChIo,\rChIoo].
Sufficient conditions for the existence of a Hopf bifurcation
in a Navier-Stokes setting are presented in [\rGaldi].

Before giving a precise statement of our result, let us replace
the vector field $u$ in the equation \equ(NavierStokes) by $\nu^{-1}u$.
The equation for the rescaled function $u$ is
$$
\alpha\partial_t u-\Delta u+\gamma(u\cdot\nabla)u+\nabla p=f\ ,\quad
\nabla\cdot u=0\quad{\rm on~}\Omega\,,
\equation(NS)
$$
where $\gamma=\nu^{-2}$.
The value of $\alpha$ that corresponds to \equ(NavierStokes) is $\nu^{-1}$,
but this can be changed to any positive value by rescaling time.

Numerically, it is possible to find stationary solutions of \equ(NS)
for a wide range of values of the parameter $\gamma$.
At a value $\gamma_0\approx 83.1733117\ldots$ we observe
a Hopf bifurcation that leads to a branch of periodic solutions for $\gamma>\gamma_0$.

For a fixed value of $\alpha$,
the time period $\tau$ of the solution varies with $\gamma$.
Instead of looking for $\tau$-periodic solution of \equ(NS)
for fixed $\alpha$, we look for $2\pi$-periodic solutions,
where $\alpha=2\pi/\tau$ has to be determined.
To simplify notation, a $2\pi$-periodic function
will be identified with a function on the circle $\torus=\real/(2\pi\integer)$.
Our main theorem is the following.

\claim Theorem(NSHopfBif)
There exists a real number $\gamma_0=83.1733117\ldots$,
an open interval $I$ including $\gamma_0$,
and a real analytic function $(\gamma,x,y)\mapsto u_\gamma(x,y)$
from $I\times\Omega$ to $\real^2$,
such that $u_\gamma$ is a stationary solution of \equ(NS) and \equ(NavierBC)
for each $\gamma\in I$.
In addition, there exists a real number $\alpha_0=4.66592275\ldots$,
an open interval $J$ centered at the origin,
two real analytic functions $\gamma$ and $\alpha$ on $J$
that satisfy $\gamma(0)=\gamma_0$ and $\alpha(0)=\alpha_0$, respectively,
as well as two real analytic functions
$(s,t,x,y)\mapsto u_{s,\even}(t,x,y)$ and $(s,t,x,y)\mapsto u_{s,\odd}(t,x,y)$
from $J\times\torus\times\Omega$
to $\real^2$, such that the following holds.
For any given $\beta\in\complex$ satisfying $\beta^2\in J$,
the vector field $u=u_{s,\even}+\beta u_{s,\odd}$ with $s=\beta^2$
is a solution of \equ(NS) and \equ(NavierBC)
with $\gamma=\gamma(s)$ and $\alpha=\alpha(s)$.
Furthermore, $u_{0,\even}(t,\bdot\,,\bdot)=u_{\gamma_0}$
and $\partial_t u_{0,\odd}(t,\bdot\,,\bdot)\ne 0$.

To our knowledge, this is the first result establishing
the existence of a Hopf bifurcation for the Navier-Stokes equation
in a stationary environment.

Our proof of this theorem is computer-assisted.
The solutions are obtained by rewriting \equ(NS) and \equ(NavierBC)
as a suitable fixed point equation for scalar vorticity of $u$.
Here we take advantage of the fact that the domain is two-dimensional.
We isolate the periodic branch from the stationary branch
by using a scaling that admits two distinct limits at the bifurcation point.
This approach is also known as the blow-up method,
which is a common tool in the study of singularities and bifurcations [\rDum].

\smallskip
Computer-assisted methods have been applied successfully
to many different problems in analysis,
mostly in the areas of dynamical systems and partial differential equations.
Here we will just mention work that
concerns the Navier-Stokes equation or Hopf bifurcations.
For the Navier-Stokes equation,
the existence of symmetry-breaking bifurcations among stationary solutions
has been established in [\rNWYNK,\rAGK].
Periodic solutions for the Navier-Stokes flow in a stationary environment
have been obtained in [\rBBLV].
In the case of periodic forcing,
the problem of existence and stability of periodic orbits
has been investigated in [\rHJNS].
Concerning the existence of Hopf bifurcations,
a computer-assisted proof was given recently in [\rBLQ]
for a finite-dimensional dynamical system;
and an extension of their method to the Kuramoto-Sivashinsky PDE
is presented in [\rBQ].
For other recent computer-assisted proofs
we refer to [\rAKxix,\rNPW,\rGoSe,\rWiZg] and references therein.

Figure 1 depicts snapshots at $t=0$ and $t=\pi$
of a solution $u:\torus\times\Omega\to\real^2$ of the
equations \equ(NS) with boundary conditions \equ(NavierBC)
and forcing \equ(specificNScurlf),
obtained numerically for the parameter value $\gamma\approx 84.00\ldots$.

%%%%%%%%%%%%%%%%%%%%%%%%%%%%%%%%%%%%%%%%%%%%%%%%%%%%%%%%%%%%%
\vskip0.15in
\hbox{\hskip 55pt
\includegraphics[height=1.5in,width=1.5in]{figures/ps1.eps}
\hskip 35pt
\includegraphics[height=1.5in,width=1.5in]{figures/ps2.eps}}
\vskip0.1in
\centerline{\eightpoint {\noindent\bf Figure 1.}
Snapshots at two distinct times of a time-periodic solution
for $\gamma\approx 84.00\ldots$
}
\vskip0.15in
%%%%%%%%%%%%%%%%%%%%%%%%%%%%%%%%%%%%%%%%%%%%%%%%%%%%%%%%%%%%%

As mentioned earlier, a system similar to the one
considered here is known to exhibit a symmetry-breaking bifurcation
within the class of stationary solutions [\rAGK].
The broken symmetry is $y\mapsto\pi/2-y$.
Based on a numerical computation of eigenvalues,
we expect an analogous bifurcation to occur here at $\gamma\approx 1450$.
Interestingly, the Hopf bifurcation described here occurs
at a significantly smaller value of $\gamma$.
We have not tried to prove the existence
of a symmetry-breaking bifurcation for the forcing \equ(specificNScurlf),
since such an analysis would duplicate the work in [\rAGK]
and go beyond the scope of the present paper.

\smallskip
The remaining part of this paper is organized as follows.
In Section 2, we first rewrite \equ(NS) as an equation for
the function $\Phi=\partial_y u_1-\partial_x u_2$,
which is the scalar vorticity of $-u$.
After a suitable scaling $\Phi=U_\beta\phi$,
the problem of constructing the solution branches
described in \clm(NSHopfBif) is reduced to three fixed point problems
for the function $\phi$.
These fixed point equations are solved in Section 3,
based on estimates described in Lemmas 3.3, 3.4, and 3.6.
Section 4 is devoted to the proof of these estimates,
which involves reducing them to a large number of trivial bounds
that can be (and have been) verified with the aid of a computer [\rProgs].

%%%%%%%%%%%%%%%%%%%%%%%%%%%%%%
%%%%%%%%%%%%%%%%%%%%%%%%%%%%%%
\section Fixed point equations
%%%%%%%%%%%%%%%%%%%%%%%%%%%%%%
%%%%%%%%%%%%%%%%%%%%%%%%%%%%%%

The goal here is to rewrite the equation \equ(NS)
with boundary conditions \equ(NavierBC) as a fixed point problem.
Applying the operator $\nsv:(u_1,u_2)\mapsto\partial_2 u_1-\partial_1 u_2$
on both sides of the equation \equ(NS), we obtain
$$
\alpha\partial_t\Phi
-\Delta\Phi+\gamma u\cdot\nabla\Phi=\nsv f\,,\qquad\Phi=\nsv u\,.
\equation(curlNS)
$$
Here, we have used that $\nsv\,(u\cdot\nabla)u=u\cdot\nabla\Phi$.
Using the divergence-free condition $\nabla\cdot u=0$,
one also finds that
$$
\Delta u=\iso\nabla\Phi\,,\qquad
\iso=\stwomat{0}{1}{-1}{0}\,.
\equation(Lapu)
$$
If $\Phi$ vanishes on the boundary of $\partial\Omega$,
then the equation \equ(Lapu) can be inverted to yield
$$
u=\nsv^{-1}\Phi\defeq\iso\nabla\Delta^{-1}\Phi\,,
\equation(curlInv)
$$
where $\Delta$ denotes the Dirichlet Laplacean on $\Omega$.

In Section 3
we will define a space of real analytic functions $\Phi$
that admit a representation
$$
\Phi(t,x,y)=\sum_{j,k\in\oldnatural_1}\Phi_{j,k}(t)\sin(jx)\sin(ky)\,,
\equation(PhixyExpansion)
$$
with the series converging uniformly on a complex
open neighborhood of $\torus^3$.
Here, and in what follows, $\oldnatural_1$ denotes the set of all positive integers.
If $\Phi$ admits such an expansion,
then the equation \equ(curlInv) yields
$$
\eqalign{
u_1(t,x,y)&=\sum_{j,k\in\oldnatural_1}\,{-k\over j^2+k^2}\Phi_{j,k}(t)\sin(jx)\cos(ky)\,,\cr
u_2(t,x,y)&=\sum_{j,k\in\oldnatural_1}\,{j\over j^2+k^2}\Phi_{j,k}(t)\cos(jx)\sin(ky)\,.\cr}
\equation(uxyExpansion)
$$
It is straightforward to check that the corresponding
vector field $u=(u_1,u_1)$ satisfies the Navier boundary conditions \equ(NavierBC).
So a solution $u$ of \equ(NS) and \equ(NavierBC)
can be obtained via \equ(uxyExpansion) from a solution $\Phi$
of the equation \equ(curlNS).
For convenience, we write \equ(curlNS) as
$$
(\alpha\partial_t-\Delta)\Phi+\thalf\gamma\LOP(\Phi)\Phi=\nsv f\,,
\equation(scNS)
$$
where $\LOP$ is the symmetric bilinear form defined by
$$
\LOP(\phi)\psi=
(\nabla\phi)\cdot\nsv^{-1}\psi
+(\nabla\psi)\cdot\nsv^{-1}\phi\,.
\equation(LOPDef)
$$
The coefficients $\Phi_{j,k}$ in the series \equ(PhixyExpansion)
are $2\pi$-periodic functions and thus admit an expansion
$$
\Phi_{j,k}=\sum_{n\in\oldinteger}\Phi_{n,j,k}\cosi_n\,,\qquad
\cosi_n(t)=\cases{\cos(nt) &if $n\ge 0$,\cr\sin(-nt) &if $n<0$.\cr}
\equation(PhijkSeries)
$$
Denote by $\oldnatural_0$ the set of all nonnegative integers.
For any subset $N\subset\oldnatural_0$ we define
$$
\mean_N\Phi=\sum_{n\in\oldinteger\atop|n|\in N}\sum_{j,k\in\oldnatural_1}
\Phi_{n,j,k}\cosi_n\times\sin_j\times\sin_k\,,
\equation(FreqProj)
$$
where $\sin_m(z)=\sin(mz)$.
In particular, the even frequency part $\Phi_\even$
(odd frequency part $\Phi_\odd$) of $\Phi$ is defined
to be the function $\mean_N\Phi$,
where $N$ is the set of all even (odd) nonnegative integers.
This leads to the decomposition $\Phi=\Phi_\even+\Phi_\odd$
that will be used below.

To simplify the discussion,
consider first non-stationary periodic solutions.
For $\gamma$ near the bifurcation point $\gamma_0$,
we expect $\Phi$ to be nearly time-independent.
So in particular, $\Phi_\odd$ is close to zero.
Consider the function $\phi=\phi_\even+\phi_\odd$
obtained by setting $\phi_\even=\Phi_\even$ and $\phi_\odd=\beta^{-1}\Phi_\odd$.
The scaling factor $\beta\ne 0$ will be chosen below,
in such a way that $\phi_\even$ and $\phi_\odd$ are of comparable size.
Substituting
$$
\Phi=U_\beta\phi\defeq\phi_\even+\beta\phi_\odd
\equation(PhiDecomp)
$$
into \equ(scNS) yields the equation
$$
(\alpha\partial_t-\Delta)\phi+\thalf\gamma\LOP_s(\phi)\phi=\nsv f\,,
\equation(scNSs)
$$
where $s=\beta^2$ and
$$
\LOP_s(\phi)\psi
=\LOP(\phi_\even)\psi_\even+\LOP(\phi_\even)\psi_\odd
+\LOP(\phi_\odd)\psi_\even+s\LOP(\phi_\odd)\psi_\odd\,.
\equation(LOPsDef)
$$
Finally, we convert \equ(scNSs) to a fixed point equation
by applying the inverse of $\alpha\partial_t-\Delta$ to both sides.
Setting $g=(-\Delta)^{-1}\nsv f$,
the resulting equation is $\tilde\phi=\phi$, where
$$
\tilde\phi=g
-\thalf\gamma|\Delta|^{1/2}(\alpha\partial_t-\Delta)^{-1}\hat\phi\,,\qquad
\hat\phi\defeq|\Delta|^{-1/2}\LOP_s(\phi)\phi\,.
\equation(scNSsFix)
$$

One of the features of the equation \equ(scNSs)
is that  the time-translate of a solution is again a solution.
We eliminate this symmetry by imposing the condition $\phi_{1,1,1}=0$.
In addition, we choose $\beta=\Beth^{-1}\Phi_{-1,1,1}$,
where $\Beth$ is some fixed constant that will be specified later.
This leads to the normalization conditions
$$
A\phi\defeq\phi_{1,1,1}=0\,,\qquad
B\phi\defeq\phi_{-1,1,1}=\Beth\,.
\equation(ABNormaliz)
$$
Notice that $\beta$ enters our main equation $\tilde\phi=\phi$
only via its square $s=\beta^2$.
It is convenient to regard $s$ to be the independent parameter
and express $\gamma$ as a function of $s$.
The functions $\gamma=\gamma(s)$ and $\alpha=\alpha(s)$
are determined by the condition that $\tilde\phi$
satisfies the normalization conditions \equ(ABNormaliz).
Applying the functionals $A$ and $B$ to both sides of \equ(scNSs),
using the identities
$A\Delta=-2A$, $A\partial_t=B$, $B\Delta=-2B$, $B\partial_t=-A$,
and imposing the conditions $A\tilde\phi=0$ and $B\tilde\phi=\Beth$,
we find that
$$
\gamma=-2^{3/2}{\Beth\over B\hat\phi}\,,\qquad
\alpha=2{A\hat\phi\over B\hat\phi}\,.
\equation(gammaalphaSol)
$$
For a fixed value of $s$, define $\FF_s(\phi)=\tilde\phi$,
where $\tilde\phi$ is given by \equ(scNSsFix),
with $\gamma=\gamma(s,\phi)$ and $\alpha=\alpha(s,\phi)$
determined by \equ(gammaalphaSol).
The fixed point equation for $\FF_s$ is used
to find non-stationary time-periodic solutions of \equ(scNSs).

\demo Remark(NotNormalized)
The choice \equ(gammaalphaSol) guarantees that
$A\tilde\phi=0$ and $B\tilde\phi=\Beth$,
even if $\phi$ does not satisfy the normalization conditions \equ(ABNormaliz).
Thus, the domain of the map $\FF_s$
can include non-normalized function $\phi$.
(The same is true for the map $\FF_\gamma$ described below.)
But a fixed point of this map will be normalized by construction.

In order to determine the bifurcation point $\gamma_0$
and the corresponding frequency $\alpha_0$,
we consider the map $\FF:\phi\mapsto\tilde\phi$ given by \equ(scNSsFix) with $s=0$.
The values of $\gamma$ and $\alpha$ are again given by \equ(gammaalphaSol),
so that $A\tilde\phi=0$ and $B\tilde\phi=\Beth$.
We will show that this map $\FF$ has a fixed point
$\phi$ with the property that $\phi_{n,j,k}=0$ whenever $|n|>1$.
The values of $\gamma$ and $\alpha$ for this fixed point
define $\gamma_0$ and $\alpha_0$.

A similar map $\FF_\gamma:\phi\mapsto\tilde\phi$,
given by \equ(scNSsFix) with $s=0$,
is used to find stationary solutions of the equation \equ(scNS).
In this case, the value of $\gamma$ is being fixed,
and $\phi_\odd$ is taken to be zero.
The goal is to show that this map $\FF_\gamma$
has a fixed point $\phi_\gamma$ that is independent of time $t$.
Then $\Phi=\phi_\gamma$ is a stationary solution of \equ(scNS).

\medskip
We finish this section by computing
the derivative of the map $\FF_s$ described after \equ(gammaalphaSol).
The resulting expressions will be needed later.
Like some of the above, the following is purely formal.
A proper formulation will be given in the next section.
For simplicity, assume that $\phi$ depends on a parameter.
The derivative of a quantity $q$ with respect to this parameter
will be denoted by $\dot q$.
Define
$$
\LL_\alpha=|\Delta|^{1/2}(\alpha\partial_t-\Delta)^{-1}\,,\qquad
\LL_\alpha'=\partial_t(\alpha\partial_t-\Delta)^{-1}\,.
\equation(LLalphaDef)
$$
Using that $\FF_s(\phi)=g-\thalf\gamma\LL_\alpha\hat\phi$
with $\hat\phi=|\Delta|^{-1/2}\LOP_s(\phi)\phi$,
the parameter-derivative of $\FF_s(\phi)$ is given by
$$
D\FF_s(\phi)\dot\phi
=-\thalf\LL_\alpha\Bigl[\bigl(\dot\gamma-\gamma\dot\alpha\LL_\alpha'\bigr)\hat\phi
+\gamma\dot{\hat\phi}\Bigr]\,,\qquad
\dot{\hat\phi}=2|\Delta|^{-1/2}\LOP_s(\phi)\dot\phi\,,
\equation(DFFsPhiDotPhi)
$$
where
$$
\dot\gamma
=2^{-3/2}{\gamma^2\over\Beth} B\dot{\hat\phi}\,,\qquad
\dot\alpha
=2^{-3/2}{\alpha\gamma\over\Beth}B\dot{\hat\phi}
-2^{-1/2}{\gamma\over\Beth}A\dot{\hat\phi}\,.
\equation(DotgammaDotalpha)
$$
The above expressions for $\dot\gamma$ and $\dot\alpha$
are obtained by differentiating \equ(gammaalphaSol).

%%%%%%%%%%%%%%%%%%%%%%%%%%%%%%%%%%%%
%%%%%%%%%%%%%%%%%%%%%%%%%%%%%%%%%%%%
\section The associated contractions
%%%%%%%%%%%%%%%%%%%%%%%%%%%%%%%%%%%%
%%%%%%%%%%%%%%%%%%%%%%%%%%%%%%%%%%%%

In this section, we formulate the fixed point problems
for the maps $\FF$, $\FF_\gamma$, and $\FF_s$
in a suitable functional setting.
The goal is to reduce the problems to a point
where we can invoke the contraction mapping theorem.
After describing the necessary estimates,
we give a proof of \clm(NSHopfBif) based on these estimates.

We start by defining suitable function spaces.
Given a real number $\rho>1$, denote by $\AA$
the space of all functions $h\in\rmL^2(\torus)$
that have a finite norm $\|h\|$, where
$$
\|h\|=|h_0|
+\sum_{n\in\oldnatural_1}\sqrt{|h_n|^2+|h_{-n}|^2}\rho^n\,,\qquad
h=\sum_{n\in\oldinteger}h_n\cosi_n\,.
\equation(AADef)
$$
Here $\cosi_n$ are the trigonometric function defined in \equ(PhijkSeries).
It is straightforward to check that $\AA$
is a Banach algebra under the pointwise product of functions.
That is, $\|gh\|\le\|g\|\|h\|$ for any two functions $g,h\in\AA$.
We also identify functions on $\torus$ with $2\pi$-periodic functions on $\real$.
In this sense, a function in $\AA$
extends analytically to the strip $T(\rho)=\{z\in\complex:|\Im z|<\log\rho\}$.

Given in addition $\varrho>1$, we denote by $\buB$
the space of all function $\Phi:\torus^2\to\AA$
that admit a representation \equ(PhixyExpansion)
and have a finite norm
$$
\|\Phi\|
=\sum_{j,k\in\oldnatural_1}\|\Phi_{j,k}\|\varrho^{j+k}\,.
\equation(GGrhoepsNorm)
$$
A function $(x,y)\mapsto(t\mapsto\Phi(t,x,y))$ in this space
will also be identified with a function
$(t,x,y)\mapsto\Phi(t,x,y)$ on $\torus^3$,
or with a function on $\real^3$ that is $2\pi$-periodic in each argument.
In this sense, every function in $\buB$
extends analytically to $T(\rho)\times T(\varrho)^2$.

We consider $\AA$ and $\buB$ to be Banach spaces
over $\field\in\{\real,\complex\}$.
In the case $\field=\real$, the functions in these spaces
are assumed to take real values for real arguments.

\smallskip
Clearly, a function $\Phi\in\buB$
admits an expansion \equ(FreqProj) with $N=\oldnatural_0$.
The sequence of Fourier coefficients $\Phi_{n,k,j}$
converges to zero exponentially as $|n|+j+k$ tends to infinity.
If all but finitely many of these coefficients vanish,
then $\Phi$ is called a Fourier polynomial.
The equation \equ(FreqProj) with $N\subset\oldnatural_0$ non-empty
defines a continuous projection $\mean_N$ on $\buB$
whose operator norm is $1$.
Using Fourier series, it is straightforward
to see that the equation \equ(LLalphaDef)
defines two bounded linear operators $\LL_\alpha$ and $\LL_\alpha'$
on $\buB$, for every $\alpha\in\complex$.
The operator $\LL_\alpha$ is in fact compact.
Specific estimates will be given in Section 4.
The following will be proved in Section 4 as well.

\claim Proposition(LOPBound)
If $\Phi$ and $\phi$ belong to $\buB$,
then so does $|\Delta|^{-1/2}\LOP(\Phi)\phi$, and
$$
\bigl\||\Delta|^{-1/2}\LOP(\Phi)\phi\bigr\|
\le\bigl\||\Delta|^{-1/2}\Phi\bigr\|\|\phi\|
+\|\Phi\|\bigl\||\Delta|^{-1/2}\phi\bigr\|\,.
\equation(LOPBound)
$$

This estimate implies e.g.~that the transformation $\phi\mapsto\tilde\phi$,
given by \equ(scNSsFix) for fixed values of $s$, $\gamma$ and $\alpha$,
is well-defined and compact as a map from $\buB$ to $\buB$.

As is common in computer-assisted proofs,
we reformulate the fixed point equation for the map $\phi\mapsto\tilde\phi$
as a fixed point problem for an associated quasi-Newton map.
Since we need three distinct versions of this map,
let us first describe a more general setting.

\medskip
Let $\FF:\DD\to\BB$ be a $\rmC^1$ map
defined on an open domain $\DD$ in a Banach space $\BB$.
Let $h\mapsto\varphi+Lh$ be a continuous affine map on $\BB$.
We define quasi-Newton map $\NN$ for $(\DD,\FF,\varphi,L)$
by setting
$$
\NN(h)=\FF(\varphi+Lh)-\varphi+(\id-L)h\,.
\equation(quasiNewtonMap)
$$
The domain of $\NN$ is defined to be the set of of all $h\in\BB$
with the property that $\varphi+Lh\in\DD$.
Notice that, if $h$ is a fixed point of $\NN$,
then $\varphi+Lh$ is a fixed point of $\FF$.
In our applications, $\varphi$ is an approximate fixed point of $\FF$
and $L$ is an approximate inverse of $\id-D\FF(\varphi)$.

The following is an immediate consequence of the contraction mapping theorem.

\claim Proposition(ContrMappingThm)
Let $\FF:\DD\to\BB$ be a $\rmC^1$ map
defined on an open domain in a Banach space $\BB$.
Let $h\mapsto\varphi+Lh$ be a continuous affine map on $\BB$.
Assume that the quasi-Newton map \equ(quasiNewtonMap)
includes a non-empty ball $B_\delta=\{h\in\BB: \|h\|<\delta\}$ in its domain,
and that
$$
\|\NN(0)\|<\eps\,,\qquad\|D\NN(h)\|<K\,,\qquad h\in B_\delta\,,
\equation(ContrMappingThm)
$$
where $\eps,K$ are positive real numbers that satisfy $\eps+K\delta<\delta$.
Then $\FF$ has a fixed point in $\varphi+LB_\delta$.
If $L$ is invertible, then this fixed point is unique in $\varphi+LB_\delta$.

In our applications below, $\BB$ is always a subspace of $\buB$.
The domain parameter $\rho$ and the constant $\Beth$
that appears in the normalization condition \equ(ABNormaliz)
are chosen to have the fixed values
$$
\rho=2^5\,,\qquad\Beth=2^{-12}\,.
\equation(ParamValues)
$$
The domain parameter $\varrho$ is defined implicitly in our proofs.
That is, the lemmas below
hold for $\varrho>1$ sufficiently close to $1$.

Consider first the problem of determining
the bifurcation point $\gamma_0$ and the associated frequency $\alpha_0$.
Let $\BB=\mean_{\sss\{0,1\}}\buB$ over $\real$.
For every $\delta>0$ define $B_\delta=\{h\in\BB: \|h\|<\delta\}$.
Let $s=0$, and denote by $\DD$ the set of all functions $\phi\in\BB$
with the property that $B\hat\phi\ne 0$.
Define $\FF:\DD\to\BB$ to be the map $\phi\mapsto\tilde\phi$ given by \equ(scNSsFix),
with $\gamma=\gamma(\phi)$ and $\alpha=\alpha(\phi)$
defined by the equation \equ(gammaalphaSol).
Clearly, $\FF$ is not only $\rmC^1$ but real analytic on $\DD$.

\claim Lemma(BifPoint)
With $\FF$ as described above,
there exists an affine isomorphism $h\mapsto\varphi+L_1h$ of $\BB$
and real numbers $\eps,\delta,K>0$ satisfying $\eps+K\delta<\delta$,
such that the following holds.
The quasi-Newton map $\NN$ associated with $(\BB,\FF,\varphi,L_1)$
includes the ball $B_\delta$ in its domain
and satisfies the bounds \equ(ContrMappingThm).
The domain of $\FF$ includes the ball in $\BB$
of radius $r=\delta\|L_1\|$, centered at $\varphi$.
For every function $\phi$ in this ball,
$\gamma(\phi)=83.1733117\ldots$ and $\alpha(\phi)=4.66592275\ldots$.

Our proof of this lemma is computer-assisted
and will be described in Section 4.

By \clm(ContrMappingThm),
the map $\FF$ has a unique fixed point $\phi^\ast\in\varphi+L_1 B_\delta$.
We define $\gamma_0=\gamma(\phi^\ast)$ and $\alpha_0=\alpha(\phi^\ast)$.

Our next goal is to construct
a branch of periodic solutions for the equation \equ(scNS).
Consider $\BB=\buB$ over $\field\in\{\real,\complex\}$.
By continuity, there exists an open ball $\JJ_0\subset\field$ centered at the origin,
and an open neighborhood $\DD$ of $\phi^\ast$ in $\BB$, such that
$B\hat\phi=B|\Delta|^{-1/2}\LOP_s(\phi)\phi$ is nonzero for all $s\in\JJ_0$
and all $\phi\in\DD$.
For every $s\in\JJ_0$,
define $\FF_s:\DD\to\BB$ to be the map $\phi\mapsto\tilde\phi$
given by \equ(scNSsFix), with $\gamma=\gamma(s,\phi)$
and $\alpha=\alpha(s,\phi)$ defined by the equation \equ(gammaalphaSol).

\claim Lemma(PerSolutions)
Let $\field=\real$.
There exists a isomorphism $L$ of $\buB$
such that the following holds.
If $\NN_0$ denotes the the quasi-Newton map
associated with $(\DD,\FF_0,\phi^\ast,L)$,
then the derivative $D\NN_0(0)$ of $\NN_0$
at the origin is a contraction.

Our proof of this lemma is computer-assisted
and will be described in Section 4.
As a consequence we have the following.

\claim Corollary(PerBranch)
Consider $\field=\complex$.
There exists an open disk $\JJ\subset\complex$,
centered at the origin,
and an analytic curve $s\mapsto\phi_s$ on $\JJ$ with values in $\DD$,
such that $\FF_s(\phi_s)=\phi_s$ for all $s\in\JJ$.
If $s$ belongs to the real interval $\JJ\cap\real$, then $\phi_s$ is real.
Furthermore, $\phi_0=\phi^\ast$.

\proof
Consider still $\field=\complex$.
For $s\in\II_0$, the derivative of $\NN_s$ on its domain is given by
$$
D\NN_s(h)=D\FF_s(\phi^\ast+Lh)L+\id-L\,.
\equation(DNNsh)
$$
Assume that some function $\psi\in\buB$ satisfies
$D\FF_0(\phi^\ast)\psi=\psi$.
We may assume that $\psi$ takes real values for real arguments.
A straightforward computation shows that
$D\NN_0(0)L^{-1}\psi=L^{-1}\psi$.
Since $D\NN_0(0)$ is a contraction in the real setting,
by \clm(PerSolutions), this implies that $\psi=0$.
So the operator $D\FF_0(\phi^\ast)$ does not have an eigenvalue $1$.
This operator is compact,
since it is the composition of a bounded linear operator
with the compact operator $\LL_\alpha$.
Thus, $D\FF_0(\phi^\ast)$ has no spectrum at $1$.
By the implicit function theorem,
there exists a complex open ball $\JJ$, centered at the origin,
such that the fixed point equation $\FF_s(\phi)=\phi$
has a solution $\phi=\phi_s$ for all $s\in\JJ$.
Furthermore, the curve $s\mapsto\phi_s$ is analytic,
passes through $\phi^\ast$ at $s=0$,
and there is a unique curve with this property.
By uniqueness, we also have $\ov{\phi_{\bar s}}=\phi_s$ for all $s\in\JJ$,
so $\phi_s$ is real for real values of $s\in\JJ$.
\qed

A branch of stationary periodic solutions for \equ(scNS)
is obtained similarly.
Consider $\BB=\mean_{\sss\{0\}}\buB$ over $\field\in\{\real,\complex\}$.
For every $\gamma\in\field$,
define $\FF_\gamma:\BB\to\BB$ to be the map $\phi\mapsto\tilde\phi$
given by \equ(scNSsFix), with $s=\alpha=0$.
Notice that $\phi^\ast_\even$ is a fixed point of $\FF_{\gamma_0}$.

\claim Lemma(StatSolutions)
Let $\field=\real$. There exists an isomorphism $L_0$ of $\BB$
such that the following holds.
If $\NN_{\gamma_0}$ denotes the the quasi-Newton map
associated with $(\BB,\FF_{\gamma_0},\phi^\ast_\even,L_0)$,
then the derivative $D\NN_{\gamma_0}(0)$ of $\NN_{\gamma_0}$
at the origin is a contraction.

Our proof of this lemma is computer-assisted
and will be described in Section 4.
As a consequence we have the following.

\claim Corollary(StatBranch)
Consider $\field=\complex$.
There exists an open disk $\II\subset\complex$,
centered at $\gamma_0$,
and an analytic curve $\gamma\mapsto\phi_\gamma$ on $\II$ with values in $\BB$,
such that $\FF_\gamma(\phi_\gamma)=\phi_\gamma$ for all $\gamma\in\II$.
If $\gamma$ belongs to the real interval $\II\cap\real$, then $\phi_\gamma$ is real.
Furthermore, $\phi_{\gamma_0}=\phi^\ast_\even$.

The proof of this corollary is analogous
to the proof of \clm(PerBranch).

We note that the disk $\II\ni\gamma_0$
is disjoint from the disk $\JJ\ni 0$ described in \clm(PerBranch).
So there is no ambiguity in using the notation
$\gamma\mapsto\phi_\gamma$ and $s\mapsto\phi_s$
for the curve of stationary and periodic solutions,
respectively, of the equation \equ(scNSs),

\smallskip
Based on the results stated in this section, we can now give a

\proofof(NSHopfBif)
As described in the preceding sections,
the curve $\gamma\mapsto\phi_\gamma$ for $\gamma\in\II$
yields a curve $\gamma\mapsto u_\gamma$ of stationary solutions of
the equation \equ(NS), where $u_\gamma=\nsv^{-1}\phi_\gamma$.
By our choice of function spaces,
the function $(\gamma,x,y)\mapsto u_\gamma(x,y)$
is real analytic on $I\times\torus^2$, where $I=\II\cap\real$.

Similarly, the curve $s\mapsto\phi_s$ for $s\in\JJ$
defines a family of of non-stationary periodic solutions
for \equ(NS), with $\gamma=\gamma_s$ and $\alpha=\alpha_s$
determined via the equation \equ(gammaalphaSol).
To be more precise,
the even frequency part $\phi_{s,\even}$ of $\phi_s$
determines a vector field $u_{s,\even}=\nsv^{-1}\phi_{s,\even}$,
and the odd frequency part $\phi_{s,\odd}$
determines a vector field $u_{s,\odd}=\nsv^{-1}\phi_{s,\odd}$.
If $\beta$ is a complex number such that $s=\beta^2\in\JJ$,
then $u=u_{s,\even}+\beta u_{s,\odd}$
is a periodic solution of \equ(NS), with $\gamma=\gamma_s$ and $\alpha=\alpha_s$.
Here, we have used the decomposition \equ(PhiDecomp).
By our choice of function spaces, the functions
$(s,t,x,y)\mapsto u_{s,\even}(t,x,y)$ and $(s,t,x,y)\mapsto u_{s,\odd}(t,x,y)$
are real analytic on $J\times\torus^3$, where $J=\JJ\cap\real$.
Clearly, $\partial_t u_{0,\odd}(t,\bdot\,,\bdot)\ne 0$,
due to the normalization condition $\phi_{-1,1,1}=\Beth$ imposed in \equ(ABNormaliz).
And by construction, we have $u=u_{\gamma_0}$ for $s=0$.
\qed

%%%%%%%%%%%%%%%%%%%%%%%%%%%%
%%%%%%%%%%%%%%%%%%%%%%%%%%%%
\section Remaining estimates
%%%%%%%%%%%%%%%%%%%%%%%%%%%%
%%%%%%%%%%%%%%%%%%%%%%%%%%%%

What remains to be proved are Lemmas
\clmno(BifPoint), \clmno(PerSolutions), and \clmno(StatSolutions).
Our method used in the proof of \clm(BifPoint)
can be considered perturbation theory
about the approximate fixed point $\varphi$ of $\FF$.
The function $\varphi$ is a Fourier polynomial with over $20000$ nonzero coefficients,
so a large number of estimates are involved.

We start by describing bounds on the bilinear function $\LOP$
and on the linear operators $\LL_\alpha$ and $\LL_\alpha'$.
These are the basic building blocks for our transformations
$\FF$, $\FF_s$, and $\FF_\gamma$.
The ``mechanical'' part of these estimates will be described in Subsection 4.4.

%%%%%%%%%%%%%%%%%%%%%%%%%%%%%%%%%%%%%%%%%%%%%%%%%%%%%%%%%%%%%%%%%%%%%%%%%%%
\subsection The bilinear form $\boldmath\LOP$ and a proof of \clm(LOPBound)
%%%%%%%%%%%%%%%%%%%%%%%%%%%%%%%%%%%%%%%%%%%%%%%%%%%%%%%%%%%%%%%%%%%%%%%%%%%

Consider the bilinear form $\LOP$ defined by \equ(LOPDef).
Using the identity \equ(curlInv), we have
$$
\eqalign{
\LOP(\Phi)\phi
&=(\nabla\Phi)\cdot\iso\nabla\Delta^{-1}\phi
+(\nabla\phi)\cdot\iso\nabla\Delta^{-1}\Phi\cr
&=\bigl[(\partial_x\Phi)\Delta^{-1}\partial_y\phi
-(\partial_y\Phi)\Delta^{-1}\partial_x\phi\bigr]
-\bigl[(\Delta^{-1}\partial_x\Phi)\partial_y\phi
-(\Delta^{-1}\partial_y\Phi)\partial_x\phi\bigr]\,.\cr}
\equation(LOPPhiphi)
$$
In order to obtain accurate estimates,
it is useful to have explicit expressions for $\LOP(\Phi)\phi$
in terms of the Fourier coefficients of $\Phi$ and $\phi$.
Given that $\LOP$ is bilinear,
and that the identity \equ(LOPPhiphi) holds pointwise in $t$,
it suffices to compute $\LOP(\Phi)\phi$
for the time-independent monomials
$$
\Phi=\sin_J\times\sin_K\,,\qquad
\phi=\sin_j\times\sin_k\,,
\equation(PhiphiModes)
$$
with $J,K,j,k>0$.
A straightforward computation shows that
$$
\eqalign{
\LOP(\Phi)\phi&=\Theta(Jk+jK)\bigl[
\sin_{J+j}\times\sin_{K-k}-\sin_{J-j}\times\sin_{K+k}
\bigr]\cr
&\quad+\Theta(Jk-jK)\bigl[
\sin_{J+j}\times\sin_{K+k}-\sin_{J-j}\times\sin_{K-k}
\bigr]\,,\cr}
\equation(LOPPhiphiOne)
$$
with $\Theta$ as defined below.
As a result we have
$$
|\Delta|^{-1/2}\LOP(\Phi)\phi
=\sum_{\sigma,\tau=\pm 1}N_{\sigma,\tau}\sin_{\sigma J+j}\times\sin_{\tau K+k}\,,
\equation(InvHalfLapLOPDecomp)
$$
where
$$
N_{\sigma,\tau}
=\Theta
{\sigma Jk-\tau Kj\over\sqrt{(\sigma J+j)^2+(\tau K+k)^2}}\,,\qquad
\Theta=\quarter\biggl({1\over J^2+K^2}-{1\over j^2+k^2}\biggr)\,.
\equation(NpmpmDef)
$$

\proofof(LOPBound)
Using the Cauchy-Schwarz inequality in $\real^2$, we find that
$$
|N_{\sigma,\tau}|
=|\Theta|{|(\sigma J+j)k-(\tau K+k)j|\over\sqrt{(\sigma J+j)^2+(\tau K+k)^2}}
\le|\Theta|\sqrt{j^2+k^2}\,.
\equation(NppBound)
$$
Since the absolute value of $N_{\sigma,\tau}$ is invariant
under an exchange of $(j,k)$ and $(J,K)$, this implies that
$$
|N_{\sigma,\tau}|\le{1/4\over\sqrt{j^2+k^2}}\vee{1/4\over\sqrt{J^2+K^2}}\,,
\equation(NpmpmBound)
$$
where $a\vee b=\max(a,b)$ for $a,b\in\real$.
As a result, we obtain the bound
$$
\bigl\||\Delta|^{-1/2}\LOP(\Phi)\phi\bigr\|
\le\bigl\||\Delta|^{-1/2}\Phi\bigr\|_{\varrho,\epsilon}\|\phi\|
+\|\Phi\|\bigl\||\Delta|^{-1/2}\phi\bigr\|\,.
\equation(LOPPhiphiBound)
$$
Using the nature of the norm \equ(GGrhoepsNorm),
and the fact that $\AA$ is a Banach algebra
for the pointwise product of functions,
this bound extends by bilinearity to arbitrary
functions $\Phi,\phi\in\buB$.
\qed

We note that the bound \equ(LOPPhiphiBound)
exploits the cancellations that lead to the expression \equ(LOPPhiphiOne).
A more straightforward estimate loses a factor of $2$
with respect to \equ(LOPPhiphiBound).
But it is not just this factor of $2$ that counts for us.
The expressions \equ(NpmpmDef) for the coefficients $N_{\sigma,\tau}$
and the bounds \equ(NpmpmBound) are used in our computations
and error estimates.
The expression on the right hand side of \equ(NpmpmBound)
is a decreasing function of the wavenumbers $j,k,J,K$,
so it can be used to estimate $\LOP(\Phi)\phi$
when $\Phi$ and/or $\phi$ are ``tails'' of Fourier series.

%%%%%%%%%%%%%%%%%%%%%%%%%%%%%%%%%%%%%%%%%%%%%%%%%%%%%%%%%%%%%%%%%%%%%%%%%%%%%%%
\subsection The linear operators $\boldmath\LL_\alpha$ and $\boldmath\LL_\alpha'$
%%%%%%%%%%%%%%%%%%%%%%%%%%%%%%%%%%%%%%%%%%%%%%%%%%%%%%%%%%%%%%%%%%%%%%%%%%%%%%%

Consider the linear operators $\LL_\alpha$
and $\LL_\alpha'$ defined in \equ(LLalphaDef), with $\alpha$ real.
A straightforward computation shows that
$$
\psi_{n,j,k}
=\sqrt{j^2+k^2}\,{(j^2+k^2)\phi_{n,j,k}-\alpha n\phi_{-n,j,k}
\over(j^2+k^2)^2+\alpha^2 n^2}\,,\qquad
\psi=\LL_\alpha\phi\,.
\equation(LLalphaphinjk)
$$
Using the Cauchy-Schwarz inequality in $\real^2$,
this yields the estimate
$$
\sqrt{|\psi_{n,j,k}|^2+|\psi_{-n,j,k}|^2}
\le C_{n,j,k}\sqrt{|\phi_{n,j,k}|^2+|\phi_{-n,j,k}|^2}\,,
\equation(LLalphaBound)
$$
with
$$
C_{n,j,k}=\sqrt{j^2+k^2\over(j^2+k^2)^2+\alpha^2 n^2}
\le{1\over\sqrt{2|\alpha n|}}\wedge {1\over\sqrt{j^2+k^2}}
\equation(LLalphaBoundConst)
$$
for $n\ne 0$, where $a\wedge b=\min(a,b)$ for $a,b\in\real$.
The last bound in \equ(LLalphaBoundConst)
is a decreasing function of $|n|,j,k$
and can be used to estimate $\LL_\alpha\phi$
when $\phi$ is the tail of a Fourier series.

For the operator $\LL_\alpha'$ we have
$$
\psi_{n,j,k}
=n\,{(j^2+k^2)\phi_{-n,j,k}+\alpha n\phi_{n,j,k}
\over(j^2+k^2)^2+\alpha^2 n^2}\,,\qquad
\psi=\LL_\alpha'\phi\,.
\equation(LLalphaPrimenjk)
$$
A bound analogous to \equ(LLalphaBound)
holds for $\psi=\LL_\alpha'\phi$, with
$$
C_{n,j,k}=\sqrt{n^2\over(j^2+k^2)^2+\alpha^2 n^2}\,.
\equation(LLalphaPrimeBoundConst)
$$
As can be seen from \equ(DFFsPhiDotPhi),
this bound is needed only for $n=\pm 1$,
since these are the only nonzero frequencies
of the function $\hat\phi=|\Delta|^{-1/2}\LOP_0(\phi)\phi$
with $\phi\in\mean_{\sss\{0,1\}}\buB$.
And for fixed $n$, the right hand side of \equ(LLalphaPrimeBoundConst)
is decreasing in $j$ and $k$.

%%%%%%%%%%%%%%%%%%%%%%%%%%%%%%%%%%%%%
\subsection Estimating operator norms
%%%%%%%%%%%%%%%%%%%%%%%%%%%%%%%%%%%%%

Recall that a function $\phi\in\buB$ admits a Fourier expansion
$$
\phi=\sum_{n\in\oldinteger}\;\sum_{j,k\in\oldnatural_1}\phi_{n,j,k}\theta_{n,j,k}\,,\qquad
\theta_{n,j,k}\defeq\cosi_n\times\sin_j\times\sin_k\,,
\equation(thetanjkDef)
$$
and that the norm of $\phi$ is given by
$$
\|\phi\|=\sum_{j,k\in\oldnatural_1}
\biggl[|\phi_{0,j,k}|+\sum_{n\in\oldnatural_1}
\sqrt{|\phi_{n,j,k}|^2+|\phi_{-n,j,k}|^2}\,\rho^n\biggr]\varrho^{j+k}\,.
\equation(GGrhoepsNormAgain)
$$
Let now $n\ge 0$.
A linear combination $c_{\sss+}\theta_{n,j,k}+c_{\sss-}\theta_{-n,j,k}$
will be referred to as a mode with frequency $n$ and wavenumbers $(j,k)$
or as a mode of type $(n,j,k)$.
We assume of course that $c_{\sss-}=0$ when $n=0$.
Since \equ(GGrhoepsNormAgain) is a weighted $\ell^1$ norm,
except for the $\ell^2$ norm used for modes,
we have a simple expression for the operator norm
of a continuous linear operator $\LL:\buB\to\buB$, namely
$$
\|\LL\|=\sup_{j,k\in\oldnatural_1}\;
\sup_{n\in\oldnatural_0}\;\sup_u\|\LL u\|/\|u\|\,,
\equation(OpNorm)
$$
where the third supremum is over all nonzero modes
$u$ of type $(n,j,k)$.

Let now $n,j,k\ge 1$ be fixed.
In computation where $\LL\theta_{\pm n,j,k}$ is known explicitly,
we use the following estimate.
Denote by $\LL_{n,j,k}$ the restriction of $\LL$
to the subspace spanned by the two functions $\theta_{\pm n,j,k}$.
For $q\ge 1$ define
$$
\|\LL_{n,j,k}\|_q
=\sup_{0\le p<q}\|\LL v_p\|\,,\qquad
v_p=\cos\Bigl({\pi p\over q}\Bigr){\theta_{n,j,k}\over\rho^n\varrho^{j+k}}
+\sin\Bigl({\pi p\over q}\Bigr){\theta_{-n,j,k}\over\rho^n\varrho^{j+k}}\,.
\equation(mNorm)
$$
Since every unit vector in the span of $\theta_{\pm n,j,k}$
lies within a distance less than ${\pi\over q}$
of one of the vectors $v_p$ or its negative, we have
$\|\LL_{n,j,k}\|\le\|\LL_{n,j,k}\|_q+{\pi\over q}\|\LL_{n,j,k}\|$.
Thus
$$
\|\LL_{n,j,k}\|\le{q\over q-\pi}\|\LL_{n,j,k}\|_m\,,\qquad q\ge 4\,.
\equation(LLmodeNormBound)
$$

Consider now the operator
$D\FF_s(\phi)$ described in \equ(DFFsPhiDotPhi),
with $\phi\in\mean_{\sss\{0,1\}}\buB$ fixed.
If $\dot\phi=u_n$ is a nonzero mode with frequency $n\ge 3$,
then $\dot{\hat\phi}=2|\Delta|^{-1/2}\LOP_0(\phi)\dot\phi$
belongs to $\mean_N\buB$ with $N=\{n-1,n,n+1\}$.
Thus, we have $\dot\gamma=\dot\alpha=0$, and
$$
D\FF_0(\phi)u_n
=-\gamma\LL_\alpha|\Delta|^{-1/2}\LOP_0(\phi)u_n\,.
\equation(DFFophiu)
$$
Due to the factor $\LL_\alpha$ in this equation,
if $u_n=c_{\sss+}\theta_{n,j,k}+c_{\sss-}\theta_{-n,j,k}$
with $(j,k)$ and $c_{\sss\pm}$ fixed,
then the ratios
$$
\|D\FF_0(\phi)u_n\|/\|u_n\|
\equation(DFFophiunRatios)
$$
are decreasing in $n$ for $n\ge 3$.
And the limit as $n\to\infty$ of this ratio is zero.

So for the operator $\LL=D\FF_0(\phi)$,
the supremum over $n\in\oldnatural_0$ in \equ(OpNorm)
reduces to a maximum over finitely many terms.
The same holds for the operator $\LL=D\NN_0(0)=D\FF_0(\phi^\ast)L+\id-L$
that is described in \clm(PerSolutions).
This is a consequence of the following choice.

\demo Remark(LMatrix)
The operator $L$ chosen in \clm(PerSolutions)
is a ``matrix perturbation'' of the identity,
in the sense that $L\theta_{n,j,k}=\theta_{n,j,k}$
for all but finitely many indices $(n,j,k)$.
The same is true for the operators $L_1$
and $L_0$ chosen in \clm(BifPoint) and \clm(StatSolutions),
respectively.

%%%%%%%%%%%%%%%%%%%%%%%%%%%%%%
\subsection Computer estimates
%%%%%%%%%%%%%%%%%%%%%%%%%%%%%%

Lemmas \clmno(BifPoint), \clmno(StatSolutions),
and \clmno(PerSolutions) assert the existence of certain objects
that satisfy a set of strict inequalities.
The goal here is to construct these objects,
and to verify the necessary inequalities
by combining the estimates that have been described so far.

The above-mentioned ``objects'' are real numbers,
real Fourier polynomials, and linear operators
that are finite-rank perturbations of the identity.
They are obtained via purely numerical computations.
Verifying the necessary inequalities is largely an organizational task,
once everything else has been set up properly.
Roughly speaking, the procedure follows that of a well-designed numerical program,
but instead of truncation Fourier series and ignoring rounding errors,
we determine rigorous enclosures at every step along the computation.
This part of the proof is written in the programming language Ada [\rAda].
The following is meant to be a rough guide for the reader who
wishes to check the correctness of our programs.
The complete details can be found in [\rProgs].

\smallskip
An enclosure for a function $\phi\in\buB$
is a set in $\buB$ that includes $\phi$
and is defined in terms of (bounds on) a Fourier polynomial
and finitely many error terms.
We define such sets hierarchically,
by first defining enclosures for elements in simpler spaces.
In this context, a ``bound'' on a map $f:\XX\to\YY$
is a function $F$ that assigns to a set $X\subset\XX$
of a given type ({\tt Xtype}) a set $Y\subset\YY$
of a given type ({\tt Ytype}), in such a way that
$y=f(x)$ belongs to $Y$ for all $x\in X$.
In Ada, such a bound $F$ can be implemented by defining
a {\tt procedure F(X:{\gapii}in Xtype; Y:{\gapii}out Ytype)}.

Our most basic enclosures are specified by pairs {\tt S=(S.C,S.R)},
where {\tt S.C} is a representable real number ({\tt Rep})
and {\tt S.R} a nonnegative representable real number ({\tt Radius}).
Given a Banach algebra $\XX$ with unit ${\bf 1}$,
such a pair {\tt S} defines a ball in $\XX$ which we denote by
$\langle{\tt S},\XX\rangle=\{x\in\XX:\|x-({\tt S.C}){\bf 1}\|\le{\tt S.R}\}$.

When $\XX=\real$,
then the data type described above is called {\tt Ball}.
Bounds on some standard functions involving the type {\tt Ball}
are defined in the package {\tt Flts\_Std\_Balls}.
Other basic functions are covered in the packages {\tt Vectors} and {\tt Matrices}.
Bounds of this type have been used in many computer-assisted proofs;
so we focus here on the more problem-specific aspects of our programs.

Consider now the space $\AA$
for a fixed domain radius $\varrho>1$ of type {\tt Radius}.
As mentioned before \dem(LMatrix),
we only need to consider Fourier polynomials in $\AA$.
Our enclosures for such polynomials are defined by
an {\tt array(-I$_{\rmc}${\gapi}..{\gapi}I$_{\rmc}$) of Ball}.
This data type is named {\tt NSPoly},
and the enclosure associated with data {\tt P} of this type is
$$
\langle{\tt P},\AA\rangle
\defeq\sum_{i=-I_{\rmc}}^{I_{\rmc}}
\bigl\langle{\tt P(i)},\real\bigr\rangle\cosi_{\nu(i)}\,,
\equation(NSPolyEnclosure)
$$
where $\nu$ is an increasing index function with the property that $\nu(-i)=-\nu(i)$.
The type {\tt NSPoly} is defined in the package {\tt NSP},
which also implements bounds on some basic operations
for Fourier polynomials in $\AA$.
Among the arguments to {\tt NSP} is a nonnegative integer $n$
(named {\tt NN}).
Our proof of \clm(StatSolutions) and \clm(BifPoint) uses
$I_c=n=0$ and $I_c=n=1$, respectively, and $\nu(i)=i$.
Values $n\ge 2$ are uses when estimating the norm of $\LL u$
for the operator $\LL=D\NN_0(0)$, with $u$ a mode of frequency $n$.
In this case, $\nu$ takes values in $\{-n,n\}$ or $\{-n-1,-n,-n+1,0,n-1,n,n+1\}$,
depending on whether $n$ is odd or even.
(The value $\nu=0$ is being used only for $n=2$.)
The package {\tt NSP} also defines a data type {\tt NSErr}
as an {\tt array(0{\gapi}..{\gapi}I$_{\rm c}$) of Radius}.
This type will be used below.

Given in addition a positive number $\varrho\ge 1$ of type {\tt Radius},
our enclosures for functions in $\buB$
are defined by pairs {\tt(F.C,F.E)},
where {\tt F.C} is an
{\tt array(1{\gapi}..{\gapi}J$_{\rmc}$,1{\gapi}..{\gapi}K$_{\rmc}$) of NSPoly}
and {\tt F.E} is an
{\tt array(1{\gapi}..{\gapi}J$_{\rme}$,1{\gapi}..{\gapi}K$_{\rme}$) of NSErr};
all for a fixed value of the parameter {\tt NN}.
This data type is named {\tt Fourier3},
and the enclosure associated with {\tt F=(F.C,F.E)} is
$$
\langle{\tt F},\buB\rangle
\defeq\sum_{j=1}^{J_\rmc}\sum_{k=1}^{K_\rmc}\bigl\langle{\tt F.C(j,k)},\AA\bigr\rangle
\times\sin_j\times\sin_k
+\sum_{J=1}^{J_\rme}\sum_{K=1}^{K_\rme}H_{\sss J,K}({\tt F.E(J,K)})\,.
\equation(FouThreeEnclosure)
$$
Here, $H_{\sss J,K}({\tt E})$ denotes the set of all functions
$\phi=\sum_{i=0}^{I_{\rmc}}\phi^i$ with $\|\phi^i\|\le{\tt E(i)}$,
where $\phi^i$ can be any function in $\buB$
whose coefficients $\phi^i_{n,j,k}$ vanish unless
$j\ge J$, $k\ge K$, and $|n|=\nu(i)$.

The type {\tt Fourier3} and bounds on some standard functions
involving this type are defined in the child package {\tt NSP.Fouriers}.
This package is a modified version of the
package {\tt Fouriers2} that was used earlier in [\rAKxii,\rAKxviii,\rAGK].
The procedure {\tt Prod} is now a bound on the bilinear map $|\Delta|^{-1/2}\LOP_0$.
The error estimates used in {\tt Prod} are based on the inequality \equ(NpmpmBound).
The package {\tt NSP.Fouriers}
also includes bounds {\tt InvLinear} and {\tt DtInvLinear}
on the linear operators $\LL_\alpha$ and $\LL_\alpha'$, respectively.
These bounds use the estimates described in Subsection 4.3.

As far as the proof of \clm(BifPoint) is concerned,
it suffices now to compose existing bounds
to obtain a bound on the map $\FF$ and its derivative $D\FF$.
This is done by the procedures {\tt GMap} and {\tt DGMap} in {\tt Hopf.Fix}.
Here we use enclosures of type {\tt NN=1}.

The type of quasi-Newton map $\NN$ defined by \equ(quasiNewtonMap)
has been used in several computer-assisted proof before.
So the process of constructing a bound on $\NN$
from a bound on $\FF$ has been automated
in the generic packages {\tt Linear} and {\tt Linear.Contr}.
(Changes compared to earlier versions are mentioned in the program text.)
This includes the computation of an approximate inverse $L_1$
for the operator $\id-D\FF(\varphi)$.
A bound on $\NN$ is defined (in essence) by the procedure {\tt Linear.Contr.Contr},
instantiated with {\tt Map => GMap}.
And a bound on $D\NN$ is defined by {\tt Linear.Contr.Contr},
with {\tt DMap => DGMap}.
Bounds on operator norms are obtained via {\tt Linear.OpNorm}.
Another problem-dependent ingredient in these procedures,
besides {\tt Map} and {\tt DMap}, are data of type {\tt Modes}.
These data are constructed by the procedure {\tt Make}
in the package {\tt Hopf}.
They define a splitting of the given space $\BB$ into a finite direct sum.
For details on how such a splitting is defined and used we refer to [\rAKxix].

If the parameter {\tt NN} has the value $0$,
then the procedures {\tt GMap} and {\tt DGMap} define
bounds on the map $\FF_\gamma$ and its derivative, respectively.
The operator $L_0$ used in \clm(StatSolutions)
has the property that $M_0=L_0-\id$
satisfies $M_0=P_0M_0P_0$, where $P_0=\mean_{\sss\{0\}}\proj_{m_0}$
for some positive integer $m_0$.
Here, and in what follows, $\proj_m$ denotes
the canonical projection in $\buB$
with the property that $\proj_m\phi$ is obtained
from $\phi$ by restricting the second sum in \equ(thetanjkDef)
to wavenumbers $j,k\le m$.

If {\tt NN} has a value $n\ge 2$, then the procedure {\tt DGMap}
defines a bound on the map $(\phi,\psi)\mapsto D\FF_0(\phi)\psi$,
restricted to the subspace $\mean_{\sss\{0,1\}}\buB\times\mean_{\sss\{n\}}\buB$.
The linear operator $L$ that is used in \clm(PerSolutions)
admits a decomposition $L=\id+M_1+M_2+\ldots+M_N$ of the following type.
After choosing a suitable sequence $n\mapsto m_n$ of positive integers,
we set $M_n=P_n(L-\id)P_n$,
where $P_1=\mean_{\sss\{0,1\}}\proj_{m_1}$
and $P_n=\mean_{\sss\{n\}}\proj_{m_n}$ for $n=2,3,\ldots,N$.
This structure of $L$ simplifies the use of \equ(OpNorm)
for estimating the norm of $\LL=D\NN_0(0)$.
Furthermore, to check that $L$ is invertible,
it suffices to verify that $\id+M_n$ is invertible
on the finite-dimensional subspace $P_n\buB$,
for each positive $n\le N$.

The linear operator $L_1$ that is used in \clm(BifPoint)
is of the form $L_1=\id+M_1$ with $M_1$ as described above.

All the steps required in the proofs
of Lemmas \clmno(BifPoint), \clmno(StatSolutions), and \clmno(PerSolutions)
are organized in the main program {\tt Check}.
As $n$ ranges from $0$ to $N=305$,
this program defines the parameters that are used in the proof
for {\tt NN} $=n$, instantiates the necessary packages,
computes the appropriate matrix $M_n$,
verifies that $\id+M_n$ is invertible,
reads $\varphi$ from the file {\tt BP.approx},
and then calls the procedure {\tt ContrFix} from the
(instantiated version of the) package {\tt Hopf.Fix} to verify the necessary inequalities.

The representable numbers ({\tt Rep}) used in our programs
are standard [\rIEEE] extended floating-point numbers (type {\tt LLFloat}).
High precision [\rMPFR] floating-point numbers (type {\tt MPFloat})
are used as well, but not in any essential way.
Both types support controlled rounding.
{\tt Radius} is always a subtype of {\tt LLFloat}.
Our programs were run successfully on a $20$-core workstation,
using a public version of the gcc/gnat compiler [\rGnat].
For further details,
including instruction on how to compile and run our programs,
we refer to [\rProgs].

\bigskip
%%%%%%%%%%%
\references
%%%%%%%%%%%

{\ninepoint

\item{[\rHopf]} E.~Hopf,
{\sl Abzweigung einer periodischen L\"osung
von einer station\"aren L\"osung eines Differentialsystems},
Ber. Math.-Phys. Kl. Siichs. Akad. Wiss. Leipzig, {\bf 94}, 3--22 (1942).

\item{[\rSerrin]} J.~Serrin,
{\sl A Note on the Existence of Periodic Solutions of the Navier-Stokes Equations},
Arch. Rational Mech. Anal. {\bf 3} 120--122, (1959).
%%% periodic forcing

\item{[\rCrAb]} M.G.~Crandall and P.H.~Rabinowitz,
{\sl The Hopf bifurcation theorem in infinite dimensions},
Arch. Rational Mech. Anal. {\bf 67}, 53--72 (1977).
%\pdfclink{0 0 1}{online here}
%{https://www.semanticscholar.org/paper/The-Hopf-Bifurcation-Theorem-in-infinite-dimensions-Crandall-Rabinowitz/c21e703fdf69882c64cff84eb976cc16ddf7345e}

\item{[\rMarMcC]} J.~Marsden, M.~McCracken,
{\sl The Hopf bifurcation and its applications},
Springer Applied Mathematical Sciences Lecture Notes Series, Vol. 19, 1976.

\item{[\rRT]} D.~Ruelle, F.~Takens, {\sl On the Nature of Turbulence},
Commun. Math. Phys. 20, 167--192 (1971)

\item{[\rKlWe]} P.~Kloeden, R.~Wells,
{\sl An explicit example of Hopf bifurcation in fluid mechanics},
Proc. Roy. Soc. London Ser. A {\bf 390}, 293--320 (1983).
%\pdfclink{0 0 1}{online here}
%{https://royalsocietypublishing.org/doi/10.1098/rspa.1983.0133}

\item{[\rChIo]} P.~Chossat, G.~Iooss,
{\sl Primary and secondary bifurcations in the Couette-Taylor problem},
Japan J. Appl. Math. {\bf 2}, 37--68 (1985).

\item{[\rDum]} F.~Dumortier,
{\sl Techniques in the theory of local bifurcations:
blow-up, normal forms, nilpotent bifurcations, singular perturbations};
in: {\sl Bifurcations and periodic orbits of vector fields},
(D.~Schlomiuk, ed., Kluwer Acad.~Pub.)
NATO ASI Ser. C Math. Phys. Sci. {\bf 408}, 10--73 (1993).

\item{[\rChIoo]} P.~Chossat, G.~Iooss,
{\sl The Couette-Taylor problem},
Applied Mathematical Sciences, 102. Springer-Verlag, New York, 1994

\item{[\rNWYNK]} M.T.~Nakao, Y.~Watanabe, N.~Yamamoto, T.~Nishida, M.-N.~Kim,
{\sl Computer assisted proofs of bifurcating solutions for nonlinear heat convection problems},
J. Sci. Comput. {\bf 43}, 388--401 (2010).
%\pdfclink{0 0 1}{online here}{https://doi.org/10.1007/s10915-009-9303-3}
%% Oberbeck-Boussinesq equation (reduction 3d --> 2d) for The Rayleigh-Benard problem.
%% (extension to 3d is discussed at the end)
%% Consider only bifurcations of stationary solutions
%% Similarly for the referencces 9 and 12 that they cite (on the same problem)

\item{[\rAKxii]} G.~Arioli, H.~Koch,
{\sl Non-symmetric low-index solutions for a symmetric boundary value problem},
J. Differ. Equations {\bf 252}, 448--458 (2012).

\item{[\rAKxv]} G.~Arioli, H.~Koch,
{\sl Some symmetric boundary value problems and non-symmetric solutions},
J. Differ. Equations {\bf 259}, 796--816 (2015).

\item{[\rGaldi]} G.P.~Galdi,
{\sl On bifurcating time-periodic flow of
a Navier-Stokes liquid past a cylinder},
Arch. Rational Mech. Anal. {\bf 222}, 285--315 (2016).
Digital Object Identifier (DOI) 10.1007/s00205-016-1001-3

\item{[\rHJNS]} C.-H.~Hsia, C.-Y.~Jung, T.B.~Nguyen, and M.-C.~Shiu,
{\sl On time periodic solutions,
asymptotic stability and bifurcations of Navier-Stokes equations},
Numer. Math. {\bf 135}, 607--638 (2017).
%%% time-periodic forcing, and reference
%%% found this ref in [\rBBLV]

\item{[\rAKxviii]} G.~Arioli, H.~Koch,
{\sl Spectral stability for the wave equation with periodic forcing},
J. Differ. Equations {\bf 265}, 2470--2501 (2018).

\item{[\rAKxix]} G.~Arioli, H.~Koch,
{\sl Non-radial solutions for some semilinear elliptic equations on the disk},
Nonlinear Analysis {\bf 179}, 294–308 (2019).

\item{[\rNPW]} M.T.~Nakao, M.~Plum, Y.~Watanabe,
{\sl Numerical verification methods and computer-assisted proofs
for partial differential equations},
Springer Series in Computational Mathematics, Vol. 53,
Springer Singapore, 2019

\item{[\rGoSe]} J.~G\'omez-Serrano,
{\sl Computer-assisted proofs in PDE: a survey},
SeMA {\bf 76}, 459--484 (2019).

\item{[\rWiZg]} D.~Wilczak, P.~Zgliczy\'nski,
{\sl A geometric method for infinite-dimensional chaos:
Symbolic dynamics for the Kuramoto-Sivashinsky PDE on the line},
J. Differ. Equations {\bf 269}, 8509--8548 (2020).

\item{[\rBBLV]} J.~B.~van den Berg, M.~Breden, J.-P.~Lessard, L.~van Veen,
{\sl Spontaneous periodic orbits in the Navier-Stokes flow},
Preprint 2019,
%\pdfclink{0 0 1}{online here}{https://arxiv.org/abs/1902.00384}

\item{[\rAGK]} G.~Arioli, F.~Gazzola, H.~Koch,
{\sl Uniqueness and bifurcation branches for planar steady
Navier-Stokes equations under Navier boundary conditions},
Preprint 2020.

\item{[\rBLQ]} J.~B.~van den Berg, J.-P.~Lessard, E.~Queirolo,
{\sl Rigorous verification of Hopf bifurcations
via desingularization and continuation},
Preprint 2020.

\item{[\rBQ]} J.~B.~van den Berg, E.~Queirolo,
{\sl Validating Hopf bifurcation in the Kuramoto-Sivashinky PDE},
in preparation.
% talk by Elena Queirolo, see
% https://researchseminars.org/seminar/CRM-CAMP

\item{[\rProgs]} G.~Arioli, H.~Koch,
{\sl Programs and data files for the proof of Lemmas \clmno(BifPoint),
\clmno(StatSolutions), \clmno(PerSolutions), and \clmno(JoinBranches)},
\pdfclink{0 0 1}{{\tt https://web.ma.utexas.edu/users/koch/papers/nshopf/}}
{https://web.ma.utexas.edu/users/koch/papers/nshopf/}

%%%%%%%%%%%%%%%%%%%%%%%%%%%%%%%%%%%%

\item{[\rAda]} Ada Reference Manual, ISO/IEC 8652:2012(E),
available e.g. at\hfil\break
\pdfclink{0 0 1}{{\tt www.ada-auth.org/arm.html}}
{http://www.ada-auth.org/arm.html}

\item{[\rGnat]}
A free-software compiler for the Ada programming language,
which is part of the GNU Compiler Collection; see
\pdfclink{0 0 1}{{\tt gnu.org/software/gnat/}}{http://gnu.org/software/gnat/}

\item{[\rIEEE]} The Institute of Electrical and Electronics Engineers, Inc.,
{\sl IEEE Standard for Binary Float\-ing--Point Arithmetic},
ANSI/IEEE Std 754--2008.

\item{[\rMPFR]} The MPFR library for multiple-precision floating-point computations
with correct rounding; see
\pdfclink{0 0 1}{{\tt www.mpfr.org/}}{http://www.mpfr.org/}

}

\bye
%%%%%%%%%%%%%%%%%%%%%%%%%%%%%%% including param.2 %%%%%%%%%%%%%%%%%%%%%%%%%%%%%%%
%%%%%%%%%%%%%%%%%%%%%%%%%%%%%%%%% end of param.2 %%%%%%%%%%%%%%%%%%%%%%%%%%%%%%%%
%%%%%%%%%%%%%%%%%%%%%%%%%%%%%%% including fonts.6 %%%%%%%%%%%%%%%%%%%%%%%%%%%%%%%
%%%%%%%%%%%%%%%%%%%%%%%%%%%%%%%%% end of fonts.6 %%%%%%%%%%%%%%%%%%%%%%%%%%%%%%%%
%%%%%%%%%%%%%%%%%%%%%%%%%%%%%%% including smallfonts.tex %%%%%%%%%%%%%%%%%%%%%%%%%%%%%%%
%%%%%%%%%%%%%%%%%%%%%%%%%%%%%%%%% end of smallfonts.tex %%%%%%%%%%%%%%%%%%%%%%%%%%%%%%%%
%%%%%%%%%%%%%%%%%%%%%%%%%%%%%%% including symbols.1 %%%%%%%%%%%%%%%%%%%%%%%%%%%%%%%
%%%%%%%%%%%%%%%%%%%%%%%%%%%%%%%%% end of symbols.1 %%%%%%%%%%%%%%%%%%%%%%%%%%%%%%%%
%%%%%%%%%%%%%%%%%%%%%%%%%%%%%%% including links.1 %%%%%%%%%%%%%%%%%%%%%%%%%%%%%%%
%%%%%%%%%%%%%%%%%%%%%%%%%%%%%%%%% end of links.1 %%%%%%%%%%%%%%%%%%%%%%%%%%%%%%%%
%%%%%%%%%%%%%%%%%%%%%%%%%%%%%%% including titles.9 %%%%%%%%%%%%%%%%%%%%%%%%%%%%%%%
%%%%%%%%%%%%%%%%%%%%%%%%%%%%%%%%% end of titles.9 %%%%%%%%%%%%%%%%%%%%%%%%%%%%%%%%
%%%%%%%%%%%%%%%%%%%%%%%%%%%%%%% including macros.21 %%%%%%%%%%%%%%%%%%%%%%%%%%%%%%%
%%%%%%%%%%%%%%%%%%%%%%%%%%%%%%%%% end of macros.21 %%%%%%%%%%%%%%%%%%%%%%%%%%%%%%%%
%%%%%%%%%%%%%%%%%%%%%%%%%%%%%%% including mygraphicx.tex %%%%%%%%%%%%%%%%%%%%%%%%%%%%%%%
%% modification of graphicx.tex by Nathan Goldschmidt
\input miniltx

\ifx\pdfoutput\undefined
  \def\Gin@driver{dvips.def}  % we are not running PDFTeX
\else
  \def\Gin@driver{pdftex.def} % we are running PDFTeX
\fi
 
\input graphicx.sty
\resetatcatcode
%%%%%%%%%%%%%%%%%%%%%%%%%%%%%%%%% end of mygraphicx.tex %%%%%%%%%%%%%%%%%%%%%%%%%%%%%%%%
%%%%%%%%%%%%%%%%%%%%%%%%%%%%%%% including opmac.1 %%%%%%%%%%%%%%%%%%%%%%%%%%%%%%%
%%%%%%%%%%%%%%%%%%%%%%%%%%%%%%%%% end of opmac.1 %%%%%%%%%%%%%%%%%%%%%%%%%%%%%%%%
%%%%%%%%%%%%%%%%%%%%%%%%%%%%%%% including boldmath.tex %%%%%%%%%%%%%%%%%%%%%%%%%%%%%%%
%%%%%%%%%%%%%%%%%%%%%%%%%%%%%%%%% end of boldmath.tex %%%%%%%%%%%%%%%%%%%%%%%%%%%%%%%%
%\input param.2
%\input fonts.6
%\input smallfonts.tex
%\input symbols.1
%\input links.1
%\input titles.9
%\input macros.21
%\input boldmath.tex
%\input mygraphicx.tex
%
\newdimen\savedparindent
\savedparindent=\parindent
\font\tenamsb=msbm10 \font\sevenamsb=msbm7 \font\fiveamsb=msbm5
\newfam\bbfam
\textfont\bbfam=\tenamsb
\scriptfont\bbfam=\sevenamsb
\scriptscriptfont\bbfam=\fiveamsb
\def\field{{\mathds F}}
\def\bbb{\fam\bbfam}
\def\oldinteger{{\bbb Z}}
\def\oldnatural{{\bbb N}}
\def\buB{{\hbox{\magnineeufm B}}}
\let\Beth\theta
\def\gapi{\hskip 2pt}
\def\gapii{\hskip 2pt}
\def\rmc{{\rm c}}
\def\rme{{\rm e}}
\def\id{{\rm I}}
\def\LOP{{\mathds L}}
\def\nsv{\mathop{\not\!\partial}\nolimits}
\def\cosi{\mathop{\rm cosi}\nolimits}
\def\even{{\rm e}}
\def\odd{{\rm o}}
\def\bdot{\hbox{\bf .}}
\def\stwomat#1#2#3#4{{\eightpoint\left[\matrix{#1&#2\cr#3&#4\cr}\right]}}
\addref{Hopf}
\addref{Serrin}
\addref{CrAb}
\addref{MarMcC}
\addref{RT}
\addref{KlWe}
\addref{ChIo}
\addref{Dum}
\addref{ChIoo}
\addref{NWYNK}
\addref{AKxii}
\addref{AKxv}
\addref{Galdi}
\addref{HJNS}
\addref{AKxviii}
\addref{AKxix}
\addref{NPW}
\addref{GoSe}
\addref{WiZg}
\addref{BBLV}
\addref{AGK}
\addref{BLQ}
\addref{BQ}
\addref{Progs}
\addref{Ada}
\addref{Gnat}
\addref{IEEE}
\addref{MPFR}
\def\leftheadline{\sixrm\hfil ARIOLI \& KOCH\hfil}
\def\rightheadline{\sevenrm\hfil Hopf bifurcation in Navier-Stokes\hfil}
%
%%%%%%%%%%%%%%%%%%%%%%%%%%%%%%%%%%%%%%%%%%%%%%%%%%%%%%%
\cl{\magtenbf A Hopf bifurcation in the planar Navier-Stokes equations}
\bigskip

\cl{
Gianni Arioli
\footnote{$^1$}
{\eightpoint\hskip-2.9em
Department of Mathematics, Politecnico di Milano,
Piazza Leonardo da Vinci 32, 20133 Milano.
}
$^{\!\!\!,\!\!}$
\footnote{$^2$}
{\eightpoint\hskip-2.6em
Supported in part by the PRIN project
``Equazioni alle derivate parziali e
disuguaglianze analitico-geometriche associate''.}
and Hans Koch
\footnote{$^3$}
{\eightpoint\hskip-2.7em
Department of Mathematics, The University of Texas at Austin,
Austin, TX 78712.}
}

\bigskip
\abstract
We consider the Navier-Stokes equation
for an incompressible viscous fluid on a square,
satisfying Navier boundary conditions
and being subjected to a time-independent force.
As the kinematic viscosity is varied,
a branch of stationary solutions is shown
to undergo a Hopf bifurcation,
where a periodic cycle branches from the stationary solution.
Our proof is constructive
and uses computer-assisted estimates.

%%%%%%%%%%%%%%%%%%%%%%%%%%%%%%%%%%%%%
%%%%%%%%%%%%%%%%%%%%%%%%%%%%%%%%%%%%%
\section Introduction and main result
%%%%%%%%%%%%%%%%%%%%%%%%%%%%%%%%%%%%%
%%%%%%%%%%%%%%%%%%%%%%%%%%%%%%%%%%%%%

We consider the Navier-Stokes equations
$$
\partial_t u-\nu\Delta u+(u\cdot\nabla)u+\nabla p=f\ ,\quad
\nabla\cdot u=0\quad{\rm on~}\Omega\,,
\equation(NavierStokes)
$$
for the velocity $u=u(t,x,y)$
of an incompressible fluid on a planar domain $\Omega$,
satisfying suitable boundary conditions for $(x,y)\in\partial\Omega$
and initial conditions at $t=0$.
Here, $p$ denotes the pressure,
and $f=f(x,y)$ is a fixed time-independent external force.

Our focus is on solution curves and bifurcations
as the kinematic velocity $\nu$ is being varied.
In order to reduce the complexity of the problem,
the domain $\Omega$ is chosen to be as simple as possible,
namely the square $\Omega=(0,\pi)^2$.
Following [\rAGK], we impose Navier boundary conditions on $\partial\Omega$,
which are given by
$$
\eqalign{
u_1&=\partial_x u_2=0\quad{\rm on~}\{0,\pi\}\times(0,\pi)\,,\cr
u_2&=\partial_y u_1=0\quad{\rm on~}(0,\pi)\times\{0,\pi\}\,.\cr}
\equation(NavierBC)
$$
A fair amount is known about the (non)uniqueness of stationary solutions
in this case [\rAGK].
This includes the existence of a bifurcation between curves of stationary solutions
with different symmetries.

Here we prove the existence of a Hopf bifurcation
for the equation \equ(NavierStokes) with boundary conditions \equ(NavierBC),
and with a forcing function $f$ that satisfies
$$
(\partial_xf_2-\partial_yf_1)(x,y)=5\sin(x)\sin(2y)-13\sin(3x)\sin(2y)\,.
\equation(specificNScurlf)
$$
In a Hopf bifurcation, a stationary solution loses stability
and a small-amplitude limit cycle branches from the stationary solution [\rHopf,\rCrAb,\rMarMcC].
Among other things, this introduces a time scale in the system
and increases its complexity.
In this capacity, Hopf bifurcations in the Navier-Stokes equation
constitute an important first step
in the transition to turbulence in fluids,
as was described in the seminal work [\rRT].

Numerically, there is plenty of evidence that Hopf bifurcations occur
in the Navier-Stokes equation, but proofs are still very scarce.
An explicit example of a Hopf bifurcation
was given in [\rKlWe] for the rotating B\'enard problem.
A proof exists also for the Couette-Taylor problem [\rChIo,\rChIoo].
Sufficient conditions for the existence of a Hopf bifurcation
in a Navier-Stokes setting are presented in [\rGaldi].

Before giving a precise statement of our result, let us replace
the vector field $u$ in the equation \equ(NavierStokes) by $\nu^{-1}u$.
The equation for the rescaled function $u$ is
$$
\alpha\partial_t u-\Delta u+\gamma(u\cdot\nabla)u+\nabla p=f\ ,\quad
\nabla\cdot u=0\quad{\rm on~}\Omega\,,
\equation(NS)
$$
where $\gamma=\nu^{-2}$.
The value of $\alpha$ that corresponds to \equ(NavierStokes) is $\nu^{-1}$,
but this can be changed to any positive value by rescaling time.

Numerically, it is possible to find stationary solutions of \equ(NS)
for a wide range of values of the parameter $\gamma$.
At a value $\gamma_0\approx 83.1733117\ldots$ we observe
a Hopf bifurcation that leads to a branch of periodic solutions for $\gamma>\gamma_0$.

For a fixed value of $\alpha$,
the time period $\tau$ of the solution varies with $\gamma$.
Instead of looking for $\tau$-periodic solution of \equ(NS)
for fixed $\alpha$, we look for $2\pi$-periodic solutions,
where $\alpha=2\pi/\tau$ has to be determined.
To simplify notation, a $2\pi$-periodic function
will be identified with a function on the circle $\torus=\real/(2\pi\integer)$.
Our main theorem is the following.

\claim Theorem(NSHopfBif)
There exists a real number $\gamma_0=83.1733117\ldots$,
an open interval $I$ including $\gamma_0$,
and a real analytic function $(\gamma,x,y)\mapsto u_\gamma(x,y)$
from $I\times\Omega$ to $\real^2$,
such that $u_\gamma$ is a stationary solution of \equ(NS) and \equ(NavierBC)
for each $\gamma\in I$.
In addition, there exists a real number $\alpha_0=4.66592275\ldots$,
an open interval $J$ centered at the origin,
two real analytic functions $\gamma$ and $\alpha$ on $J$
that satisfy $\gamma(0)=\gamma_0$ and $\alpha(0)=\alpha_0$, respectively,
as well as two real analytic functions
$(s,t,x,y)\mapsto u_{s,\even}(t,x,y)$ and $(s,t,x,y)\mapsto u_{s,\odd}(t,x,y)$
from $J\times\torus\times\Omega$
to $\real^2$, such that the following holds.
For any given $\beta\in\complex$ satisfying $\beta^2\in J$,
the vector field $u=u_{s,\even}+\beta u_{s,\odd}$ with $s=\beta^2$
is a solution of \equ(NS) and \equ(NavierBC)
with $\gamma=\gamma(s)$ and $\alpha=\alpha(s)$.
Furthermore, $u_{0,\even}(t,\bdot\,,\bdot)=u_{\gamma_0}$
and $\partial_t u_{0,\odd}(t,\bdot\,,\bdot)\ne 0$.

To our knowledge, this is the first result establishing
the existence of a Hopf bifurcation for the Navier-Stokes equation
in a stationary environment.

Our proof of this theorem is computer-assisted.
The solutions are obtained by rewriting \equ(NS) and \equ(NavierBC)
as a suitable fixed point equation for scalar vorticity of $u$.
Here we take advantage of the fact that the domain is two-dimensional.
We isolate the periodic branch from the stationary branch
by using a scaling that admits two distinct limits at the bifurcation point.
This approach is also known as the blow-up method,
which is a common tool in the study of singularities and bifurcations [\rDum].

\smallskip
Computer-assisted methods have been applied successfully
to many different problems in analysis,
mostly in the areas of dynamical systems and partial differential equations.
Here we will just mention work that
concerns the Navier-Stokes equation or Hopf bifurcations.
For the Navier-Stokes equation,
the existence of symmetry-breaking bifurcations among stationary solutions
has been established in [\rNWYNK,\rAGK].
Periodic solutions for the Navier-Stokes flow in a stationary environment
have been obtained in [\rBBLV].
In the case of periodic forcing,
the problem of existence and stability of periodic orbits
has been investigated in [\rHJNS].
Concerning the existence of Hopf bifurcations,
a computer-assisted proof was given recently in [\rBLQ]
for a finite-dimensional dynamical system;
and an extension of their method to the Kuramoto-Sivashinsky PDE
is presented in [\rBQ].
For other recent computer-assisted proofs
we refer to [\rAKxix,\rNPW,\rGoSe,\rWiZg] and references therein.

Figure 1 depicts snapshots at $t=0$ and $t=\pi$
of a solution $u:\torus\times\Omega\to\real^2$ of the
equations \equ(NS) with boundary conditions \equ(NavierBC)
and forcing \equ(specificNScurlf),
obtained numerically for the parameter value $\gamma\approx 84.00\ldots$.

%%%%%%%%%%%%%%%%%%%%%%%%%%%%%%%%%%%%%%%%%%%%%%%%%%%%%%%%%%%%%
\vskip0.15in
\hbox{\hskip 55pt
\includegraphics[height=1.5in,width=1.5in]{figures/ps1.eps}
\hskip 35pt
\includegraphics[height=1.5in,width=1.5in]{figures/ps2.eps}}
\vskip0.1in
\centerline{\eightpoint {\noindent\bf Figure 1.}
Snapshots at two distinct times of a time-periodic solution
for $\gamma\approx 84.00\ldots$
}
\vskip0.15in
%%%%%%%%%%%%%%%%%%%%%%%%%%%%%%%%%%%%%%%%%%%%%%%%%%%%%%%%%%%%%

As mentioned earlier, a system similar to the one
considered here is known to exhibit a symmetry-breaking bifurcation
within the class of stationary solutions [\rAGK].
The broken symmetry is $y\mapsto\pi/2-y$.
Based on a numerical computation of eigenvalues,
we expect an analogous bifurcation to occur here at $\gamma\approx 1450$.
Interestingly, the Hopf bifurcation described here occurs
at a significantly smaller value of $\gamma$.
We have not tried to prove the existence
of a symmetry-breaking bifurcation for the forcing \equ(specificNScurlf),
since such an analysis would duplicate the work in [\rAGK]
and go beyond the scope of the present paper.

\smallskip
The remaining part of this paper is organized as follows.
In Section 2, we first rewrite \equ(NS) as an equation for
the function $\Phi=\partial_y u_1-\partial_x u_2$,
which is the scalar vorticity of $-u$.
After a suitable scaling $\Phi=U_\beta\phi$,
the problem of constructing the solution branches
described in \clm(NSHopfBif) is reduced to three fixed point problems
for the function $\phi$.
These fixed point equations are solved in Section 3,
based on estimates described in Lemmas 3.3, 3.4, and 3.6.
Section 4 is devoted to the proof of these estimates,
which involves reducing them to a large number of trivial bounds
that can be (and have been) verified with the aid of a computer [\rProgs].

%%%%%%%%%%%%%%%%%%%%%%%%%%%%%%
%%%%%%%%%%%%%%%%%%%%%%%%%%%%%%
\section Fixed point equations
%%%%%%%%%%%%%%%%%%%%%%%%%%%%%%
%%%%%%%%%%%%%%%%%%%%%%%%%%%%%%

The goal here is to rewrite the equation \equ(NS)
with boundary conditions \equ(NavierBC) as a fixed point problem.
Applying the operator $\nsv:(u_1,u_2)\mapsto\partial_2 u_1-\partial_1 u_2$
on both sides of the equation \equ(NS), we obtain
$$
\alpha\partial_t\Phi
-\Delta\Phi+\gamma u\cdot\nabla\Phi=\nsv f\,,\qquad\Phi=\nsv u\,.
\equation(curlNS)
$$
Here, we have used that $\nsv\,(u\cdot\nabla)u=u\cdot\nabla\Phi$.
Using the divergence-free condition $\nabla\cdot u=0$,
one also finds that
$$
\Delta u=\iso\nabla\Phi\,,\qquad
\iso=\stwomat{0}{1}{-1}{0}\,.
\equation(Lapu)
$$
If $\Phi$ vanishes on the boundary of $\partial\Omega$,
then the equation \equ(Lapu) can be inverted to yield
$$
u=\nsv^{-1}\Phi\defeq\iso\nabla\Delta^{-1}\Phi\,,
\equation(curlInv)
$$
where $\Delta$ denotes the Dirichlet Laplacean on $\Omega$.

In Section 3
we will define a space of real analytic functions $\Phi$
that admit a representation
$$
\Phi(t,x,y)=\sum_{j,k\in\oldnatural_1}\Phi_{j,k}(t)\sin(jx)\sin(ky)\,,
\equation(PhixyExpansion)
$$
with the series converging uniformly on a complex
open neighborhood of $\torus^3$.
Here, and in what follows, $\oldnatural_1$ denotes the set of all positive integers.
If $\Phi$ admits such an expansion,
then the equation \equ(curlInv) yields
$$
\eqalign{
u_1(t,x,y)&=\sum_{j,k\in\oldnatural_1}\,{-k\over j^2+k^2}\Phi_{j,k}(t)\sin(jx)\cos(ky)\,,\cr
u_2(t,x,y)&=\sum_{j,k\in\oldnatural_1}\,{j\over j^2+k^2}\Phi_{j,k}(t)\cos(jx)\sin(ky)\,.\cr}
\equation(uxyExpansion)
$$
It is straightforward to check that the corresponding
vector field $u=(u_1,u_1)$ satisfies the Navier boundary conditions \equ(NavierBC).
So a solution $u$ of \equ(NS) and \equ(NavierBC)
can be obtained via \equ(uxyExpansion) from a solution $\Phi$
of the equation \equ(curlNS).
For convenience, we write \equ(curlNS) as
$$
(\alpha\partial_t-\Delta)\Phi+\thalf\gamma\LOP(\Phi)\Phi=\nsv f\,,
\equation(scNS)
$$
where $\LOP$ is the symmetric bilinear form defined by
$$
\LOP(\phi)\psi=
(\nabla\phi)\cdot\nsv^{-1}\psi
+(\nabla\psi)\cdot\nsv^{-1}\phi\,.
\equation(LOPDef)
$$
The coefficients $\Phi_{j,k}$ in the series \equ(PhixyExpansion)
are $2\pi$-periodic functions and thus admit an expansion
$$
\Phi_{j,k}=\sum_{n\in\oldinteger}\Phi_{n,j,k}\cosi_n\,,\qquad
\cosi_n(t)=\cases{\cos(nt) &if $n\ge 0$,\cr\sin(-nt) &if $n<0$.\cr}
\equation(PhijkSeries)
$$
Denote by $\oldnatural_0$ the set of all nonnegative integers.
For any subset $N\subset\oldnatural_0$ we define
$$
\mean_N\Phi=\sum_{n\in\oldinteger\atop|n|\in N}\sum_{j,k\in\oldnatural_1}
\Phi_{n,j,k}\cosi_n\times\sin_j\times\sin_k\,,
\equation(FreqProj)
$$
where $\sin_m(z)=\sin(mz)$.
In particular, the even frequency part $\Phi_\even$
(odd frequency part $\Phi_\odd$) of $\Phi$ is defined
to be the function $\mean_N\Phi$,
where $N$ is the set of all even (odd) nonnegative integers.
This leads to the decomposition $\Phi=\Phi_\even+\Phi_\odd$
that will be used below.

To simplify the discussion,
consider first non-stationary periodic solutions.
For $\gamma$ near the bifurcation point $\gamma_0$,
we expect $\Phi$ to be nearly time-independent.
So in particular, $\Phi_\odd$ is close to zero.
Consider the function $\phi=\phi_\even+\phi_\odd$
obtained by setting $\phi_\even=\Phi_\even$ and $\phi_\odd=\beta^{-1}\Phi_\odd$.
The scaling factor $\beta\ne 0$ will be chosen below,
in such a way that $\phi_\even$ and $\phi_\odd$ are of comparable size.
Substituting
$$
\Phi=U_\beta\phi\defeq\phi_\even+\beta\phi_\odd
\equation(PhiDecomp)
$$
into \equ(scNS) yields the equation
$$
(\alpha\partial_t-\Delta)\phi+\thalf\gamma\LOP_s(\phi)\phi=\nsv f\,,
\equation(scNSs)
$$
where $s=\beta^2$ and
$$
\LOP_s(\phi)\psi
=\LOP(\phi_\even)\psi_\even+\LOP(\phi_\even)\psi_\odd
+\LOP(\phi_\odd)\psi_\even+s\LOP(\phi_\odd)\psi_\odd\,.
\equation(LOPsDef)
$$
Finally, we convert \equ(scNSs) to a fixed point equation
by applying the inverse of $\alpha\partial_t-\Delta$ to both sides.
Setting $g=(-\Delta)^{-1}\nsv f$,
the resulting equation is $\tilde\phi=\phi$, where
$$
\tilde\phi=g
-\thalf\gamma|\Delta|^{1/2}(\alpha\partial_t-\Delta)^{-1}\hat\phi\,,\qquad
\hat\phi\defeq|\Delta|^{-1/2}\LOP_s(\phi)\phi\,.
\equation(scNSsFix)
$$

One of the features of the equation \equ(scNSs)
is that  the time-translate of a solution is again a solution.
We eliminate this symmetry by imposing the condition $\phi_{1,1,1}=0$.
In addition, we choose $\beta=\Beth^{-1}\Phi_{-1,1,1}$,
where $\Beth$ is some fixed constant that will be specified later.
This leads to the normalization conditions
$$
A\phi\defeq\phi_{1,1,1}=0\,,\qquad
B\phi\defeq\phi_{-1,1,1}=\Beth\,.
\equation(ABNormaliz)
$$
Notice that $\beta$ enters our main equation $\tilde\phi=\phi$
only via its square $s=\beta^2$.
It is convenient to regard $s$ to be the independent parameter
and express $\gamma$ as a function of $s$.
The functions $\gamma=\gamma(s)$ and $\alpha=\alpha(s)$
are determined by the condition that $\tilde\phi$
satisfies the normalization conditions \equ(ABNormaliz).
Applying the functionals $A$ and $B$ to both sides of \equ(scNSs),
using the identities
$A\Delta=-2A$, $A\partial_t=B$, $B\Delta=-2B$, $B\partial_t=-A$,
and imposing the conditions $A\tilde\phi=0$ and $B\tilde\phi=\Beth$,
we find that
$$
\gamma=-2^{3/2}{\Beth\over B\hat\phi}\,,\qquad
\alpha=2{A\hat\phi\over B\hat\phi}\,.
\equation(gammaalphaSol)
$$
For a fixed value of $s$, define $\FF_s(\phi)=\tilde\phi$,
where $\tilde\phi$ is given by \equ(scNSsFix),
with $\gamma=\gamma(s,\phi)$ and $\alpha=\alpha(s,\phi)$
determined by \equ(gammaalphaSol).
The fixed point equation for $\FF_s$ is used
to find non-stationary time-periodic solutions of \equ(scNSs).

\demo Remark(NotNormalized)
The choice \equ(gammaalphaSol) guarantees that
$A\tilde\phi=0$ and $B\tilde\phi=\Beth$,
even if $\phi$ does not satisfy the normalization conditions \equ(ABNormaliz).
Thus, the domain of the map $\FF_s$
can include non-normalized function $\phi$.
(The same is true for the map $\FF_\gamma$ described below.)
But a fixed point of this map will be normalized by construction.

In order to determine the bifurcation point $\gamma_0$
and the corresponding frequency $\alpha_0$,
we consider the map $\FF:\phi\mapsto\tilde\phi$ given by \equ(scNSsFix) with $s=0$.
The values of $\gamma$ and $\alpha$ are again given by \equ(gammaalphaSol),
so that $A\tilde\phi=0$ and $B\tilde\phi=\Beth$.
We will show that this map $\FF$ has a fixed point
$\phi$ with the property that $\phi_{n,j,k}=0$ whenever $|n|>1$.
The values of $\gamma$ and $\alpha$ for this fixed point
define $\gamma_0$ and $\alpha_0$.

A similar map $\FF_\gamma:\phi\mapsto\tilde\phi$,
given by \equ(scNSsFix) with $s=0$,
is used to find stationary solutions of the equation \equ(scNS).
In this case, the value of $\gamma$ is being fixed,
and $\phi_\odd$ is taken to be zero.
The goal is to show that this map $\FF_\gamma$
has a fixed point $\phi_\gamma$ that is independent of time $t$.
Then $\Phi=\phi_\gamma$ is a stationary solution of \equ(scNS).

\medskip
We finish this section by computing
the derivative of the map $\FF_s$ described after \equ(gammaalphaSol).
The resulting expressions will be needed later.
Like some of the above, the following is purely formal.
A proper formulation will be given in the next section.
For simplicity, assume that $\phi$ depends on a parameter.
The derivative of a quantity $q$ with respect to this parameter
will be denoted by $\dot q$.
Define
$$
\LL_\alpha=|\Delta|^{1/2}(\alpha\partial_t-\Delta)^{-1}\,,\qquad
\LL_\alpha'=\partial_t(\alpha\partial_t-\Delta)^{-1}\,.
\equation(LLalphaDef)
$$
Using that $\FF_s(\phi)=g-\thalf\gamma\LL_\alpha\hat\phi$
with $\hat\phi=|\Delta|^{-1/2}\LOP_s(\phi)\phi$,
the parameter-derivative of $\FF_s(\phi)$ is given by
$$
D\FF_s(\phi)\dot\phi
=-\thalf\LL_\alpha\Bigl[\bigl(\dot\gamma-\gamma\dot\alpha\LL_\alpha'\bigr)\hat\phi
+\gamma\dot{\hat\phi}\Bigr]\,,\qquad
\dot{\hat\phi}=2|\Delta|^{-1/2}\LOP_s(\phi)\dot\phi\,,
\equation(DFFsPhiDotPhi)
$$
where
$$
\dot\gamma
=2^{-3/2}{\gamma^2\over\Beth} B\dot{\hat\phi}\,,\qquad
\dot\alpha
=2^{-3/2}{\alpha\gamma\over\Beth}B\dot{\hat\phi}
-2^{-1/2}{\gamma\over\Beth}A\dot{\hat\phi}\,.
\equation(DotgammaDotalpha)
$$
The above expressions for $\dot\gamma$ and $\dot\alpha$
are obtained by differentiating \equ(gammaalphaSol).

%%%%%%%%%%%%%%%%%%%%%%%%%%%%%%%%%%%%
%%%%%%%%%%%%%%%%%%%%%%%%%%%%%%%%%%%%
\section The associated contractions
%%%%%%%%%%%%%%%%%%%%%%%%%%%%%%%%%%%%
%%%%%%%%%%%%%%%%%%%%%%%%%%%%%%%%%%%%

In this section, we formulate the fixed point problems
for the maps $\FF$, $\FF_\gamma$, and $\FF_s$
in a suitable functional setting.
The goal is to reduce the problems to a point
where we can invoke the contraction mapping theorem.
After describing the necessary estimates,
we give a proof of \clm(NSHopfBif) based on these estimates.

We start by defining suitable function spaces.
Given a real number $\rho>1$, denote by $\AA$
the space of all functions $h\in\rmL^2(\torus)$
that have a finite norm $\|h\|$, where
$$
\|h\|=|h_0|
+\sum_{n\in\oldnatural_1}\sqrt{|h_n|^2+|h_{-n}|^2}\rho^n\,,\qquad
h=\sum_{n\in\oldinteger}h_n\cosi_n\,.
\equation(AADef)
$$
Here $\cosi_n$ are the trigonometric function defined in \equ(PhijkSeries).
It is straightforward to check that $\AA$
is a Banach algebra under the pointwise product of functions.
That is, $\|gh\|\le\|g\|\|h\|$ for any two functions $g,h\in\AA$.
We also identify functions on $\torus$ with $2\pi$-periodic functions on $\real$.
In this sense, a function in $\AA$
extends analytically to the strip $T(\rho)=\{z\in\complex:|\Im z|<\log\rho\}$.

Given in addition $\varrho>1$, we denote by $\buB$
the space of all function $\Phi:\torus^2\to\AA$
that admit a representation \equ(PhixyExpansion)
and have a finite norm
$$
\|\Phi\|
=\sum_{j,k\in\oldnatural_1}\|\Phi_{j,k}\|\varrho^{j+k}\,.
\equation(GGrhoepsNorm)
$$
A function $(x,y)\mapsto(t\mapsto\Phi(t,x,y))$ in this space
will also be identified with a function
$(t,x,y)\mapsto\Phi(t,x,y)$ on $\torus^3$,
or with a function on $\real^3$ that is $2\pi$-periodic in each argument.
In this sense, every function in $\buB$
extends analytically to $T(\rho)\times T(\varrho)^2$.

We consider $\AA$ and $\buB$ to be Banach spaces
over $\field\in\{\real,\complex\}$.
In the case $\field=\real$, the functions in these spaces
are assumed to take real values for real arguments.

\smallskip
Clearly, a function $\Phi\in\buB$
admits an expansion \equ(FreqProj) with $N=\oldnatural_0$.
The sequence of Fourier coefficients $\Phi_{n,k,j}$
converges to zero exponentially as $|n|+j+k$ tends to infinity.
If all but finitely many of these coefficients vanish,
then $\Phi$ is called a Fourier polynomial.
The equation \equ(FreqProj) with $N\subset\oldnatural_0$ non-empty
defines a continuous projection $\mean_N$ on $\buB$
whose operator norm is $1$.
Using Fourier series, it is straightforward
to see that the equation \equ(LLalphaDef)
defines two bounded linear operators $\LL_\alpha$ and $\LL_\alpha'$
on $\buB$, for every $\alpha\in\complex$.
The operator $\LL_\alpha$ is in fact compact.
Specific estimates will be given in Section 4.
The following will be proved in Section 4 as well.

\claim Proposition(LOPBound)
If $\Phi$ and $\phi$ belong to $\buB$,
then so does $|\Delta|^{-1/2}\LOP(\Phi)\phi$, and
$$
\bigl\||\Delta|^{-1/2}\LOP(\Phi)\phi\bigr\|
\le\bigl\||\Delta|^{-1/2}\Phi\bigr\|\|\phi\|
+\|\Phi\|\bigl\||\Delta|^{-1/2}\phi\bigr\|\,.
\equation(LOPBound)
$$

This estimate implies e.g.~that the transformation $\phi\mapsto\tilde\phi$,
given by \equ(scNSsFix) for fixed values of $s$, $\gamma$ and $\alpha$,
is well-defined and compact as a map from $\buB$ to $\buB$.

As is common in computer-assisted proofs,
we reformulate the fixed point equation for the map $\phi\mapsto\tilde\phi$
as a fixed point problem for an associated quasi-Newton map.
Since we need three distinct versions of this map,
let us first describe a more general setting.

\medskip
Let $\FF:\DD\to\BB$ be a $\rmC^1$ map
defined on an open domain $\DD$ in a Banach space $\BB$.
Let $h\mapsto\varphi+Lh$ be a continuous affine map on $\BB$.
We define quasi-Newton map $\NN$ for $(\DD,\FF,\varphi,L)$
by setting
$$
\NN(h)=\FF(\varphi+Lh)-\varphi+(\id-L)h\,.
\equation(quasiNewtonMap)
$$
The domain of $\NN$ is defined to be the set of of all $h\in\BB$
with the property that $\varphi+Lh\in\DD$.
Notice that, if $h$ is a fixed point of $\NN$,
then $\varphi+Lh$ is a fixed point of $\FF$.
In our applications, $\varphi$ is an approximate fixed point of $\FF$
and $L$ is an approximate inverse of $\id-D\FF(\varphi)$.

The following is an immediate consequence of the contraction mapping theorem.

\claim Proposition(ContrMappingThm)
Let $\FF:\DD\to\BB$ be a $\rmC^1$ map
defined on an open domain in a Banach space $\BB$.
Let $h\mapsto\varphi+Lh$ be a continuous affine map on $\BB$.
Assume that the quasi-Newton map \equ(quasiNewtonMap)
includes a non-empty ball $B_\delta=\{h\in\BB: \|h\|<\delta\}$ in its domain,
and that
$$
\|\NN(0)\|<\eps\,,\qquad\|D\NN(h)\|<K\,,\qquad h\in B_\delta\,,
\equation(ContrMappingThm)
$$
where $\eps,K$ are positive real numbers that satisfy $\eps+K\delta<\delta$.
Then $\FF$ has a fixed point in $\varphi+LB_\delta$.
If $L$ is invertible, then this fixed point is unique in $\varphi+LB_\delta$.

In our applications below, $\BB$ is always a subspace of $\buB$.
The domain parameter $\rho$ and the constant $\Beth$
that appears in the normalization condition \equ(ABNormaliz)
are chosen to have the fixed values
$$
\rho=2^5\,,\qquad\Beth=2^{-12}\,.
\equation(ParamValues)
$$
The domain parameter $\varrho$ is defined implicitly in our proofs.
That is, the lemmas below
hold for $\varrho>1$ sufficiently close to $1$.

Consider first the problem of determining
the bifurcation point $\gamma_0$ and the associated frequency $\alpha_0$.
Let $\BB=\mean_{\sss\{0,1\}}\buB$ over $\real$.
For every $\delta>0$ define $B_\delta=\{h\in\BB: \|h\|<\delta\}$.
Let $s=0$, and denote by $\DD$ the set of all functions $\phi\in\BB$
with the property that $B\hat\phi\ne 0$.
Define $\FF:\DD\to\BB$ to be the map $\phi\mapsto\tilde\phi$ given by \equ(scNSsFix),
with $\gamma=\gamma(\phi)$ and $\alpha=\alpha(\phi)$
defined by the equation \equ(gammaalphaSol).
Clearly, $\FF$ is not only $\rmC^1$ but real analytic on $\DD$.

\claim Lemma(BifPoint)
With $\FF$ as described above,
there exists an affine isomorphism $h\mapsto\varphi+L_1h$ of $\BB$
and real numbers $\eps,\delta,K>0$ satisfying $\eps+K\delta<\delta$,
such that the following holds.
The quasi-Newton map $\NN$ associated with $(\BB,\FF,\varphi,L_1)$
includes the ball $B_\delta$ in its domain
and satisfies the bounds \equ(ContrMappingThm).
The domain of $\FF$ includes the ball in $\BB$
of radius $r=\delta\|L_1\|$, centered at $\varphi$.
For every function $\phi$ in this ball,
$\gamma(\phi)=83.1733117\ldots$ and $\alpha(\phi)=4.66592275\ldots$.

Our proof of this lemma is computer-assisted
and will be described in Section 4.

By \clm(ContrMappingThm),
the map $\FF$ has a unique fixed point $\phi^\ast\in\varphi+L_1 B_\delta$.
We define $\gamma_0=\gamma(\phi^\ast)$ and $\alpha_0=\alpha(\phi^\ast)$.

Our next goal is to construct
a branch of periodic solutions for the equation \equ(scNS).
Consider $\BB=\buB$ over $\field\in\{\real,\complex\}$.
By continuity, there exists an open ball $\JJ_0\subset\field$ centered at the origin,
and an open neighborhood $\DD$ of $\phi^\ast$ in $\BB$, such that
$B\hat\phi=B|\Delta|^{-1/2}\LOP_s(\phi)\phi$ is nonzero for all $s\in\JJ_0$
and all $\phi\in\DD$.
For every $s\in\JJ_0$,
define $\FF_s:\DD\to\BB$ to be the map $\phi\mapsto\tilde\phi$
given by \equ(scNSsFix), with $\gamma=\gamma(s,\phi)$
and $\alpha=\alpha(s,\phi)$ defined by the equation \equ(gammaalphaSol).

\claim Lemma(PerSolutions)
Let $\field=\real$.
There exists a isomorphism $L$ of $\buB$
such that the following holds.
If $\NN_0$ denotes the the quasi-Newton map
associated with $(\DD,\FF_0,\phi^\ast,L)$,
then the derivative $D\NN_0(0)$ of $\NN_0$
at the origin is a contraction.

Our proof of this lemma is computer-assisted
and will be described in Section 4.
As a consequence we have the following.

\claim Corollary(PerBranch)
Consider $\field=\complex$.
There exists an open disk $\JJ\subset\complex$,
centered at the origin,
and an analytic curve $s\mapsto\phi_s$ on $\JJ$ with values in $\DD$,
such that $\FF_s(\phi_s)=\phi_s$ for all $s\in\JJ$.
If $s$ belongs to the real interval $\JJ\cap\real$, then $\phi_s$ is real.
Furthermore, $\phi_0=\phi^\ast$.

\proof
Consider still $\field=\complex$.
For $s\in\II_0$, the derivative of $\NN_s$ on its domain is given by
$$
D\NN_s(h)=D\FF_s(\phi^\ast+Lh)L+\id-L\,.
\equation(DNNsh)
$$
Assume that some function $\psi\in\buB$ satisfies
$D\FF_0(\phi^\ast)\psi=\psi$.
We may assume that $\psi$ takes real values for real arguments.
A straightforward computation shows that
$D\NN_0(0)L^{-1}\psi=L^{-1}\psi$.
Since $D\NN_0(0)$ is a contraction in the real setting,
by \clm(PerSolutions), this implies that $\psi=0$.
So the operator $D\FF_0(\phi^\ast)$ does not have an eigenvalue $1$.
This operator is compact,
since it is the composition of a bounded linear operator
with the compact operator $\LL_\alpha$.
Thus, $D\FF_0(\phi^\ast)$ has no spectrum at $1$.
By the implicit function theorem,
there exists a complex open ball $\JJ$, centered at the origin,
such that the fixed point equation $\FF_s(\phi)=\phi$
has a solution $\phi=\phi_s$ for all $s\in\JJ$.
Furthermore, the curve $s\mapsto\phi_s$ is analytic,
passes through $\phi^\ast$ at $s=0$,
and there is a unique curve with this property.
By uniqueness, we also have $\ov{\phi_{\bar s}}=\phi_s$ for all $s\in\JJ$,
so $\phi_s$ is real for real values of $s\in\JJ$.
\qed

A branch of stationary periodic solutions for \equ(scNS)
is obtained similarly.
Consider $\BB=\mean_{\sss\{0\}}\buB$ over $\field\in\{\real,\complex\}$.
For every $\gamma\in\field$,
define $\FF_\gamma:\BB\to\BB$ to be the map $\phi\mapsto\tilde\phi$
given by \equ(scNSsFix), with $s=\alpha=0$.
Notice that $\phi^\ast_\even$ is a fixed point of $\FF_{\gamma_0}$.

\claim Lemma(StatSolutions)
Let $\field=\real$. There exists an isomorphism $L_0$ of $\BB$
such that the following holds.
If $\NN_{\gamma_0}$ denotes the the quasi-Newton map
associated with $(\BB,\FF_{\gamma_0},\phi^\ast_\even,L_0)$,
then the derivative $D\NN_{\gamma_0}(0)$ of $\NN_{\gamma_0}$
at the origin is a contraction.

Our proof of this lemma is computer-assisted
and will be described in Section 4.
As a consequence we have the following.

\claim Corollary(StatBranch)
Consider $\field=\complex$.
There exists an open disk $\II\subset\complex$,
centered at $\gamma_0$,
and an analytic curve $\gamma\mapsto\phi_\gamma$ on $\II$ with values in $\BB$,
such that $\FF_\gamma(\phi_\gamma)=\phi_\gamma$ for all $\gamma\in\II$.
If $\gamma$ belongs to the real interval $\II\cap\real$, then $\phi_\gamma$ is real.
Furthermore, $\phi_{\gamma_0}=\phi^\ast_\even$.

The proof of this corollary is analogous
to the proof of \clm(PerBranch).

We note that the disk $\II\ni\gamma_0$
is disjoint from the disk $\JJ\ni 0$ described in \clm(PerBranch).
So there is no ambiguity in using the notation
$\gamma\mapsto\phi_\gamma$ and $s\mapsto\phi_s$
for the curve of stationary and periodic solutions,
respectively, of the equation \equ(scNSs),

\smallskip
Based on the results stated in this section, we can now give a

\proofof(NSHopfBif)
As described in the preceding sections,
the curve $\gamma\mapsto\phi_\gamma$ for $\gamma\in\II$
yields a curve $\gamma\mapsto u_\gamma$ of stationary solutions of
the equation \equ(NS), where $u_\gamma=\nsv^{-1}\phi_\gamma$.
By our choice of function spaces,
the function $(\gamma,x,y)\mapsto u_\gamma(x,y)$
is real analytic on $I\times\torus^2$, where $I=\II\cap\real$.

Similarly, the curve $s\mapsto\phi_s$ for $s\in\JJ$
defines a family of of non-stationary periodic solutions
for \equ(NS), with $\gamma=\gamma_s$ and $\alpha=\alpha_s$
determined via the equation \equ(gammaalphaSol).
To be more precise,
the even frequency part $\phi_{s,\even}$ of $\phi_s$
determines a vector field $u_{s,\even}=\nsv^{-1}\phi_{s,\even}$,
and the odd frequency part $\phi_{s,\odd}$
determines a vector field $u_{s,\odd}=\nsv^{-1}\phi_{s,\odd}$.
If $\beta$ is a complex number such that $s=\beta^2\in\JJ$,
then $u=u_{s,\even}+\beta u_{s,\odd}$
is a periodic solution of \equ(NS), with $\gamma=\gamma_s$ and $\alpha=\alpha_s$.
Here, we have used the decomposition \equ(PhiDecomp).
By our choice of function spaces, the functions
$(s,t,x,y)\mapsto u_{s,\even}(t,x,y)$ and $(s,t,x,y)\mapsto u_{s,\odd}(t,x,y)$
are real analytic on $J\times\torus^3$, where $J=\JJ\cap\real$.
Clearly, $\partial_t u_{0,\odd}(t,\bdot\,,\bdot)\ne 0$,
due to the normalization condition $\phi_{-1,1,1}=\Beth$ imposed in \equ(ABNormaliz).
And by construction, we have $u=u_{\gamma_0}$ for $s=0$.
\qed

%%%%%%%%%%%%%%%%%%%%%%%%%%%%
%%%%%%%%%%%%%%%%%%%%%%%%%%%%
\section Remaining estimates
%%%%%%%%%%%%%%%%%%%%%%%%%%%%
%%%%%%%%%%%%%%%%%%%%%%%%%%%%

What remains to be proved are Lemmas
\clmno(BifPoint), \clmno(PerSolutions), and \clmno(StatSolutions).
Our method used in the proof of \clm(BifPoint)
can be considered perturbation theory
about the approximate fixed point $\varphi$ of $\FF$.
The function $\varphi$ is a Fourier polynomial with over $11000$ nonzero coefficients,
so a large number of estimates are involved.

We start by describing bounds on the bilinear function $\LOP$
and on the linear operators $\LL_\alpha$ and $\LL_\alpha'$.
These are the basic building blocks for our transformations
$\FF$, $\FF_s$, and $\FF_\gamma$.
The ``mechanical'' part of these estimates will be described in Subsection 4.4.

%%%%%%%%%%%%%%%%%%%%%%%%%%%%%%%%%%%%%%%%%%%%%%%%%%%%%%%%%%%%%%%%%%%%%%%%%%%
\subsection The bilinear form $\boldmath\LOP$ and a proof of \clm(LOPBound)
%%%%%%%%%%%%%%%%%%%%%%%%%%%%%%%%%%%%%%%%%%%%%%%%%%%%%%%%%%%%%%%%%%%%%%%%%%%

Consider the bilinear form $\LOP$ defined by \equ(LOPDef).
Using the identity \equ(curlInv), we have
$$
\eqalign{
\LOP(\Phi)\phi
&=(\nabla\Phi)\cdot\iso\nabla\Delta^{-1}\phi
+(\nabla\phi)\cdot\iso\nabla\Delta^{-1}\Phi\cr
&=\bigl[(\partial_x\Phi)\Delta^{-1}\partial_y\phi
-(\partial_y\Phi)\Delta^{-1}\partial_x\phi\bigr]
-\bigl[(\Delta^{-1}\partial_x\Phi)\partial_y\phi
-(\Delta^{-1}\partial_y\Phi)\partial_x\phi\bigr]\,.\cr}
\equation(LOPPhiphi)
$$
In order to obtain accurate estimates,
it is useful to have explicit expressions for $\LOP(\Phi)\phi$
in terms of the Fourier coefficients of $\Phi$ and $\phi$.
Given that $\LOP$ is bilinear,
and that the identity \equ(LOPPhiphi) holds pointwise in $t$,
it suffices to compute $\LOP(\Phi)\phi$
for the time-independent monomials
$$
\Phi=\sin_J\times\sin_K\,,\qquad
\phi=\sin_j\times\sin_k\,,
\equation(PhiphiModes)
$$
with $J,K,j,k>0$.
A straightforward computation shows that
$$
\eqalign{
\LOP(\Phi)\phi&=\Theta(Jk+jK)\bigl[
\sin_{J+j}\times\sin_{K-k}-\sin_{J-j}\times\sin_{K+k}
\bigr]\cr
&\quad+\Theta(Jk-jK)\bigl[
\sin_{J+j}\times\sin_{K+k}-\sin_{J-j}\times\sin_{K-k}
\bigr]\,,\cr}
\equation(LOPPhiphiOne)
$$
with $\Theta$ as defined below.
As a result we have
$$
|\Delta|^{-1/2}\LOP(\Phi)\phi
=\sum_{\sigma,\tau=\pm 1}N_{\sigma,\tau}\sin_{\sigma J+j}\times\sin_{\tau K+k}\,,
\equation(InvHalfLapLOPDecomp)
$$
where
$$
N_{\sigma,\tau}
=\Theta
{\sigma Jk-\tau Kj\over\sqrt{(\sigma J+j)^2+(\tau K+k)^2}}\,,\qquad
\Theta=\quarter\biggl({1\over J^2+K^2}-{1\over j^2+k^2}\biggr)\,.
\equation(NpmpmDef)
$$

\proofof(LOPBound)
Using the Cauchy-Schwarz inequality in $\real^2$, we find that
$$
|N_{\sigma,\tau}|
=|\Theta|{|(\sigma J+j)k-(\tau K+k)j|\over\sqrt{(\sigma J+j)^2+(\tau K+k)^2}}
\le|\Theta|\sqrt{j^2+k^2}\,.
\equation(NppBound)
$$
Since the absolute value of $N_{\sigma,\tau}$ is invariant
under an exchange of $(j,k)$ and $(J,K)$, this implies that
$$
|N_{\sigma,\tau}|\le{1/4\over\sqrt{j^2+k^2}}\vee{1/4\over\sqrt{J^2+K^2}}\,,
\equation(NpmpmBound)
$$
where $a\vee b=\max(a,b)$ for $a,b\in\real$.
As a result, we obtain the bound
$$
\bigl\||\Delta|^{-1/2}\LOP(\Phi)\phi\bigr\|
\le\bigl\||\Delta|^{-1/2}\Phi\bigr\|_{\varrho,\epsilon}\|\phi\|
+\|\Phi\|\bigl\||\Delta|^{-1/2}\phi\bigr\|\,.
\equation(LOPPhiphiBound)
$$
Using the nature of the norm \equ(GGrhoepsNorm),
and the fact that $\AA$ is a Banach algebra
for the pointwise product of functions,
this bound extends by bilinearity to arbitrary
functions $\Phi,\phi\in\buB$.
\qed

We note that the bound \equ(LOPPhiphiBound)
exploits the cancellations that lead to the expression \equ(LOPPhiphiOne).
A more straightforward estimate loses a factor of $2$
with respect to \equ(LOPPhiphiBound).
But it is not just this factor of $2$ that counts for us.
The expressions \equ(NpmpmDef) for the coefficients $N_{\sigma,\tau}$
and the bounds \equ(NpmpmBound) are used in our computations
and error estimates.
The expression on the right hand side of \equ(NpmpmBound)
is a decreasing function of the wavenumbers $j,k,J,K$,
so it can be used to estimate $\LOP(\Phi)\phi$
when $\Phi$ and/or $\phi$ are ``tails'' of Fourier series.

%%%%%%%%%%%%%%%%%%%%%%%%%%%%%%%%%%%%%%%%%%%%%%%%%%%%%%%%%%%%%%%%%%%%%%%%%%%%%%%
\subsection The linear operators $\boldmath\LL_\alpha$ and $\boldmath\LL_\alpha'$
%%%%%%%%%%%%%%%%%%%%%%%%%%%%%%%%%%%%%%%%%%%%%%%%%%%%%%%%%%%%%%%%%%%%%%%%%%%%%%%

Consider the linear operators $\LL_\alpha$
and $\LL_\alpha'$ defined in \equ(LLalphaDef), with $\alpha$ real.
A straightforward computation shows that
$$
\psi_{n,j,k}
=\sqrt{j^2+k^2}\,{(j^2+k^2)\phi_{n,j,k}-\alpha n\phi_{-n,j,k}
\over(j^2+k^2)^2+\alpha^2 n^2}\,,\qquad
\psi=\LL_\alpha\phi\,.
\equation(LLalphaphinjk)
$$
Using the Cauchy-Schwarz inequality in $\real^2$,
this yields the estimate
$$
\sqrt{|\psi_{n,j,k}|^2+|\psi_{-n,j,k}|^2}
\le C_{n,j,k}\sqrt{|\phi_{n,j,k}|^2+|\phi_{-n,j,k}|^2}\,,
\equation(LLalphaBound)
$$
with
$$
C_{n,j,k}=\sqrt{j^2+k^2\over(j^2+k^2)^2+\alpha^2 n^2}
\le{1\over\sqrt{2|\alpha n|}}\wedge {1\over\sqrt{j^2+k^2}}
\equation(LLalphaBoundConst)
$$
for $n\ne 0$, where $a\wedge b=\min(a,b)$ for $a,b\in\real$.
The last bound in \equ(LLalphaBoundConst)
is a decreasing function of $|n|,j,k$
and can be used to estimate $\LL_\alpha\phi$
when $\phi$ is the tail of a Fourier series.

For the operator $\LL_\alpha'$ we have
$$
\psi_{n,j,k}
=n\,{(j^2+k^2)\phi_{-n,j,k}+\alpha n\phi_{n,j,k}
\over(j^2+k^2)^2+\alpha^2 n^2}\,,\qquad
\psi=\LL_\alpha'\phi\,.
\equation(LLalphaPrimenjk)
$$
A bound analogous to \equ(LLalphaBound)
holds for $\psi=\LL_\alpha'\phi$, with
$$
C_{n,j,k}=\sqrt{n^2\over(j^2+k^2)^2+\alpha^2 n^2}\,.
\equation(LLalphaPrimeBoundConst)
$$
As can be seen from \equ(DFFsPhiDotPhi),
this bound is needed only for $n=\pm 1$,
since these are the only nonzero frequencies
of the function $\hat\phi=|\Delta|^{-1/2}\LOP_0(\phi)\phi$
with $\phi\in\mean_{\sss\{0,1\}}\buB$.
And for fixed $n$, the right hand side of \equ(LLalphaPrimeBoundConst)
is decreasing in $j$ and $k$.

%%%%%%%%%%%%%%%%%%%%%%%%%%%%%%%%%%%%%
\subsection Estimating operator norms
%%%%%%%%%%%%%%%%%%%%%%%%%%%%%%%%%%%%%

Recall that a function $\phi\in\buB$ admits a Fourier expansion
$$
\phi=\sum_{n\in\oldinteger}\;\sum_{j,k\in\oldnatural_1}\phi_{n,j,k}\theta_{n,j,k}\,,\qquad
\theta_{n,j,k}\defeq\cosi_n\times\sin_j\times\sin_k\,,
\equation(thetanjkDef)
$$
and that the norm of $\phi$ is given by
$$
\|\phi\|=\sum_{j,k\in\oldnatural_1}
\biggl[|\phi_{0,j,k}|+\sum_{n\in\oldnatural_1}
\sqrt{|\phi_{n,j,k}|^2+|\phi_{-n,j,k}|^2}\,\rho^n\biggr]\varrho^{j+k}\,.
\equation(GGrhoepsNormAgain)
$$
Let now $n\ge 0$.
A linear combination $c_{\sss+}\theta_{n,j,k}+c_{\sss-}\theta_{-n,j,k}$
will be referred to as a mode with frequency $n$ and wavenumbers $(j,k)$
or as a mode of type $(n,j,k)$.
We assume of course that $c_{\sss-}=0$ when $n=0$.
Since \equ(GGrhoepsNormAgain) is a weighted $\ell^1$ norm,
except for the $\ell^2$ norm used for modes,
we have a simple expression for the operator norm
of a continuous linear operator $\LL:\buB\to\buB$, namely
$$
\|\LL\|=\sup_{j,k\in\oldnatural_1}\;
\sup_{n\in\oldnatural_0}\;\sup_u\|\LL u\|/\|u\|\,,
\equation(OpNorm)
$$
where the third supremum is over all nonzero modes
$u$ of type $(n,j,k)$.

Let now $n,j,k\ge 1$ be fixed.
In computation where $\LL\theta_{\pm n,j,k}$ is known explicitly,
we use the following estimate.
Denote by $\LL_{n,j,k}$ the restriction of $\LL$
to the subspace spanned by the two functions $\theta_{\pm n,j,k}$.
For $q\ge 1$ define
$$
\|\LL_{n,j,k}\|_q
=\sup_{0\le p<q}\|\LL v_p\|\,,\qquad
v_p=\cos\Bigl({\pi p\over q}\Bigr){\theta_{n,j,k}\over\rho^n\varrho^{j+k}}
+\sin\Bigl({\pi p\over q}\Bigr){\theta_{-n,j,k}\over\rho^n\varrho^{j+k}}\,.
\equation(mNorm)
$$
Since every unit vector in the span of $\theta_{\pm n,j,k}$
lies within a distance less than ${\pi\over q}$
of one of the vectors $v_p$ or its negative, we have
$\|\LL_{n,j,k}\|\le\|\LL_{n,j,k}\|_q+{\pi\over q}\|\LL_{n,j,k}\|$.
Thus
$$
\|\LL_{n,j,k}\|\le{q\over q-\pi}\|\LL_{n,j,k}\|_m\,,\qquad q\ge 4\,.
\equation(LLmodeNormBound)
$$

Consider now the operator
$D\FF_s(\phi)$ described in \equ(DFFsPhiDotPhi),
with $\phi\in\mean_{\sss\{0,1\}}\buB$ fixed.
If $\dot\phi=u_n$ is a nonzero mode with frequency $n\ge 3$,
then $\dot{\hat\phi}=2|\Delta|^{-1/2}\LOP_0(\phi)\dot\phi$
belongs to $\mean_N\buB$ with $N=\{n-1,n,n+1\}$.
Thus, we have $\dot\gamma=\dot\alpha=0$, and
$$
D\FF_0(\phi)u_n
=-\gamma\LL_\alpha|\Delta|^{-1/2}\LOP_0(\phi)u_n\,.
\equation(DFFophiu)
$$
Due to the factor $\LL_\alpha$ in this equation,
if $u_n=c_{\sss+}\theta_{n,j,k}+c_{\sss-}\theta_{-n,j,k}$
with $(j,k)$ and $c_{\sss\pm}$ fixed,
then the ratios
$$
\|D\FF_0(\phi)u_n\|/\|u_n\|
\equation(DFFophiunRatios)
$$
are decreasing in $n$ for $n\ge 3$.
And the limit as $n\to\infty$ of this ratio is zero.

So for the operator $\LL=D\FF_0(\phi)$,
the supremum over $n\in\oldnatural_0$ in \equ(OpNorm)
reduces to a maximum over finitely many terms.
The same holds for the operator $\LL=D\NN_0(0)=D\FF_0(\phi^\ast)L+\id-L$
that is described in \clm(PerSolutions).
This is a consequence of the following choice.

\demo Remark(LMatrix)
The operator $L$ chosen in \clm(PerSolutions)
is a ``matrix perturbation'' of the identity,
in the sense that $L\theta_{n,j,k}=\theta_{n,j,k}$
for all but finitely many indices $(n,j,k)$.
The same is true for the operators $L_1$
and $L_0$ chosen in \clm(BifPoint) and \clm(StatSolutions),
respectively.

%%%%%%%%%%%%%%%%%%%%%%%%%%%%%%
\subsection Computer estimates
%%%%%%%%%%%%%%%%%%%%%%%%%%%%%%

Lemmas \clmno(BifPoint), \clmno(StatSolutions),
and \clmno(PerSolutions) assert the existence of certain objects
that satisfy a set of strict inequalities.
The goal here is to construct these objects,
and to verify the necessary inequalities
by combining the estimates that have been described so far.

The above-mentioned ``objects'' are real numbers,
real Fourier polynomials, and linear operators
that are finite-rank perturbations of the identity.
They are obtained via purely numerical computations.
Verifying the necessary inequalities is largely an organizational task,
once everything else has been set up properly.
Roughly speaking, the procedure follows that of a well-designed numerical program,
but instead of truncation Fourier series and ignoring rounding errors,
we determine rigorous enclosures at every step along the computation.
This part of the proof is written in the programming language Ada [\rAda].
The following is meant to be a rough guide for the reader who
wishes to check the correctness of our programs.
The complete details can be found in [\rProgs].

\smallskip
An enclosure for a function $\phi\in\buB$
is a set in $\buB$ that includes $\phi$
and is defined in terms of (bounds on) a Fourier polynomial
and finitely many error terms.
We define such sets hierarchically,
by first defining enclosures for elements in simpler spaces.
In this context, a ``bound'' on a map $f:\XX\to\YY$
is a function $F$ that assigns to a set $X\subset\XX$
of a given type ({\tt Xtype}) a set $Y\subset\YY$
of a given type ({\tt Ytype}), in such a way that
$y=f(x)$ belongs to $Y$ for all $x\in X$.
In Ada, such a bound $F$ can be implemented by defining
a {\tt procedure F(X:{\gapii}in Xtype; Y:{\gapii}out Ytype)}.

Our most basic enclosures are specified by pairs {\tt S=(S.C,S.R)},
where {\tt S.C} is a representable real number ({\tt Rep})
and {\tt S.R} a nonnegative representable real number ({\tt Radius}).
Given a Banach algebra $\XX$ with unit ${\bf 1}$,
such a pair {\tt S} defines a ball in $\XX$ which we denote by
$\langle{\tt S},\XX\rangle=\{x\in\XX:\|x-({\tt S.C}){\bf 1}\|\le{\tt S.R}\}$.

When $\XX=\real$,
then the data type described above is called {\tt Ball}.
Bounds on some standard functions involving the type {\tt Ball}
are defined in the package {\tt Flts\_Std\_Balls}.
Other basic functions are covered in the packages {\tt Vectors} and {\tt Matrices}.
Bounds of this type have been used in many computer-assisted proofs;
so we focus here on the more problem-specific aspects of our programs.

Consider now the space $\AA$
for a fixed domain radius $\varrho>1$ of type {\tt Radius}.
As mentioned before \dem(LMatrix),
we only need to consider Fourier polynomials in $\AA$.
Our enclosures for such polynomials are defined by
an {\tt array(-I$_{\rmc}${\gapi}..{\gapi}I$_{\rmc}$) of Ball}.
This data type is named {\tt NSPoly},
and the enclosure associated with data {\tt P} of this type is
$$
\langle{\tt P},\AA\rangle
\defeq\sum_{i=-I_{\rmc}}^{I_{\rmc}}
\bigl\langle{\tt P(i)},\real\bigr\rangle\cosi_{\nu(i)}\,,
\equation(NSPolyEnclosure)
$$
where $\nu$ is an increasing index function with the property that $\nu(-i)=-\nu(i)$.
The type {\tt NSPoly} is defined in the package {\tt NSP},
which also implements bounds on some basic operations
for Fourier polynomials in $\AA$.
Among the arguments to {\tt NSP} is a nonnegative integer $n$
(named {\tt NN}).
Our proof of \clm(StatSolutions) and \clm(BifPoint) uses
$I_c=n=0$ and $I_c=n=1$, respectively, and $\nu(i)=i$.
Values $n\ge 2$ are uses when estimating the norm of $\LL u$
for the operator $\LL=D\NN_0(0)$, with $u$ a mode of frequency $n$.
In this case, $\nu$ takes values in $\{-n,n\}$ or $\{-n-1,-n,-n+1,0,n-1,n,n+1\}$,
depending on whether $n$ is odd or even.
(The value $\nu=0$ is being used only for $n=2$.)
The package {\tt NSP} also defines a data type {\tt NSErr}
as an {\tt array(0{\gapi}..{\gapi}I$_{\rm c}$) of Radius}.
This type will be used below.

Given in addition a positive number $\varrho\ge 1$ of type {\tt Radius},
our enclosures for functions in $\buB$
are defined by pairs {\tt(F.C,F.E)},
where {\tt F.C} is an
{\tt array(1{\gapi}..{\gapi}J$_{\rmc}$,1{\gapi}..{\gapi}K$_{\rmc}$) of NSPoly}
and {\tt F.E} is an
{\tt array(1{\gapi}..{\gapi}J$_{\rme}$,1{\gapi}..{\gapi}K$_{\rme}$) of NSErr};
all for a fixed value of the parameter {\tt NN}.
This data type is named {\tt Fourier3},
and the enclosure associated with {\tt F=(F.C,F.E)} is
$$
\langle{\tt F},\buB\rangle
\defeq\sum_{j=1}^{J_\rmc}\sum_{k=1}^{K_\rmc}\bigl\langle{\tt F.C(j,k)},\AA\bigr\rangle
\times\sin_j\times\sin_k
+\sum_{J=1}^{J_\rme}\sum_{K=1}^{K_\rme}H_{\sss J,K}({\tt F.E(J,K)})\,.
\equation(FouThreeEnclosure)
$$
Here, $H_{\sss J,K}({\tt E})$ denotes the set of all functions
$\phi=\sum_{i=0}^{I_{\rmc}}\phi^i$ with $\|\phi^i\|\le{\tt E(i)}$,
where $\phi^i$ can be any function in $\buB$
whose coefficients $\phi^i_{n,j,k}$ vanish unless
$j\ge J$, $k\ge K$, and $|n|=\nu(i)$.

The type {\tt Fourier3} and bounds on some standard functions
involving this type are defined in the child package {\tt NSP.Fouriers}.
This package is a modified version of the
package {\tt Fouriers2} that was used earlier in [\rAKxii,\rAKxviii,\rAGK].
The procedure {\tt Prod} is now a bound on the bilinear map $|\Delta|^{-1/2}\LOP_0$.
The error estimates used in {\tt Prod} are based on the inequality \equ(NpmpmBound).
The package {\tt NSP.Fouriers}
also includes bounds {\tt InvLinear} and {\tt DtInvLinear}
on the linear operators $\LL_\alpha$ and $\LL_\alpha'$, respectively.
These bounds use the estimates described in Subsection 4.3.

As far as the proof of \clm(BifPoint) is concerned,
it suffices now to compose existing bounds
to obtain a bound on the map $\FF$ and its derivative $D\FF$.
This is done by the procedures {\tt GMap} and {\tt DGMap} in {\tt Hopf.Fix}.
Here we use enclosures of type {\tt NN=1}.

The type of quasi-Newton map $\NN$ defined by \equ(quasiNewtonMap)
has been used in several computer-assisted proof before.
So the process of constructing a bound on $\NN$
from a bound on $\FF$ has been automated
in the generic packages {\tt Linear} and {\tt Linear.Contr}.
(Changes compared to earlier versions are mentioned in the program text.)
This includes the computation of an approximate inverse $L_1$
for the operator $\id-D\FF(\varphi)$.
A bound on $\NN$ is defined (in essence) by the procedure {\tt Linear.Contr.Contr},
instantiated with {\tt Map => GMap}.
And a bound on $D\NN$ is defined by {\tt Linear.Contr.DContr},
with {\tt DMap => DGMap}.
Bounds on operator norms are obtained via {\tt Linear.OpNorm}.
Another problem-dependent ingredient in these procedures,
besides {\tt Map} and {\tt DMap}, are data of type {\tt Modes}.
These data are constructed by the procedure {\tt Make}
in the package {\tt Hopf}.
They define a splitting of the given space $\BB$ into a finite direct sum.
For details on how such a splitting is defined and used we refer to [\rAKxix].

If the parameter {\tt NN} has the value $0$,
then the procedures {\tt GMap} and {\tt DGMap} define
bounds on the map $\FF_\gamma$ and its derivative, respectively.
The operator $L_0$ used in \clm(StatSolutions)
has the property that $M_0=L_0-\id$
satisfies $M_0=P_0M_0P_0$, where $P_0=\mean_{\sss\{0\}}\proj_{m_0}$
for some positive integer $m_0$.
Here, and in what follows, $\proj_m$ denotes
the canonical projection in $\buB$
with the property that $\proj_m\phi$ is obtained
from $\phi$ by restricting the second sum in \equ(thetanjkDef)
to wavenumbers $j,k\le m$.

If {\tt NN} has a value $n\ge 2$, then the procedure {\tt DGMap}
defines a bound on the map $(\phi,\psi)\mapsto D\FF_0(\phi)\psi$,
restricted to the subspace $\mean_{\sss\{0,1\}}\buB\times\mean_{\sss\{n\}}\buB$.
The linear operator $L$ that is used in \clm(PerSolutions)
admits a decomposition $L=\id+M_1+M_2+\ldots+M_N$ of the following type.
After choosing a suitable sequence $n\mapsto m_n$ of positive integers,
we set $M_n=P_n(L-\id)P_n$,
where $P_1=\mean_{\sss\{0,1\}}\proj_{m_1}$
and $P_n=\mean_{\sss\{n\}}\proj_{m_n}$ for $n=2,3,\ldots,N$.
This structure of $L$ simplifies the use of \equ(OpNorm)
for estimating the norm of $\LL=D\NN_0(0)$.
Furthermore, to check that $L$ is invertible,
it suffices to verify that $\id+M_n$ is invertible
on the finite-dimensional subspace $P_n\buB$,
for each positive $n\le N$.

The linear operator $L_1$ that is used in \clm(BifPoint)
is of the form $L_1=\id+M_1$ with $M_1$ as described above.

All the steps required in the proofs
of Lemmas \clmno(BifPoint), \clmno(StatSolutions), and \clmno(PerSolutions)
are organized in the main program {\tt Check}.
As $n$ ranges from $0$ to $N=305$,
this program defines the parameters that are used in the proof
for {\tt NN} $=n$, instantiates the necessary packages,
computes the appropriate matrix $M_n$,
verifies that $\id+M_n$ is invertible,
reads $\varphi$ from the file {\tt BP.approx},
and then calls the procedure {\tt ContrFix} from the
(instantiated version of the) package {\tt Hopf.Fix} to verify the necessary inequalities.

The representable numbers ({\tt Rep}) used in our programs
are standard [\rIEEE] extended floating-point numbers (type {\tt LLFloat}).
High precision [\rMPFR] floating-point numbers (type {\tt MPFloat})
are used as well, but not in any essential way.
Both types support controlled rounding.
{\tt Radius} is always a subtype of {\tt LLFloat}.
Our programs were run successfully on a $20$-core workstation,
using a public version of the gcc/gnat compiler [\rGnat].
For further details,
including instruction on how to compile and run our programs,
we refer to [\rProgs].

\bigskip
%%%%%%%%%%%
\references
%%%%%%%%%%%

{\ninepoint

\item{[\rHopf]} E.~Hopf,
{\sl Abzweigung einer periodischen L\"osung
von einer station\"aren L\"osung eines Differentialsystems},
Ber. Math.-Phys. Kl. Siichs. Akad. Wiss. Leipzig, {\bf 94}, 3--22 (1942).

\item{[\rSerrin]} J.~Serrin,
{\sl A Note on the Existence of Periodic Solutions of the Navier-Stokes Equations},
Arch. Rational Mech. Anal. {\bf 3} 120--122, (1959).
%%% periodic forcing

\item{[\rCrAb]} M.G.~Crandall and P.H.~Rabinowitz,
{\sl The Hopf bifurcation theorem in infinite dimensions},
Arch. Rational Mech. Anal. {\bf 67}, 53--72 (1977).
%\pdfclink{0 0 1}{online here}
%{https://www.semanticscholar.org/paper/The-Hopf-Bifurcation-Theorem-in-infinite-dimensions-Crandall-Rabinowitz/c21e703fdf69882c64cff84eb976cc16ddf7345e}

\item{[\rMarMcC]} J.~Marsden, M.~McCracken,
{\sl The Hopf bifurcation and its applications},
Springer Applied Mathematical Sciences Lecture Notes Series, Vol. 19, 1976.

\item{[\rRT]} D.~Ruelle, F.~Takens, {\sl On the Nature of Turbulence},
Commun. Math. Phys. 20, 167--192 (1971)

\item{[\rKlWe]} P.~Kloeden, R.~Wells,
{\sl An explicit example of Hopf bifurcation in fluid mechanics},
Proc. Roy. Soc. London Ser. A {\bf 390}, 293--320 (1983).
%\pdfclink{0 0 1}{online here}
%{https://royalsocietypublishing.org/doi/10.1098/rspa.1983.0133}

\item{[\rChIo]} P.~Chossat, G.~Iooss,
{\sl Primary and secondary bifurcations in the Couette-Taylor problem},
Japan J. Appl. Math. {\bf 2}, 37--68 (1985).

\item{[\rDum]} F.~Dumortier,
{\sl Techniques in the theory of local bifurcations:
blow-up, normal forms, nilpotent bifurcations, singular perturbations};
in: {\sl Bifurcations and periodic orbits of vector fields},
(D.~Schlomiuk, ed., Kluwer Acad.~Pub.)
NATO ASI Ser. C Math. Phys. Sci. {\bf 408}, 10--73 (1993).

\item{[\rChIoo]} P.~Chossat, G.~Iooss,
{\sl The Couette-Taylor problem},
Applied Mathematical Sciences, 102. Springer-Verlag, New York, 1994

\item{[\rNWYNK]} M.T.~Nakao, Y.~Watanabe, N.~Yamamoto, T.~Nishida, M.-N.~Kim,
{\sl Computer assisted proofs of bifurcating solutions for nonlinear heat convection problems},
J. Sci. Comput. {\bf 43}, 388--401 (2010).
%\pdfclink{0 0 1}{online here}{https://doi.org/10.1007/s10915-009-9303-3}
%% Oberbeck-Boussinesq equation (reduction 3d --> 2d) for The Rayleigh-Benard problem.
%% (extension to 3d is discussed at the end)
%% Consider only bifurcations of stationary solutions
%% Similarly for the referencces 9 and 12 that they cite (on the same problem)

\item{[\rAKxii]} G.~Arioli, H.~Koch,
{\sl Non-symmetric low-index solutions for a symmetric boundary value problem},
J. Differ. Equations {\bf 252}, 448--458 (2012).

\item{[\rAKxv]} G.~Arioli, H.~Koch,
{\sl Some symmetric boundary value problems and non-symmetric solutions},
J. Differ. Equations {\bf 259}, 796--816 (2015).

\item{[\rGaldi]} G.P.~Galdi,
{\sl On bifurcating time-periodic flow of
a Navier-Stokes liquid past a cylinder},
Arch. Rational Mech. Anal. {\bf 222}, 285--315 (2016).
Digital Object Identifier (DOI) 10.1007/s00205-016-1001-3

\item{[\rHJNS]} C.-H.~Hsia, C.-Y.~Jung, T.B.~Nguyen, and M.-C.~Shiu,
{\sl On time periodic solutions,
asymptotic stability and bifurcations of Navier-Stokes equations},
Numer. Math. {\bf 135}, 607--638 (2017).
%%% time-periodic forcing, and reference
%%% found this ref in [\rBBLV]

\item{[\rAKxviii]} G.~Arioli, H.~Koch,
{\sl Spectral stability for the wave equation with periodic forcing},
J. Differ. Equations {\bf 265}, 2470--2501 (2018).

\item{[\rAKxix]} G.~Arioli, H.~Koch,
{\sl Non-radial solutions for some semilinear elliptic equations on the disk},
Nonlinear Analysis {\bf 179}, 294–308 (2019).

\item{[\rNPW]} M.T.~Nakao, M.~Plum, Y.~Watanabe,
{\sl Numerical verification methods and computer-assisted proofs
for partial differential equations},
Springer Series in Computational Mathematics, Vol. 53,
Springer Singapore, 2019

\item{[\rGoSe]} J.~G\'omez-Serrano,
{\sl Computer-assisted proofs in PDE: a survey},
SeMA {\bf 76}, 459--484 (2019).

\item{[\rWiZg]} D.~Wilczak, P.~Zgliczy\'nski,
{\sl A geometric method for infinite-dimensional chaos:
Symbolic dynamics for the Kuramoto-Sivashinsky PDE on the line},
J. Differ. Equations {\bf 269}, 8509--8548 (2020).

\item{[\rBBLV]} J.~B.~van den Berg, M.~Breden, J.-P.~Lessard, L.~van Veen,
{\sl Spontaneous periodic orbits in the Navier-Stokes flow},
Preprint 2019,
%\pdfclink{0 0 1}{online here}{https://arxiv.org/abs/1902.00384}

\item{[\rAGK]} G.~Arioli, F.~Gazzola, H.~Koch,
{\sl Uniqueness and bifurcation branches for planar steady
Navier-Stokes equations under Navier boundary conditions},
Preprint 2020.

\item{[\rBLQ]} J.~B.~van den Berg, J.-P.~Lessard, E.~Queirolo,
{\sl Rigorous verification of Hopf bifurcations
via desingularization and continuation},
Preprint 2020.

\item{[\rBQ]} J.~B.~van den Berg, E.~Queirolo,
{\sl Validating Hopf bifurcation in the Kuramoto-Sivashinky PDE},
in preparation.
% talk by Elena Queirolo, see
% https://researchseminars.org/seminar/CRM-CAMP

\item{[\rProgs]} G.~Arioli, H.~Koch,
{\sl Programs and data files for the proof of Lemmas \clmno(BifPoint),
\clmno(StatSolutions), and \clmno(PerSolutions)
\pdfclink{0 0 1}{{\tt https://web.ma.utexas.edu/users/koch/papers/nshopf/}}
{https://web.ma.utexas.edu/users/koch/papers/nshopf/}

%%%%%%%%%%%%%%%%%%%%%%%%%%%%%%%%%%%%

\item{[\rAda]} Ada Reference Manual, ISO/IEC 8652:2012(E),
available e.g. at\hfil\break
\pdfclink{0 0 1}{{\tt www.ada-auth.org/arm.html}}
{http://www.ada-auth.org/arm.html}

\item{[\rGnat]}
A free-software compiler for the Ada programming language,
which is part of the GNU Compiler Collection; see
\pdfclink{0 0 1}{{\tt gnu.org/software/gnat/}}{http://gnu.org/software/gnat/}

\item{[\rIEEE]} The Institute of Electrical and Electronics Engineers, Inc.,
{\sl IEEE Standard for Binary Float\-ing--Point Arithmetic},
ANSI/IEEE Std 754--2008.

\item{[\rMPFR]} The MPFR library for multiple-precision floating-point computations
with correct rounding; see
\pdfclink{0 0 1}{{\tt www.mpfr.org/}}{http://www.mpfr.org/}

}

\bye